\documentclass[a4paper,reqno,oneside]{amsart}
\usepackage{amsmath,geometry}
\usepackage{graphicx}
\usepackage{mathrsfs}
\usepackage{amssymb,amsmath, amsfonts, amsthm}
\usepackage{latexsym}
\usepackage{color}
\usepackage[dvips]{epsfig}
\usepackage{graphicx}

\allowdisplaybreaks[1]

\numberwithin{equation}{section}

\renewcommand{\(}{\left(}
\renewcommand{\)}{\right)}
\renewcommand{\[}{\left[}
\renewcommand{\]}{\right]}

\newtheorem{theorem}{Theorem}[section]
\newtheorem{proposition}[theorem]{Proposition}
\newtheorem{lemma}[theorem]{Lemma}
\newtheorem{remark}[theorem]{Remark}

\renewcommand{\a }{\alpha }
\renewcommand{\b }{\beta }
\renewcommand{\d}{\delta }

\renewcommand{\le}{\leqslant}
\renewcommand{\ge}{\geqslant}
\renewcommand{\a }{\alpha }

\renewcommand{\b }{\beta }

\renewcommand{\d }{\delta }

\newcommand{\g }{\gamma }

\newcommand{\G}{\Gamma}

\newcommand{\D}{\mathcal{D}}
\newcommand{\M}{\mathcal{M}}

\newcommand{\beq}{\begin{equation}}
\newcommand{\eeq}{\end{equation}}
\newcommand{\beqs}{\begin{equation*}}
\newcommand{\eeqs}{\end{equation*}}
\newcommand{\beqn}{\begin{eqnarray}}
\newcommand{\eeqn}{\end{eqnarray}}
\newcommand{\beqns}{\begin{eqnarray*}}
\newcommand{\eeqns}{\end{eqnarray*}}
\newcommand{\bdoc}{\begin{document}}
\newcommand{\edoc}{\end{document}}
\newcommand{\be}{\begin{enumerate}}
\newcommand{\ee}{\end{enumerate}}
\newcommand{\bdescr}{\begin{description}}
\newcommand{\edescr}{\end{description}}
\newcommand{\ba}{\begin{array}}
\newcommand{\ea}{\end{array}}
\newcommand{\intR}{\int_{\mathbb R^N}}
\newcommand{\R}{\mathbb R}
\newcommand{\e}{\epsilon}
\newcommand{\ve}{\epsilon}
 \renewcommand{\(}{\left(}
\renewcommand{\)}{\right)}
\renewcommand{\[}{\left[}
\renewcommand{\]}{\right]}

\newcommand{\pa}{\partial}

\begin{document}
\title[Concentration on minimal submanifolds for a Yamabe type problem]{Concentration on minimal submanifolds for a Yamabe type problem}

\author{Shengbing Deng}
\address[Shengbing Deng]{School of Mathematics and Statistics, Southwest University,
Chongqing 400715, People's Republic of China,
 and
Departamento de Matem\'atica, Pontificia Universidad Catolica de
Chile, Santiago, Chile}
\email{shbdeng65@gmail.com}
\author{Monica Musso}
\address[Monica Musso]{Departamento de Matem\'aticas, Pontificia Universidad Cat\'olica de Chile, Santiago, Chile}
\email{mmusso@mat.puc.cl}
\author{Angela Pistoia}
\address[Angela Pistoia] {Dipartimento SBAI, Universt\`{a} di Roma ``La Sapienza", via Antonio Scarpa 16, 00161 Roma, Italy}
\email{pistoia@dmmm.uniroma1.it}

\begin{abstract} We construct solutions to a Yamabe type problem on a Riemannian manifold $M$ without boundary and
of dimension greater than $2$, with nonlinearity close to higher critical Sobolev exponents. These solutions concentrate their mass around a non degenerate minimal submanifold of $M$, provided a certain geometric condition involving the sectional curvatures is satisfied. A connection  with the solution of a class of P.D.E.'s   on the submanifold with a singular term of attractive or repulsive type is established.

  \end{abstract}

\date{\today}

\maketitle

\footnotetext{The first author was supported by Fondecyt grant 1130360 and Fondo Basal CMM. The second and the third authors have been partially supported by the Gruppo Nazionale per l'Analisi Matematica, la Probabilit\'a  e le loro Applicazioni (GNAMPA) of the Istituto Nazionale di Alta Matematica (INdAM).
}

\section{Introduction and statement of main results}\label{intro}

Let $(M,g)$ be a compact Riemannian manifold of dimension $m \geq 3$ without boundary.
This paper deals with the semilinear elliptic problem
 \begin{equation}\label{p} -\Delta _gu+h u=u^{q-1},\ u>0,\ \hbox{in}\ (M,g),\end{equation}
where the potential
$h\in C^2(M)$ is   such that $-\Delta_g+h$  is coercive and the exponent $q>2$.\\

Existence of non-trivial solutions to problem \eqref{p} is strictly related to the position of $q$ with respect to the critical Sobolev exponent
$2^*_m:={2m\over m-2}$.
   Indeed, in the subcritical case, i.e. $q<2^*_m$,
the Sobolev embedding $H^1_g(M)\hookrightarrow L^q_g(M)$ is compact for any $q\in(2,2^*_m)$
 and so
\begin{equation}\label{inf}\inf\limits_{u\in H^1_g(M)\atop u\not=0}{\int\limits_M\(|\nabla_g u|^2+h u^2\)d\sigma_g\over \(\int\limits_M
| u|^qd\sigma_g\)^{2/q}}\end{equation}
is achieved
 and problem \eqref{p} has a non-trivial solution.\\

 The critical case, i.e. $q=2^*_m$, has important links with the well known  Yamabe problem \cite{y}, namely
 find a metric $\widetilde g$ in the conformal class $[g]=\{\phi g:\, \phi \in C^\infty(M),\, \phi>0 \}$
 with constant  scalar curvature $\kappa.$ This is equivalent to set
$\widetilde g= u^{4/\(m-2\)} g$ and to find a solution $u$ to the Yamabe problem
\begin{equation}\label{yam}
-\Delta _gu+{m-2\over 4(m-1)}S_g  u=   \kappa u^{m+2\over m-2},\ u>0,\ \hbox{in}\ (M,g),\end{equation}
where $S_g$ is the scalar curvature of $(M,g). $
  The Yamabe problem on the round sphere $(\mathbb S^{m},g_0)$, equipped with the standard metric $g_0$,  plays a crucial role in solving problem \eqref{yam}.
 Since the scalar curvature of the round sphere is  $m(m-1)$ equation \eqref{yam} reduces to
 $$-\Delta _{g_0}u+ {m(m-2)\over 4} u=u^{m+2\over m-2},\ u>0,\ \hbox{in}\ (\mathbb S^m,g_0), $$
 which is equivalent   (via the stereographic projection) to the problem in the Euclidean space
 \begin{equation}\label{rm}-\Delta w=w^{m+2\over m-2},\ w>0,\ \hbox{in}\  \mathbb R^m.\end{equation}
Problem \eqref{rm} has infinitely many solutions (see \cite{cgs}),
\begin{equation}\label{bbs}  w_{\delta,y}(x):= \alpha_m\({\delta \over  \delta^2+|x-y|^2}\)^{m-2\over2},\  x,y\in\mathbb R^m,\ \delta>0,\end{equation}
where $\alpha_m:=\(m(m-2)\)^{m-2\over4}.$
In the general case, Aubin in \cite{a} proved that
if
$$ \mu_g(M,h):=\inf\limits_{u\in H^1_g(M)\atop u\not=0}\displaystyle{\int\limits_M\(|\nabla_g u|^2+h u^2\)d\sigma_g\over \(\int\limits_M | u|^{2m\over m-2}d\sigma_g\)^{m-2\over m}}$$
is such that
\begin{equation}\label{aubin}\mu_g(M,h)<\mu_{g_0}\(\mathbb S^m,{m(m-2)\over 4}\)\end{equation}
then $\mu_g(M,h)$ is achieved and so the corresponding critical problem \eqref{p} with $q={2^*_m}$ has a non-trivial solution.
The validity of \eqref{aubin} turns out to be strictly related to the position of the potential $h$ with respect to the geometric potential
\begin{equation}\label{geopot}\omega(\xi):={m-2\over 4(m-1)}S_g(\xi),\  \xi\in M.\end{equation}
Indeed, if
\begin{equation}\label{sotto}h(\xi)<\omega(\xi)\ \hbox{for any}\ \xi\in M,\end{equation}
it is not difficult to check that condition \eqref{aubin} holds.
In the case of the Yamabe problem, i.e. $h\equiv \omega$, condition \eqref{aubin} is also true, but the proof is a delicate issue.
It was proved by Trudinger \cite{t} when  $\mu_g(M,\omega) \le 0$, by
Aubin \cite{a} when $\mu_g(M,\omega) > 0$ and $(M,g)$ is not locally conformally flat and $m\ge6$ and by Schoen \cite{schoen}
when   $\mu_g (M,\omega)> 0$ and either $(M,g)$ is locally conformally flat or $3\le m\le 5$.\\

We can summarize known results just saying that problem \eqref{p} has a non-trivial solution if either $q<2^*_m$ and no extra assumptions on $h$,
or $q=2^*_m$ and $h$ has to satisfy $h\le \omega$ on $M.$
Therefore, it is natural to ask what happens when
$$ h(\xi)>\omega(\xi)\ \hbox{for some}\ \xi\in M\ \hbox{or}\
\hbox{$p$ is supercritical, i.e.}\ p>2^*_m.$$\\

A first partial answer was  given by Micheletti, Pistoia and  V\'etois \cite{mpv}
in a perturbative setting. They consider the almost critical problem
\begin{equation}\label{ac} -\Delta _gu+h u=u^{q_\epsilon-1},\ u>0,\ \hbox{in}\ (M,g),\ \hbox{with}\  q_\epsilon:=2^*_m\pm\epsilon\end{equation}
where $\epsilon$ is a positive small parameter.
If $q_\epsilon =2^*_m-\epsilon$ problem \eqref{ac} is said to be slightly subcritical,
   while if   $q_\epsilon =2^*_m+\epsilon$ it is said to be slightly supercritical.
They proved the following result.

\begin{theorem}\label{mpv-teo} [Theorem 1.1, Theorem 1.2, \cite{mpv}]
Let $m\ge6$ and $\xi_0\in M$ be a non degenerate critical point of $h-\omega$ (see \eqref{geopot}).
\begin{itemize}
\item[(i)] If $h(\xi_0)>\omega(\xi_0)$ then  there exists $\epsilon_0>0$ such that for any $\epsilon\in(0,\epsilon_0)$ the slightly subcritical problem \eqref{ac} with $q_\epsilon={2^*_m}-\epsilon$ has a solution   $u_\epsilon$ such that
$$|\nabla u_\epsilon|^2\rightharpoonup c_m\delta_{\xi_0}\ \hbox{as}\ \epsilon\to0.$$
\item[(ii)] If $h(\xi_0)<\omega(\xi_0)$ then there exists $\epsilon_0>0$ such that for any $\epsilon\in(0,\epsilon_0)$ the  slightly supercritical problem \eqref{ac} with $q_\epsilon={2^*_m}+\epsilon$ has a solution   $u_\epsilon$ such that
$$|\nabla u_\epsilon|^2\rightharpoonup c_m\delta_{\xi_0}\ \hbox{as}\ \epsilon\to0.$$
\end{itemize}
Here $\delta_{\xi_0}$
stands for the Dirac measure supported on $\xi_0$ and $c_m$ is an explicit positive
constant depending only on $m.$\end{theorem}
The profile of $u_\epsilon$  close to the concentration point $\xi_0$ is given by (see \eqref{bbs})
$$u_\epsilon(x)\approx \alpha_m \({\delta_\epsilon\over \delta_\epsilon^2+|x-\xi_0|^2}\)^{m-2\over2}$$
where the concentration parameter $\delta_\epsilon\sim d \sqrt\epsilon$, as $\epsilon \to 0$,
and the positive number $d$ solves
\begin{equation}\label{com1}
a_m\underbrace{\[h(\xi_0)-\omega(\xi_0)\]}_{>0}d -{b_m\over d }=0 
\quad \hbox{in the slightly sub-critical case}\end{equation}
 or
\begin{equation}\label{com2}a_m\underbrace{\[h(\xi_0)-\omega S_g(\xi_0)\]}_{<0}d +{b_m\over d }=0\quad \hbox{in the slightly super-critical case.}\end{equation}
 Here $a_m$ and $b_m$ are positive constants which only depend on $m.$\\

This result suggests to explore what happens when the exponent $q$ is close to higher critical exponents.
More precisely, for any integer $0\le k\le m-3$ we consider the $(k+1)$-st critical exponent $2^*_{m,k}={2(m-k)\over m-k-2} =2^*_{m-k,0}$, namely    the critical exponent for the Sobolev embedding  ${\mathrm H}^ 1_h(N)\hookrightarrow {\mathrm L}^{q}_h(N)$
  where $(N,h)$ is a $(m-k)-$dimensional Riemannian manifold. In particular, $ 2^*_{m,0}={2 m\over m- 2}$    is the usual Sobolev critical exponent. We know  by Theorem \ref{mpv-teo} that   problem \eqref{p} when the exponent $q$ approaches the first critical Sobolev exponent $2^*_{m,0}$ has  solutions   which blow-up at a single point. A set consisting of a single point is a $0-$dimensional submanifold of $M$. We ask:
{\it if $q$ approaches the $(k+1)$-st critical exponent $2^*_{m,k}$, do positive solutions blowing-up at $k-$dimensional submanifolds of $M$ exist?}
\\
Recently, a partial answer has been given by D\'avila, Pistoia and Vaira in \cite{dpv} when $k=1$ and by  Ghimenti, Micheletti and Pistoia in \cite{gmp} in a symmetric setting. Here we deal with the general case. Let us consider the   almost $k$-critical problem
\begin{equation}\label{pb} -\Delta _gu+h u=u^{q_\epsilon-1},\ u>0,\ \hbox{in}\ (M,g),\ \hbox{with}\ q_\epsilon:=2^*_{m,k}\pm\epsilon\end{equation}
where  $\epsilon$ is a positive small parameter.
If $q=2^*_{m,k}-\epsilon$ problem \eqref{pb} is said to be slightly $k$-th   subcritical,
   while if   $q=2^*_{m,k}+\epsilon$ it is said to be slightly  $k$-th supercritical.
\\

To state our result we need to introduce some geometric background. Let $K\subset M$ be a $k-$dimensional submanifold.
Set $N:=m-k.$
Let us introduce Fermi coordinates in $\M$ near the submanifold $K.$\\
Let $((E_a)_{a=1,\cdots, k},(E_i)_{i=1,\cdots, N })$ be  a local oriented and orthonormal frame field
along $K$. At
points $\xi$ of $K$, $T_\xi \M$ splits  as $T_\xi K \oplus N_\xi K$,
where $T _\xi K$ is the tangent bundle to $K$ with orthonormal basis
$(E_a)_a$ and $N_\xi K$ is the normal bundle, which is spanned
by the orthonormal basis  $(E_j)_j$.
We assume that  the normal vectors
$(E_i)_i$, $i = 1, \dots, N $, are parallel transported along $K$,
namely
\begin{equation}\label{eq:parall}
     g\left(  \nabla_{E_a}E_j\,,E_i\right)=0  \ \hbox{ at }\ \xi,
    \ \hbox{for any}\  i,j = 1, \dots, N ,\ a = 1, \dots, k.
\end{equation}
Here $\nabla$ is the connection associated with the metric $g$.
We denote by $\Gamma_a^b(\cdot)$ the 1-forms defined
on the normal bundle of $K$ by
\begin{equation}\label{eq:Gab}
    \Gamma_{ai}^b := \Gamma_a^b(E_i)=  g(\nabla_{E_a}E_b,E_i).
\end{equation}
The {\it minimal condition on $K$ } translates precisely into
\begin{equation}
\label{minimality}
\sum\limits_{a=1}^k\Gamma_{ai}^a  = 0\ \hbox{for any}\  i=1, \ldots
N.
\end{equation}
In a neighborhood of $\xi$ in $K$, we consider normal geodesic
coordinates
\begin{equation}\label{defoff}
f(y) : = \exp^K_\xi (\sum\limits_{a=1}^ky_a\, E_a ), \qquad y := (y_{1}, \ldots, y_{k}),
\end{equation}
where $\exp^K$ is the exponential map on $K$.
In a neighborhood of $\xi$ in $\M$, we introduce
\begin{equation}\label{eq:fc}
\mathfrak F ( y,  x) =\exp _{f(y)}(
\sum\limits_{i=1}^{N}\,x_i\,E_i); \qquad \quad (y,x) =
\left(( y_a)_a,(x_i)_i\right),
\end{equation}
where $\exp_{f(y) }$ is the exponential map at $f(y) $ in $\M$.
It holds true that $f(y)=\mathfrak F(y,0)\in K.$ Let $\tilde g_{ab}$ be the coefficients of the  induced metric on $K$ and
let
 $R_{\a\b\g\d}$  be the components of the curvature tensor computed at the point $\xi$ of
$K.$
{\it Non-degeneracy of $K$ } translates into the fact that the linear system
\begin{equation}\label{jacobi}
-\Delta_K \Phi_\ell+\sum\limits_{m=1}^N\sum\limits_{a,b=1}^k\(\tilde g^{ab}R_{mab\ell}-\Gamma^b_{am}\Gamma^a_{b\ell}\)\Phi_m=0,\ \ell=1,\dots,N,
\end{equation}
has only the trivial solution $\Phi=\(\Phi_1,\dots,\Phi_N\)\equiv 0.$\\

\textcolor{blue}{}
The Levi-Civita connection $\nabla$ for $g$ induces a connection $\nabla^N$ on the normal bundle $N_\xi K.$ We denote by
$$\mathcal R^N:=\sum\limits_{a=1}^k\(R(E_a,\cdot) E_a\)^N, $$
the curvature operator for this connection.
The second fundamental form
$$B:T_\xi K\times T_\xi K\to N_\xi K,\ B(X,Y)=\(\nabla _X Y\)^N$$
defines a symmetric operator
$\mathcal B^N:=B^t\cdot B,$
in terms of the coefficients $\Gamma^b_a:=B(E_b,E_a),$
$$g(\mathcal B^NX,Y)=\sum\limits_{a,b=1}^k\Gamma^b_a(X)\Gamma^a_b(Y),\ X,Y\in T_\xi M.$$
We also use the Ricci tensor
$$Ric(X,Y)=\sum\limits_{j=1}^N g\(R(X,E_j)E_j,Y\)+\sum\limits_{a=1}^k g\(R(X,E_a)E_a,Y\),\ X,Y\in T_\xi M.$$
We introduce the quadratic form
$$\mathcal Q(X,Y):=\frac 13 Ric(X,Y)-\frac23 g\(\mathcal R^N X,Y\)+g(\mathcal B^NX,Y),\ X,Y\in N_\xi K.$$
Recall that $N=m-k$. We set
\begin{equation}\label{hgk}
\Omega(\xi):={3(m-k-2)\over4(m-k-1)}\sum\limits_{i=1}^{m-k}\mathcal Q(E_i,E_i),\ \xi\in K.
\end{equation}
The expression of $\Omega$ in Fermi coordinates is given  by
$$
\hat \Omega (y) =  {3(N-2)\over 4(N-1)}\[\sum\limits_{i,j=1}^{N}
\frac13 R_{jiji}(y)+ \sum\limits_{i=1}^{N}\sum\limits_{a,b=1}^k
 \(\tilde g^{ab}R_{iaib}(y)+ \Gamma_{ai}^b(y)
\Gamma_{bi}^a(y) \)\] .
$$
This function appears in our construction in \eqref{defgy}.
Observe  that if $k=0$ the function $\Omega$ is nothing but the geometric potential ${m-2\over 4(m-1)} S_g(\xi)$ introduced in
\eqref{geopot}. Surprisingly enough, the function defined in \eqref{hgk} is not new in the literature, and it appears in a completely different context in \cite{mmpac} (Section 4). In fact, Mahmoudi, Mazzeo and Pacard in  \cite{mmpac}  deal with the existence of a family of Constant Mean Curvature submanifolds condensating to a fixed submanifold of a given Riemannian manifold. The limit submanifold has to be a closed non degenerate minimal submanifold. It would be interesting to further investigate the relation between our construction and the results in \cite{mmpac}. \\

Existence of solutions to the almost $k$-th critical problem  \eqref{pb} is strictly related to the existence of solutions to the elliptic PDE's with a singularity of   {\it attractive} type\begin{equation}\label{choiceofmusub}
- \Delta_K d +a_N\underbrace{\left[h(\xi)-\Omega(\xi)\right]}_{>0\ \hbox{in}\ K}d -\frac{b_N}{d}=0,\ d>0\ \hbox{in} \ K
\ \hbox{in the slightly $k$-th sub-critical case}
\end{equation}
or  with a singularity of   {\it repulsive} type 
\begin{equation}\label{choiceofmusup}
- \Delta_Kd +a_N\underbrace{\left[h(\xi)-\Omega(\xi)\right]}_{<0\ \hbox{in}\ K}d +\frac{b_N}{d}=0,\ \mu>0\ \hbox{in} \ K \ \hbox{in the slightly $k$-th super-critical case.}
\end{equation}
  Here \begin{equation}\label{abn}
a_N=\frac{4(N-1)}{(N-2)(N+2)}\ \hbox{and}\ b_N =   {(N-2)^2(N-4) \over 2 (N+2) }
.\end{equation}
 This relation is new and unexpected. It seems to be the natural extension of conditions \eqref{com1} and \eqref{com2} to higher critical problems.
 More precisely,
our main result reads as follows.

  \begin{theorem}\label{main}
Let $K\subset M$ be a closed non-degenerate minimal $k-$dimensional submanifold. Assume $m-k\ge 7.$
\begin{itemize}
\item[(i)] Assume  \eqref{choiceofmusub} has a non-degenerate solution. Then  there exists a sequence $\epsilon=\epsilon_ n\to0$ such that  the slightly $k$-th sub-critical problem \eqref{pb} with $q_\epsilon={2^*_{m,k}}-\epsilon$ has a solution   $u_\epsilon$ such that
$$|\nabla u_\epsilon|^2\rightharpoonup c_{m,k}\delta_{K}\ \hbox{as}\ \epsilon\to0.$$
\item [(ii)] Assume \eqref{choiceofmusup} has a non-degenerate solution. Then  there exists a sequence $\epsilon=\epsilon_ n\to0$ such that  the slightly $k$-th super-critical problem \eqref{pb} with $q_\epsilon={2^*_{m,k}}+\epsilon$ has a solution   $u_\epsilon$ such that
$$|\nabla u_\epsilon|^2\rightharpoonup c_{m,k}\delta_{K}\ \hbox{as}\ \epsilon\to0.$$
\end{itemize}
Here $\delta_{K}$
stands for the Dirac measure supported on $K$ and $c_{m,k}$ is an explicit positive
constant depending only on $m$ and $k.$\end{theorem}
The profile of $u_\epsilon$  close to the submanifold $K$ is given in Fermi coordinate by (see \eqref{bbs})
$$u_\epsilon(y,x)\approx \alpha_{m-k} \({\delta_\epsilon\over \delta_\epsilon^2+|x|^2}\)^{m-k-2\over2}$$
where the concentration parameter $\delta_\epsilon=\delta_\epsilon(y)$ satisfies $\delta_\epsilon\sim d \sqrt\epsilon$
and the positive function $d=d(y),$ defined on  $K,$ solves either the attractive singular PDE \eqref{choiceofmusub}
in the slightly sub-critical case or the repulsive singular PDE \eqref{choiceofmusup}
in the slightly super-critical case. \\

It is important to point out that if $\min\limits_{\xi\in K}\[h(\xi)-\Omega(\xi)\]>0$ then problem \eqref{choiceofmusub} has a non-degenerate solution as proved in Theorem \ref{th1}.
On the other hand, existence of solutions to problem \eqref{choiceofmusup} is a difficult issue, unless we deal with constant function $h(\xi)-\Omega(\xi)$ (see Remark \ref{rm4}). Indeed, as far as we know, there is only one result in the literature, which was proved by del Pino, Man\'asevich and Montero in \cite{demamo} in the case $k=1$, when
$\max\limits_{\xi\in K}\[h(\xi)-\Omega(\xi)\]<0$ (see Theorem \ref{geode}).
We would like to stress the fact that existence of solutions to problem  \eqref{choiceofmusup} is an interesting open question by itself, which as a by product allows to find solutions to the supercritical problem \eqref{pb}.\\

Another remark is that the result we find suggests that the natural extension to higher critical exponent of the classical Yamabe equation
is
\begin{equation}\label{hyp}
-\Delta_g u+\Omega(\xi) u=u^{m-k+2\over m-k-2},\ u>0\ \hbox{in}\ (M,g).
\end{equation}
where $\Omega$ is the function defined in \eqref{hgk}. If $k=0$ Problem \eqref{hyp} reduces to  the classical Yamabe equation since ${m-k+2\over m-k-2} = { m+2\over m-2}$ and $\Omega (\xi ) = {m-2\over 4(m-1)} S_g(\xi)$, as we already mentioned.
A natural open question is thus: {\it does problem \eqref{hyp} have a solution?}\\

Finally, we point out some interesting problems, whose solutions could help in understanding equation \eqref{hyp}.

\begin{itemize}
\item[(i)] Theorem \ref{main} holds true when $m\ge k+7.$ The question is:
 {\it does problem \eqref{pb} have any blowing-up solutions when $3\le m\le k+6$?} The case $k=0$ was completely studied by Druet \cite{do1,do2}.
\item[(ii)]  Theorem \ref{main} holds true when $h(\xi)\not= \Omega(\xi)$ for any $\xi\in K.$ The question is:
 {\it does problem \eqref{pb} have   any blowing-up solutions if $h(\xi)= \Omega(\xi)$ at some $\xi\in K$?} The case $k=0$ was extensively studied by Esposito, Pistoia and Vet\'ois \cite{epv} and by Esposito and  Pistoia in \cite{ep}.
\item[(iii)]  Theorem \ref{main} holds true when $q\to 2^*_{m,k}$.  The question is:
 {\it does problem \eqref{pb} have any solutions if  $q= 2^*_{m,k}$?} The case $k=0$ is nothing but the well known Yamabe problem.
\end{itemize}

\medskip
\noindent
In the last few years several investigations have been carried out around the possibility of constructing singular limit
solutions to non linear elliptic PDEs or problems in geometric analysis, depending on some parameters, whose mass or energy concentrate on sets of high dimension, like curves, surfaces, or higher dimensional sets. We refer the readers to \cite{mmpac,dmp,mmah,mm1,mm2,m,mussoyang,dkwei,demamu,dww,dengmamu} for instance, and the references therein.
First contributions on concentration at higher dimensional set for problems involving higher critical Sobolev exponents
are contained in the papers \cite{dmp,mussoyang,demamu}.

The general strategy used to prove all the above results is the so-called infinite dimensional version of the Liapunov-Schmidit reduction method. A main ingredient is to construct
an approximate solution   with arbitrary degree of accuracy in
powers of $\ve$, in a neighborhood of the submanifold manifold $K$. This approximation is, at main order,  a solution of some limit problem, which is independent of some of the variables. After this is done, one
builds the desired solution by linearizing the equation
around the approximation. The associated linear operator turns out
to be invertible with inverse controlled in a suitable norm by
certain large negative power of $\ve$, provided that $\ve$ remains
away from certain critical values where resonance occurs. The
interplay of the size of the error and that of the inverse of the
linearization then makes it possible a fixed point scheme.

\medskip
The rest of the paper is organized as follows.

We first discuss solvability and non-degeneracy of solutions to problems  \eqref{choiceofmusub} and  \eqref{choiceofmusup}. This is done in Section \ref{one}. In Section \ref{sec4} we introduce some scaled variables around the submanifold $K$ and we describe the Laplace Beltrami operator in these new variables.  Section \ref{aprsol} is
devoted to the construction of the approximate solution to our
problem using the local coordinates around the sub-manifold $K$
introduced before.  To perform this construction we need to invert a linear operator and to estimate the inverse. The  proof of this result is postponed to Section \ref{luigi}. In Section \ref{s:linear} we define globally the approximation and we write the solution to our problem as the sum
 of the global approximation plus a remaining term. Thus we express our original problem as a non linear problem in the remaining term and we prove our Theorem.
 To solve such problem, we need to understand the invertibility properties of another linear operator. To do so we start expanding a quadratic functional
 associated to the linear problem.
 This is done in Section \ref{linearas}.

\section{Some remarks on  PDE's with a singular term}\label{one}

First let us consider the attractive case, i.e. problem \eqref{choiceofmusub}. We can deal with a more general situation.
 \begin{theorem}\label{th1}
 Let $(M,g)$ be a smooth Riemannian compact manifold without boundary.
 Assume $\alpha,\beta\in C^0(M) $ and $\min\limits_M\alpha,\min\limits_M\beta>0.$
 Then there exists a non-degenerate solution to   
 \begin{equation}\label{p1}
\left\{\begin{aligned}
&-\Delta_g u+\alpha u-{\beta\over u}=0 \quad \hbox{in}\ M\\
&  u>0 \quad \hbox{in}\ M.\\
\end{aligned}\right.
\end{equation}
\end{theorem}
\begin{proof}

Let us prove that \eqref{p1} does have a solution.
Set $L(u):=-\Delta_g u+\alpha u.$
Let us rewrite problem \eqref{p1} in the following way
\begin{equation}\label{p11}L(u)=f(x,u)\ ,\ u>0\ \hbox{in}\ M,\end{equation}
 where $L(u)=-\Delta_g u+\alpha u$ and $f(x,u):={\beta\over u}.$ The linear operator $L$ is coercive.
 \\
First of all, we prove that problem \eqref{p1} has a lower solution $\underline u$ and an upper solution $\overline u,$  i.e.
$$L(\underline u)\le f(x,\underline u)\ \hbox{in}\ M\quad\hbox{and}\ L(\overline u)\ge f(x,\overline u)\ \hbox{in}\ M$$
such that
$$0<\underline u(x)\le \overline u(x)\ \hbox{for any}\ x\in M.$$
 It is enough to consider $\underline u$ and $\overline u$ as positive constant functions   and to observe that
 $L(c)-f(x,c)<0$ if $c$ is small enough and $L(C)-f(x,C)>0$ if $C$ is large enough.
 \\
 As a second step, we consider the modified problem
 \begin{equation}\label{p12}L(u)=\tilde f(x,u)\ ,\ u>0\ \hbox{in}\ M,\end{equation}
where
$$\tilde f(x,u):=\left\{\begin{aligned}
&{\beta(x)\over \underline u(x)}\ &\hbox{if}\ u(x)<\underline u(x)\\
&{\beta(x)\over  u(x)}\ &\hbox{if}\ \underline u(x)\le u(x)\le\overline u(x)\\
&{\beta(x)\over \overline u(x)}\ &\hbox{if}\ u(x)>\overline u(x)\\
\end{aligned}\right.
 $$
We point out that any solutions of the modified problem \eqref{p12} is a solution to the problem \eqref{p11}. Indeed, assume $u$ solves \eqref{p12}. We want to show that $\underline u(x)\le u(x)\le\overline u(x)$ for any $x\in M.$ Suppose, by contradiction that $\max\limits_M\(\underline u-u\)>0.$ Then there exists a point $x_0\in M$ such that $\(\underline u-u\)(x_0):=\max\limits_M\(\underline u-u\)>0$ and an open   set $\Gamma\subset M$ such that $x_0\in \Gamma $ and $\(\underline u-u\)(x)\ge 0$ for any $x\in\Gamma.$  Moreover, the function $ \underline u-u $ solves
$$L\(\underline u-u\)\le 0\ \hbox{in}\ \Gamma.$$ Since it achieves a maximum at the point $x_0$ which is in $\Gamma,$ by remark \ref{rm2} we immediately get a contradiction.
\\
As a final step, we prove that problem \eqref{p12} has a solution.
We remark that $u$ solves problem \eqref{p12} if $u$ is a fixed point of the operator $K(u):= T\(\tilde f(x,u)\),$ $u\in C^0(M)$ where $T$ is defined in Remark \ref{rm3}.
By Remark \ref{rm3} and by the definition of $\tilde f(x,u)$ we deduce that $K:C^0(M)\to C^0(M)$ is a compact operator and moreover that there exists $R>0$ such that
$\|K(u)\|_{C^0(M)}<R$ for any $u\in C^0(M).$ Hence for any $t\in[0,1]$, using the homotopy invariance of the Leray-Schauder degree, we get
$$\textrm{deg}\(I-K,B(0,R)\)=\textrm{deg}\(I-tK,B(0,R)\)=\textrm{deg}\(I,B(0,R)\)=1$$
and so problem \eqref{p12} has a solution. In order to prove that it is non degenerate, we point out that
the linearized equation
$$-\Delta_g v+\alpha v+{\beta\over u^2}v=0 \ , \ \hbox{in}\ M.$$
has only the trivial solution, since $\alpha$ and $\beta$ are strictly positive functions on $M.$\\

That concludes the proof.
\end{proof}

\begin{remark}\label{rm1}
 $L$ satisfies the maximum principle, namely
if $u\in H^1(M)$ is such that $Lu\le 0$ in $M$ then $u\le 0$ in $M.$
\end{remark}

\begin{remark}\label{rm2}
Assume that for some open  set $\Gamma\subset M$ the function
 $u\in H^1(M)$ solves $Lu\le 0$ in $\Gamma.$   Then if $u$ achieves its maximum  at a point $x\in \Gamma $ then $u(x)\le 0.$
\end{remark}

\begin{remark}\label{rm3}
For any $h\in C^0(M)$ there exists a unique $w\in C^2(M)$ such that
$Lu=h$ in $M.$ The linear map $T:C^0(M)\to C^0(M)$ defined by $Th=w$ is continuous and  compact.\end{remark}
\begin{proof}
It is enough to remark that the linear map $T:C^0(M)\to C^2(M)$ is continuos, because by standard elliptic regularity theory there exists a constant $c$ which only depends on $M$ and $g$ such that
$$\|Th\|_{C^2(M)}\le c\|h\|_{C^0(M)}\ \hbox{for any}\ h\in C^0(M).$$
Moreover, the embedding $C^2(M) \hookrightarrow C^0(M)$ is compact because of the Ascoli-Arzel\'a Theorem.
\end{proof}

As far as it concerns the repulsive case, i.e. equation \eqref{choiceofmusup}, we quote the results obtained  by del Pino, Man\'asevich and Montero in \cite{demamo} in the case $k=1$ and by  D\'avila, Pistoia and Vaira in \cite{dpv}.
 Let us consider the more general problem
 \begin{equation}\label{p1bis}
\left\{\begin{aligned}
&-\Delta_g u+\alpha u+{\beta\over u}=0 \quad \hbox{in}\ M\\
&  u>0 \quad \hbox{in}\ M,\\
\end{aligned}\right.
\end{equation}
where $\alpha,\beta\in C^0(M) ,$ with $\min\limits_M \beta>0.$

\begin{theorem}\label{geode}
Let $M$ be a one-dimensional  manifold whose length is $\ell.$ Assume
\begin{equation}\label{mamamo}-\({(\kappa+1)\pi\over2\ell}\)^2<\max\limits_{\xi\in M}\alpha(\xi)<-\({\kappa\pi\over2\ell}\)^2<0,
\end{equation}
for some integer $\kappa\ge1$.
Then problem  \eqref{p1bis} has a solution (see \cite{demamo}).\\
Moreover, it is non-degenerate for most functions $\alpha$'s (see \cite{dpv}).
\end{theorem}
In the general case, we can only make a few remarks.

 \begin{remark}\label{rm4}
 \begin{itemize}
 \item[(i)] If \eqref{p1bis} has a solution, then
$\min\limits_{\xi\in M}\alpha(\xi)<0.$
\item[(ii)] Let $\alpha=a$ and $\beta=b$ be constants. If $a<0$,
 then problem \eqref{p1bis}  has a constant solution,
which is non-degenerate if in addition $-{\lambda_{\kappa+1}}<2a <-{\lambda_\kappa}<0$ holds for some $\kappa.$
Here $\(\lambda_\kappa\)_{\kappa\ge1}$ denotes the sequence of eigenvalues of $-\Delta_Mu=\lambda_\kappa u$ on $M.$
\end{itemize}
\end{remark}
\begin{proof} To prove (i) it is enough to integrate equation \eqref{p1bis}  on $M$, so we get
$$ \int\limits_M \alpha (\xi)u(\xi) d\xi +\int\limits_M{\beta(\xi)\over u(\xi)}d\xi=0$$
which implies that  $\alpha$ has to be negative somewhere in $M.$
The proof of (ii) follows by straightforward computations.

\end{proof}

It would be really interesting to find conditions on $\alpha$ and $\beta$ which ensure the existence of a solution to problem \eqref{p1bis} in a more general setting.

 \section{Laplace-Beltami operator in scaled variables}\label{sec4}

In this section  we describe the
Laplace-Beltrami operator in some scaled variables, by means of the Fermi coordinates introduced in \eqref{eq:fc}.\\

Let $\mu_\ve$ be a positive smooth function $\mu_\e
= \mu_\e (y) $  defined on $K$ which we assume to be uniformly bounded, as $\e \to 0$, along $K$. Let also $\Phi_\e$ be a smooth normal section (in
$M$) $\Phi_\e\,:K\longrightarrow\,NK$ defined by $
  \Phi_\e (y)= \Phi_\e^j(y)E_j $, and we assume that $ \Phi_\e^j(y)$, $  j=1,\cdots,N$, are functions uniformly bounded, as $\e \to 0$, in $K$.
Having introduced the above function, we define the following change of variables
\begin{equation}
\label{ba}
 u(\mathfrak F (y, x))=(1+\alpha_\e)(\sqrt{\e}\,\mu_\e (y))^{-\frac{N-2}{2}} \, v\left(\frac{y}{\sqrt{\e} }, \frac{ x-\e \Phi_\e (y)}{\sqrt{\e}\,\mu_\e (y)} \right),
\end{equation}
where $\mathfrak F (y, x)$ is the change of variables defined in  (\ref{eq:fc}).
and
\begin{equation}
\label{b0}
  v=v(z,\xi), \quad z= {y \over \sqrt{\e}}, \quad \xi=\frac{ x-\e \Phi_\e}{\sqrt{\e}\,\mu_\e}.
\end{equation}
In (\ref{ba})  $\alpha_\e$ is a number  defined so that $(1+\alpha_\e)^{p\pm\e-1}\e^{\mp\frac{N-2}{4}\e }=1$, that is,
\begin{equation}\label{alphaepsilon}
\alpha_\e=\e^{\pm\frac{(N-2)^2}{16\pm4(N-2)\e}\e}-1.
\end{equation}
To emphasize the dependence of the above change of variables on
$\mu_\e$ and $\Phi_\e$, we will use the notation
\begin{equation} \label{defTT}
u={\mathcal T}_{\mu_\e , \Phi_\e} (v)  \quad \Longleftrightarrow
\quad u\ \ {\mbox {and}} \quad v \quad {\mbox {satisfy (\ref{ba})}}.
\end{equation}

Recall that   the original variables
$(y,x)\in \R^{k+N}$ are {\it  local} coordinates along $K$. Thus  we let the
variables $(z, \xi )$ vary in the set $\D$ defined by
\begin{equation}
\label{defD} \D = \left\{ (z,   \xi  ) \, : \, \sqrt{\e} z \in K,
\quad | \xi | <{\eta \over \sqrt{\e}} \right\}
\end{equation}
for some small and fixed positive number $\eta$ that will be fixed in the sequel. We will also use
the notation $ \D = K_\e \times \hat \D$, where $K_\e = {K\over
\sqrt{\e}}$ and
\begin{equation}
\label{hatD}
\hat \D = \left\{  \xi  \, : \, | \xi | < {\eta \over \sqrt{\e}} \right\}.
\end{equation}
We note that $\partial\hat{ \D} = \left\{  \xi  \in\hat\D \, : \, | \xi | = {\eta \over \sqrt{\e} } \right\}$.

We are interested in computing  the Laplace Beltrami
operator in the new variables $(z, \xi )$ in terms of the parameter
$\e$, of the function $\mu_\e (y)$ and of the normal section
$\Phi_\e$.

\medskip
We have the validity of the following

\begin{lemma} \label{scaledlaplacian}
Given  the change of variables defined  in \eqref{ba}, the following
expansion for the Laplace Beltrami operator holds true
\begin{equation}
\label{lap1}
 (1+\alpha_\e)^{-1} \e^{{N+2 \over 4}}\mu_\e^{{N+2 \over 2}} \Delta_g u=  {\mathcal A}_{\mu_\e , \Phi_\e } (v) := {\mu_\e^2}  \Delta_{K_\e} v +
\Delta_{\xi} v + \sum_{\ell=0}^2 {\mathcal A}_\ell v + B(v).
\end{equation}
Above, the expression ${\mathcal A}_k$ denotes  specific
differential operators, respectively defined as follows
\begin{equation} \label{D0}
\begin{array}{rlllll}
{\mathcal A}_0 v
&= &   -
\e\,\mu_\e \,\Delta_K
 \mu_\e \,\left( \gamma
v + D_\xi v \, [\xi] \right) \\[3mm]
 & + &   \e \,| \nabla_K\mu_\e|^2 \left[ D_{\xi\xi} v \, [
\xi]^2 +
2 (1+\gamma ) D_\xi v [\xi] + \gamma (1+ \gamma ) v \right]   \\[3mm]
& - &  \e^{3\over 2} \mu_\e D_{ \xi}\, v \,[\Delta_K \Phi_\e ] + \e^{3\over 2} \nabla_K \mu_\e \,\cdot\,\left\{ 2D_{\xi\,\xi} v[\xi] + N
D_{\xi} v
\right\}\, [\nabla_K  \Phi_\e ] \\[3mm]
&+ & \e^2 D_{\xi\,\xi} v \,[\nabla_K \Phi_\e]^2\\[3mm]
& - & 2\, \e^2 \mu_\e \, g^{ab}\,\left[  D_\xi (\partial_{\bar a} v ) [\partial_b
\mu_\e \xi] +\e^{1\over 2} D_{\xi} (\partial_{\bar a} v )[ \partial_b
\Phi_\e ] +
 \gamma  \partial_a\mu_\e\, \partial_{\bar b} v \right],
\end{array}
\end{equation}
where we have set $\gamma=\frac{N-2}{2}$,
\smallskip
\begin{equation} \label{D1}
 {\mathcal A}_1 \, v  =  \,-{\e\over 3}\,\, \sum\limits_{i,j} \bigg[\sum\limits_{m,l} R_{mijl}
(\mu_\e \xi_m +\sqrt{\e} \Phi_\e^m )
(\mu_\e \xi_l + \sqrt{\e} \Phi_\e^l ) \bigg]  \partial^2_{ij} v,
\end{equation}
and
\smallskip
\begin{equation}
\label{D4} {\mathcal A}_2 v =\e \mu_\e \sum\limits_j \bigg[ \sum_s
\frac23 R_{mssj}+\sum\limits_{m,a,b}
 \big(  {\tilde g}^{ab}\, R_{mabj}- \G_{am}^b
\G_{bj}^a \big)\bigg](\mu_\e\xi_m+\sqrt{\e} \Phi_\e^m)
 \partial_j v.
\end{equation}
Finally,
the operator $B(v)$ can be described as follows: $B(v) = \e^2 \hat B (v)$
\begin{eqnarray*}
\hat {\mathcal{B}}(v)&=&O \left(  |\mu_\ve \xi + \sqrt{\e} \Phi |^3   \right) \partial^2_{ij} v\\
&+&
 O(| \mu_\e  + \sqrt{\e} \Phi_\e + \partial_{z_a} ( \mu_\e + \sqrt{\e} \Phi_\e ) |^2 ) O ( v + \xi_i  \partial_{\xi_i }  v ).
\end{eqnarray*}

\smallskip
We recall that the symbols $\partial_{a}$, $\partial_{\over a}$ and $\partial_i$
denote the derivatives with respect to $\partial_{y_a}$, $\partial_{z_a}$ and
$\partial_{\xi_i}$ respectively.

\end{lemma}

\medskip
\begin{proof}
Recall that the Laplace-Beltrami
operator is defined by
$$ \Delta_{
g}=\frac{1}{\sqrt{\det g}}\,\partial_A(\,\sqrt{\det g}\,(
g)^{AB}\,\partial_B\,)\,,
$$
where indices $A$ and $B$ run between
$1$ and $n=N+k$. In other words
\begin{equation}
\label{zero0}\Delta_{
g}
=({g})^{AB}\,\partial^2_{AB}+\partial_A\,({
g})^{AB}\,\partial_B+\partial_A(\,\log{\sqrt{\det
g}}\,)\,({g})^{AB}\,\partial_B
\end{equation}
If now $u$ and $v$ are defined as in (\ref{ba}), we have
$$
 (1+\alpha_\e)^{-1} \e^{N \over 4} \mu_\e^{N \over 2} \partial_{x_j} u = \partial_{\xi_j } v,
\quad   (1+\alpha_\e)^{-1} \e^{N+2 \over 4} \mu_\e^{N+2 \over 2} \partial^2_{x_j , x_i} u = \partial^2_{\xi_j , \xi_i } v
$$
and
\begin{eqnarray*}
 (1+\alpha_\e)^{-1}  \e^{N+2 \over 4} \mu_\e^{N+2 \over 2}  \partial^2_{y_a , y_b} u &=& \mu_\e^2 \partial^2_{z_a z_b} v + \mu_\e^{N-2 \over 2} \partial^2_{z_a z_b} (\mu_\e^{-{N-2 \over 2}} ) v  \\
&+& 2  \mu_\e^{N-2 \over 2} \partial^2_{z_a } (\mu_\e^{-{N-2 \over 2}} ) [ \partial_{z_b} v + \partial_{z_b} (\mu_\e^{-1} ) \nabla v \cdot \xi - \sqrt{\e} \nabla v \cdot \partial_{z_b} \Phi_\e ]  \\
&+& \partial_{z_a} [ \partial_{z_b} (\mu_\e^{-1} ) \nabla v \cdot \xi - \sqrt{\e} \nabla v \cdot \partial_{z_b} \Phi_\e ]
\end{eqnarray*}
On the other hand, by our choice of coordinates (\ref{defoff}), on $K$ the metric $g$ splits in
the following way
\begin{equation}\label{eq:splitovg}
     g(q) =   g_{ab}(q)\,d y_a\otimes d y_b+
g_{ij}(q)\,dx_i\otimes dx_j, \qquad \quad q \in K.
\end{equation}
If we denote by $r$ the distance function from $K$, at any point $\mathfrak F (y,x)$ (see (\ref{eq:fc}),
 we have
$$
\begin{array}{rllll}  g_{ij}(y, x)&=\delta_{ij}+\frac{1}{3}\,R_{istj}\,x_s\,x_t\,
+\,{\mathcal O}(r^3);\\[3mm]
  g_{aj}(y,  x)&={\mathcal O}(r^2);\\[3mm]
  g_{ab}(y,
x)&={\tilde g}_{ab}-\, [  {\tilde g}_{ac}\,\Gamma_{bi}^c+{\tilde g}_{bc}\,\Gamma_{ai}^c ] \,x_i+\left[R_{sabl}+ {\tilde g}_{cd} \Gamma_{a s}^c\,
\Gamma_{d l}^b \right]x_s x_l+{\mathcal O}(r^3).
\end{array}
$$
Here  $a=1,...,k$,  $i,j=1,...,N$. See \cite{demamu}.
Let now $g^\e$ be the scaled metric on $\M_\e=\e^{-1/2}\M,$ whose coefficients are defined by
$$g_{\alpha , \beta}^\e(z,x)=g_{\alpha,\beta}(\sqrt\e z,\sqrt \e x).$$
For the metric $g^\e$ in the above coordinates $(z,x )$
we have the expansions
\begin{eqnarray*}
&&g^\e_{ij}=\delta_{ij} +
\frac{ \e }{3} \,R_{istj}\,x_s\,x_t
\,+\,{\mathcal O}(\e^{3\over 2} (|x|^3),  \quad 1\le i,j\le N;\\[3mm]
 &&g_{aj}^\e={\mathcal O}(\e |x|^2)
 \quad 1\le a\le k , \, 1\le j\le  N;\\[3mm]
 &&g_{ab}^\e={\tilde g}^\e_{ab}-\sqrt{\e} \bigg\{{\tilde g}^\e_{ac}\,\Gamma_{bi}^c+{\tilde g}^\e_{bc}\,\Gamma_{ai}^c\bigg\}\,x_i
+
\e\,\bigg[R_{sabl}+{\tilde g}^\e_{cd}\Gamma_{as}^c\,
\Gamma_{dl}^b \bigg]x_s x_l + {\mathcal O}(\e^{3 \over 2} |x|^3), \\
&&\quad 1\le a,b\le k.
\end{eqnarray*}
Thus we first conclude that
\begin{eqnarray}\label{one1}
 (1+\alpha_\e)^{-1}  \e^{N+2 \over 4} \mu_\e^{N+2 \over 2} ({g})^{AB}\,\partial^2_{AB} &=& \mu_\e^2 g^{ab} \partial^2_{z_a z_b} v + \partial^2_{jj} v  \nonumber \\
&-& {\e \over 3} R_{mijl} (\mu_\e \xi_l + \sqrt{\e} \Phi_\e^l )  (\mu_\e \xi_m + \sqrt{\e} \Phi_\e^m ) \, \partial_{\xi_i \xi_j} v \nonumber \\
&+& {\mathcal B}_0 (v) +  {\mathcal B}_1 (v)
\end{eqnarray}
where
\begin{eqnarray*}
{\mathcal B}_0 (v) &=& -
\e\,\mu_\e \, g^{ab} \partial_{z_a z_b}^2
 \mu_\e \,\left( {N-2 \over 2}
v + D_\xi v \, [\xi] \right) \\[3mm]
 & + &   \e \, g^{ab} \partial_{z_a} \mu_\e \partial_{z_b} \mu_\e  \left[ D_{\xi\xi} v \, [
\xi]^2 +
N D_\xi v [\xi] + {(N-2) N \over 4} v \right]   \\[3mm]
& - &  \e^{3\over 2} \mu_\e g^{ab} \partial^2_{z_a z_b} \Phi_\e^ j \partial_{\xi_j} v+ \e^{3\over 2} g^{ab} \partial_{z_a} \mu_\e \partial_{z_b } \Phi_\e^l \left\{ 2 \partial^2_{\xi_j \xi_l} v \xi_j + N
\partial_{\xi_l} v
\right\}\\[3mm]
&+ & \e^2 g^{ab} \partial_{z_a} \Phi_\e^j \partial_{z_b} \Phi_\e^l \partial^2_{\xi_j \xi_l} v \\[3mm]
& - & 2\, \e^2 \mu_\e \, g^{ab}\,\left[  D_\xi (\partial_{\bar a} v ) [\partial_b
\mu_\e \xi] +\e^{1\over 2} D_{\xi} (\partial_{\bar a} v )[ \partial_b
\Phi_\e ] +
 \gamma  \partial_a\mu_\e\, \partial_{\bar b} v \right],
\end{eqnarray*}
and
$$
{\mathcal B}_1 (v) =O(\e^2 |\mu_\e \xi + \sqrt{\e} \Phi_\e |^3 ) \partial^2_{\xi_i \xi_j}  v.
$$

\medskip
\noindent
Moreover
\begin{eqnarray} \label{two2}
 (1+\alpha_\e)^{-1}  \e^{N+2 \over 4} \mu_\e^{N+2 \over 2} \partial_A\,({
g})^{AB}\,\partial_B u &=& \mu_\e^2 \partial_a (g^{ab} ) \partial_{z_b} v \nonumber \\
&+& {\e \over 3} \, \mu_\e R_{liij} (\mu_\e \xi_l +\sqrt{\e} \Phi_\e^l ) \, \partial_{\xi_j} v
\nonumber \\
&+& {\mathcal B}_2 (v)
\end{eqnarray}
where
\begin{eqnarray*}
{\mathcal B}_2 (v) &=&  -
\e\,\mu_\e \, \partial_{z_a} (g^{ab}) \partial_{z_b}
 \mu_\e \,\left( {N-2 \over 2}
v + D_\xi v \, [\xi] \right) \\[3mm]
& - &  \e^{3\over 2} \mu_\e  \partial_{z_a }( g^{ab}) \partial_{ z_b} \Phi_\e^ j \partial_{\xi_j} v
\\[3mm]
&+&\e^2   O(| \mu_\e  + \sqrt{\e} \Phi_\e + \partial_{z_a} ( \mu_\e + \sqrt{\e} \Phi_\e )|^2  ) O ( v + \xi_i  \partial_{\xi_i }  v ).
\end{eqnarray*}
Finally,  we recall that we have the validity of the following expansions for the
square root of the determinant of $g^\e$ and the log of determinant of $g^\e$
\begin{eqnarray}\label{expdeterminante}
 \sqrt{\det g^\e}  &=&   \sqrt{\det  g^\e}   \times \\
& & \bigg \{1+ \frac{\e}6  R_{miil}  x_m x_l + \frac{\e}2
 \bigg(  {\tilde g}^{ab}\,R_{mabl}-\G_{am}^c
\G_{cl}^a  \bigg) x_m   x_l
 +    \e^{3 \over 2} \mathcal{O}(|x|^3) \biggl\} \nonumber
\end{eqnarray}
and
\begin{eqnarray*}
 \log\big(\det g^\e\big) &=&  \log\big(\det  g^\e\big)+ \frac{\e}3  R_{miil}\,x_m x_l \\
& + &
 \e\bigg( {\tilde g}^{ab}\, R_{mabl}-\G_{am}^{c}
\G_{cl}^{a}  \bigg)x_m x_l + \mathcal{O}(\e^{3\over 2}|x|^3).
\end{eqnarray*}
See for instance \cite{demamu}.

So we get
\begin{eqnarray} \label{three3}
 (1+\alpha_\e)^{-1}  \e^{N+2 \over 4} &&\mu_\e^{N+2 \over 2}
\partial_A(\,\log{\sqrt{\det
g}}\,)\,({g})^{AB}\,\partial_B u =  \partial_a(\,\log{\sqrt{\det
g}}\,)\,({g})^{ab}\,\partial_b v \nonumber \\
&+& \e \left( {R_{mssj} \over 3} + ({\tilde g}^{ab} R_{mabj} - \Gamma_{am}^c \Gamma_{cj}^a ) \right) (\mu_\e \xi_m + \sqrt{\e} \Phi_\e^m ) \partial_{\xi_j } v \nonumber \\
&+& {\mathcal B}_3 (v).
\end{eqnarray}
Here ${\mathcal B}_3 (v)$ is a function that can be described as follows
\begin{eqnarray*}
{\mathcal B}_3 (v) &=& \e^2   O(| \mu_\e  + \sqrt{\e} \Phi_\e + \partial_{z_a} ( \mu_\e + \sqrt{\e} \Phi_\e )  |^2 ) O ( v + \xi_i  \partial_{\xi_i }  v ).
\end{eqnarray*}
\medskip
\noindent
Colleting (\ref{one1}), (\ref{two2}) and (\ref{three3}) in (\ref{zero0}), we get the proof of the Lemma.

\end{proof}

\medskip

 \setcounter{equation}{0}
\section{Construction of an approximate solution}\label{aprsol}

Using the local coordinates along the submanifold $K$ introduced in
Section \ref{sec4}, after performing the change of variables in
\eqref{ba}, the original equation in $u$  reduces locally close to
$K_\e $ to the following equation in $v$
\begin{equation}
\label{adesso}
 -{\mathcal A}_{\mu_\e , \Phi_\e } v + \e\, \mu_\e^2 \, h v-   \mu_\e^{\mp\e\frac{N-2}{2}}v^{p\pm\e} =0 ,
\end{equation}
where ${\mathcal A}_{\mu_\e , \Phi_\e }$ is defined in \eqref{lap1}
and $p={N+2 \over N-2}$. Let us denote by $\Xi_\e$ the
operator given by
\begin{equation}
\label{Sep} \Xi_\e (v ) := -{\mathcal A}_{\mu_\e ,
\Phi_\e } v + \e\mu_\e^2 \, hv-  \mu_\e^{\mp\e\frac{N-2}{2}}v^{p\pm\e}.
\end{equation}

\medskip
\noindent
This section is devoted to build an approximate solution to Problem \eqref{adesso} locally around $K_\e$, in the set $\D=K_\e\times\hat\D$,
(see Section \ref{sec4}).

Let $r$ be an
integer. For a function $w$ defined in $\D=K_\e \times \hat\D$, we define
\begin{equation}\label{eqinftynu}
\|w\|_{\e,r}:=\sup_{(z,\xi)\in K_\e \times \hat\D}\left( \,(1+|\xi|^2)^{r \over 2}|w(z,\xi)|
\right).
\end{equation}
Let $\sigma \in (0,1)$. We define
\begin{equation}
\label{normsigma} \| w \|_{\e, r, \sigma} := \| w \|_{\e , r} + \sup_{(z,\xi)\in K_\e
\times \hat\D}\left( \,(1+|\xi|^2)^{r
+\sigma\over 2} [w]_{\sigma, B(\xi , 1)} \right)
\end{equation}
where we have denoted
\begin{equation}\label{eqnorms}
[w]_{\sigma, B(\xi , 1)} := \sup_{ \xi_1, \xi_2\in B(\xi , 1)}
\frac{|w(z,\xi_2)-w(z,\xi_1)|}{|\xi_1-\xi_2|^\sigma}
\end{equation}

The main result of this section is as follows.

\medskip
\begin{lemma} \label{Construction}
There exist $\e_0 >0$, $\eta >0$ in the definition of ${\mathcal D}$ in (\ref{defD}), and a constant $C>0$, such that, for any integer $I $ and for all $\e \in (0, \e_0)$ there exist  a smooth function $\mu_{I+1,\e} :K
\to \R$,  a smooth normal section $\Phi_{I+1,\e} : K \to NK$, of the form $\Phi_{I+1, \e} (y) = \Phi_{I+1, \e}^j (y) E_j$
\begin{equation}
\label{bf2}
\|   \mu_{I+1,\e}\|_{ \infty  } +\| \partial_a
\mu_{I+1,\e}\|_{\infty } +\|\partial^2_a   \mu_{I+1,\e}\|_{ \infty
} \leq C
\end{equation}
\begin{equation}
\label{bf3}
\|  \Phi_{I+1 ,\e} \|_{ \infty } +\| \partial_a
\Phi_{I+1 , \e} \|_{ \infty } +\|\partial^2_a
\Phi_{I+1 , \e} \|_{ \infty } \leq C,
\end{equation}
and a positive
function $v_{I+1, \e} :K_\e \times \hat\D \to \R$ such that
$$
-{\mathcal A}_{\mu_{I+1,\e} , \Phi_{I+1 , \e}}   (v_{I+1 , \e} ) +\e  \mu_{I+1,\e}^2 h v_{I+1
, \e} - \mu_{I+1,\e}^{\mp \frac{N-2}{2}\e} v_{I+1, \e}^{p\pm \e} = {\mathcal E}_{I+1 , \e} \quad {\mbox {in}}
\quad  \D
$$
with
\begin{equation}
\label{boh1}
\| v_{I+1 , \e} - v_{I , \e} \|_{\e , N-4 , \sigma } \leq C
\e^{I+\frac{1}{2}}
\end{equation}
and
\begin{equation} \label{bf4}
\| {\mathcal E}_{I+1 , \e} \|_{\e , N-2 , \sigma} \leq C \e^{I+\frac{1}{2}}.
\end{equation}
We refer to (\ref{lap1}) for the definition of ${\mathcal A}_{\mu_\e , \Phi_\e} $, to (\ref{defD}) for $ K_\e \times \hat\D$.
\end{lemma}

The proof of Lemma \ref{Construction} is based on an explicit construction of the functions $\mu_{\e , I+1}$,
$\Phi_{\e , I+1}$ and $v_{\e , I+1}$, via an iterative scheme, in the spirit developed in \cite{demamu}.
Fix an integer $I>1$, we will define the functions $\mu_{\e , I}$ and $\Phi_{\e, I}$ respectively of the form
\begin{equation}
\label{muep}
\mu_{I,\e}:=   \mu_0 + \e  \mu_1 + \e^2 \mu_2+ \ldots + \e^{I-1} \mu_{I-1},
\end{equation}
and
\begin{equation}
\label{phiep}
\Phi_{\e}: =  \Phi_{1,\e} + \Phi_{2,\e}+\ldots + \Phi_{I-1,\e}.
\end{equation}
to be solutions of  certain linear elliptic PDEs on the sub manifold $K$. The solvability of these equations is
related to the result contained in Section \ref{one}.
At each step $I$, we also define
\begin{equation}
\label{roma1}
v_{I} (z, \xi ) : =  w_0 (\xi ) + w_{1,\e} (z, \xi )+ w_{2,\e} (z, \xi ) + w_{3,\e} (z, \xi )+ \ldots + w_{I,\e} (z, \xi ),
\end{equation}
where each term $w_{j,\e}$ will also be solution of a linear problem, this time defined on $\D$. The function $w_0$ has been already defined as solution to \begin{equation}
\label{w0}
\Delta u + u^{N+2 \over N-2} =0\quad \mbox{in}\  \R^N,
\end{equation}
given explicitely by
\begin{equation}
\label{defw0} w_0 (\xi ) = \alpha_N (1+ |\xi|^2 )^{-{N-2 \over 2}} .
\end{equation}
We consider the domain $\D$ defined as (\ref{defD}) and for function $\phi$ defined on $\D$, an operator of the form
$$
L(\phi):=- \Delta_\xi \phi - p w_0^{p-1} \phi +\e\, a(\e z) \phi,
$$
where $a$ is a given smooth function  $a:K \to \R$ with $a(y) \geq
\lambda >0$ for all $y \in K$.

Let us introduce the functions
\begin{equation}
 \label{lezetas}
Z_j (\xi ) = {\partial w_0 \over \partial \xi_j} , \quad j=1, \ldots
, N\quad {\mbox {and}} \quad Z_0 (\xi ) = \xi \cdot \nabla w_0
(\xi )  + \frac {N-2}2 w_0 (\xi )
\end{equation}
that are known to be the only bounded solutions to the linearized equation around
$w_0$ of problem \eqref{w0}
$$
-\Delta \phi - p w_0^{p-1} \phi =0 \quad {\mbox {in}} \quad \R^{N}.
$$
See \cite{bianchiengel}.

Given a function $g:K \times \hat\D\to \R$ that depends smoothly on the variable $y \in K$,
we want to find a linear theory for the following linear problem
\begin{equation}\label{eq:eqwd}
 \left\{
    \begin{array}{ll}
    L(\phi)=h,  &\  \hbox{ in } \D\\
\phi = 0 & \hbox{ on } \partial \hat \D \\
    \int_{\hat\D} \phi (\e z, \xi ) Z_j (\xi ) \, d\xi = 0 &\  \forall z \in K_\e, \quad  j=0, \ldots N.
    \end{array}
  \right.
\end{equation}

We have the validity of the following result.

\begin{proposition}\label{linear}
Let $r$ be an integer such that $4<r < N$.

Let $a : K \to \R$ be a smooth function, such that $a (y) \geq \lambda
>0$ for all $y \in K$. Then there exist $\e_0 >0$,  $\eta>0$, that depends only on $\sup_{y\in K} |a(y)|$, in the definition of $\D$ in (\ref{defD}), and $C>0$ such that, for any $\e \in (0, \e_0 )$ and  for any
function  $h : K \times \hat\D
\to \R $  that depends smoothly on the variable $y \in
K$, such that $ \| h \|_{\e , r} $ is bounded, uniformly in $\e$,
and
$$
\int_{\hat\D}h(\e z,\xi)Z_j(\xi)d\xi=0\quad \mbox{for\ all}\ z\in K_\e,\ \ j=0,1,\ldots,N,
$$
then there exists a solution $\phi$ of problem (\ref{eq:eqwd}) such that
\begin{equation}
\label{est0a} \| D^2_\xi \phi \|_{\e , r , \sigma } + \| D_\xi \phi
\|_{\e , r -1 , \sigma } +\|\phi \|_{\e, r- 2 , \sigma}\le C
\|h\|_{\e,r ,\sigma}
\end{equation}
Furthermore, the function $\phi$ depends smoothly on the variable
$\sqrt{\e} z$, and the following estimates hold true: for any integer $l$
there exists a positive constant $C_l$ such that
\begin{equation}
\label{est1a} \| D^l_z \phi \|_{\e , r- 2 , \sigma }  \le C_l \left(
\sum_{k\leq l} \|D^k_z h\|_{\e,r , \sigma}\right).
\end{equation}
\end{proposition}

We postpone the proof of Proposition \ref{linear} to Section \ref{luigi}.
We devote the rest of the section to the Proof of Proposition \ref{Construction}.

\begin{proof}[Proof of Proposition \ref{Construction}]

Define
\begin{equation}
\label{roma2}
\hat {\mathcal A} =\Delta_{\R^N} v + \sum_{\ell=0}^2+ {\mathcal A}_\ell v + B(v)
\end{equation}
referring to  Lemma \ref{scaledlaplacian}, and
\begin{equation}
\label{esse1}
\Xi_\e (u) = - \hat {\mathcal A} u+ \e \mu_\e^2 \, h \,  u -   \mu_\e^{\mp\e\frac{N-2}{2}}\,u^{p\pm\e}
\end{equation}

\medskip
\noindent
\noindent  We start with $I=1$ and the construction of $w_{1,\e}$ and $\mu_0$ .\ \
A direct computation gives
\begin{eqnarray*}
\Xi_\e(v_1)%&=&- \hat {\mathcal A} v_{1} + \e \mu_0^2 \, h \,  v_{1} -   \mu_0^{\mp\e\frac{N-2}{2}}\,v_{1}^{p\pm\e}\\
& = &
- \hat {\mathcal A} ( w_0+w_{1,\e}) + \e \mu_0^2 \, h \,   w_0+ \e \mu_0^2 \, h \,w_{1,\e} -  \mu_0^{\mp\frac{N-2}{2}\e}\,( w_0+w_{1,\e}) ^{p+\e}\\
& = &
- \hat {\mathcal A} ( w_0+w_{1,\e}) + \e  \mu_0^2 \, h \,  w_0+ \e \mu_0^2 \, h \,w_{1,\e}   \\
&&-   \mu_0^{\mp\frac{N-2}{2}\e}\, w_0^{p\pm\e}- (p\pm\e) \mu_0^{\mp\frac{N-2}{2}\e}\,w_0^{p-1\pm\e}w_{1,\e}\\
&&\underbrace{-  \mu_0^{\mp\frac{N-2}{2}\e}\,\left[ (w_0+w_{1,\e}) ^{p+\e}-w_0^{p\pm\e}-(p\pm\e) w_0^{p-1\pm\e}w_{1,\e}\right]}\limits_{Q_\e (w_{1 })}\\
& = & - \Delta w_{1,\e}-pw_0^{p-1}w_{1,\e}-\sum_{\ell=0}^2 {\mathcal A}_\ell w_{1,\e} - B(w_{1,\e})-  \sum_{\ell=0}^2 {\mathcal A}_\ell w_0 - B(w_0) \\
&&- \underbrace{\left\{ \mu_0^{\mp\frac{N-2}{2}\e}\, w_0^{p\pm\e}-w_0^p\right\} }\limits_{I_1}+ \e \mu_0^2 \, h \,  w_0\\
&&- \underbrace{\left\{(p\pm\e)  \mu_0^{\mp\frac{N-2}{2}\e}\,w_0^{p-1\pm\e}w_{1,\e}-pw_0^{p-1}w_{1,\e}\right\}}\limits_{I_2}+ \e \mu_0^2 \, h \,w_{1,\e}+Q_\e (w_{1 })\\
& = &- \Delta w_{1,\e}-pw_0^{p-1}w_{1,\e}+ \e \mu_0^2 \, h \,w_{1,\e} + \e \mu_0^2 \, h \,  w_0-\underbrace{ \mu_{0}^{\pm\frac{(N-2)^2}{8}\e} \sum_{\ell=0}^2 {\mathcal A}_\ell w_0}\limits_{I_0}\\
&&- \underbrace{  \left\{ \mu_0^{\mp\frac{N-2}{2}\e}\,  w_0^{p\pm\e}-w_0^p\right\}}\limits_{I_1}- \underbrace{\left\{(p\pm\e)  \mu_0^{\mp\frac{N-2}{2}\e}\,w_0^{p-1\pm\e}w_{1,\e}-pw_0^{p-1}w_{1,\e}\right\}}\limits_{I_2}\\
&& - B( w_0) -\sum_{\ell=0}^2 {\mathcal A}_\ell w_{1,\e} - B(w_{1,\e})+Q_\e (w_{1 }).
\end{eqnarray*}
We next analyze each one of the above terms.
Using the expression of the operators $ {\mathcal A}_\ell $, $\ell = 0 , \ldots , 2$,  given by Lemma \ref{scaledlaplacian}, we get
\begin{eqnarray*}
I_0 &=& \e\, \left\{-\mu_0 \,\Delta_K
 (\mu_0) \, Z_0+\,| \nabla_K\mu_0|^2 \mathcal{T}_1(w_0)  -\mu_0^2(\mathcal{T}_2(w_0)-\mathcal{T}_3(w_0))\right\}
\\ &+& O(\e^2) b(\xi )
\end{eqnarray*}
where $b(\xi )$ is a smooth function such that $\| (1+ |\xi|^{N-2} ) b(\xi ) \|_\infty \leq C$, for some constant $C$ independent of $\e$. Furthermore,
we recall that
$
Z_0=\gamma
w_0 + D_\xi w_0 \, [\xi].
$
Also we denoted
\begin{eqnarray}\label{mathcala}
\mathcal{T}_1(w_0)=D_{\xi\xi} w_0 \, [
\xi]^2 +
2 (1+\gamma ) D_\xi w_0 [\xi] + \gamma (1+ \gamma ) w_0,
\end{eqnarray}
\begin{eqnarray}\label{mathcalb}
\mathcal{T}_2(w_0)= {1 \over 3}\,\, \sum\limits_{i,j} \bigg[\sum\limits_{m,l} R_{mijl}
 \xi_m
\xi_l \bigg]  \partial^2_{ij} w_0,\end{eqnarray}
\begin{eqnarray}\label{mathcalc}
\mathcal{T}_3(w_0)= \sum\limits_j \bigg[ \sum_s
\frac23 R_{mssj}+\sum\limits_{m,a,b}
 \big( \,\tilde g^{ab} R_{mabj}- \G_{am}^b
\G_{bj}^a \big)\bigg] \xi_m
 \partial_j w_0.
\end{eqnarray}
On the other hand, a direct computation shows that
\begin{eqnarray*}
 I_1=w_0^p\left [\mu_0^{\mp\frac{N-2}{2}\e}\, w_0^{ \pm\e}-1\right]=\pm \e\, w_0^p\ln w_0 \left( 1 +O(\e^2) \right).
\end{eqnarray*}
Thus we can write
\begin{eqnarray*}
\Xi_\e(v_1)
& = &
  - \Delta_{\R^{N}} w_{1}
  - p w_0^{p-1} w_{1}+\e \mu_0^2 \, h \, w_{1}
+\e\, H_1 (z, \xi) + \e L_\e (w_1 )
  + Q_\e (w_{1 }),
\end{eqnarray*}
where
\begin{eqnarray*}
H_1 (z,\xi ) &=& \mu_0 \,\Delta_K
 (\mu_0) \, Z_0-\,| \nabla_K\mu_0|^2 \mathcal{T}_1(w_0)  +\mu_0^2(\mathcal{T}_2(w_0)-\mathcal{T}_3(w_0))\\
& & \\
& &+  \mu_0^2 h w_0\mp w_0^p\ln (w_0)+ \mathcal{E}_{1,\e},
\end{eqnarray*}
with $\mathcal{E}_{1, \e}$ is a sum of
functions of the form
$$
\e \mu_0 \left( \e \mu_0 + \e \partial_a \mu_0 + \e \partial^2_a
\mu_0 \right) a (z) b (\xi )
$$
and $a(\e z)$ is a smooth function uniformly bounded, together
with its derivatives,  as $ \e \to 0$, while the function $b$ is
such that
$$
\sup_{\xi } (1+|\xi|^{N-2} ) |b (\xi ) |< \infty.
$$
The term $ \e L_\e (w_1 )$ is linear in $w_1$, infact it is explicitely given by
\begin{equation} \label{mi1}
\e L_\e (w_1) =I_2 -\sum_{\ell=0}^2 {\mathcal A}_\ell w_{1,\e} - B(w_{1,\e})
\end{equation}
The term $Q_\ve (w_{1,\e} )$ is quadratic in $w_{1,\e}$, in fact it
is explicitly given by
\begin{equation}\label{mi2}
\mu_0^{\mp \frac{N-2}{2}\e}\left[(w_0 + w_{1,\e} )^{p\pm \e} - w_0^{p\pm \e} - p w_0^{p-1\pm \e} w_{1,\e}\right].
\end{equation}
We ask the function  $w_{1,\e}$ to
satisfy the following equation
\begin{equation}\label{eq:eqw1}
    - \Delta_{\R^{N}} w_{1,\e}
  - p w_0^{p-1} w_{1,\e}+\e \mu_0^2 \, h \, w_{1,\e}
=-\e\, H_1 (z, \xi),   \  \hbox{ in } \D , \quad \phi = 0 \quad {\mbox {on}} \quad \partial \hat \D.
\end{equation}
Using Proposition \ref{linear}, we see that equation \eqref{eq:eqw1} is
solvable if the right-hand side satisfies the orthogonality conditions in (\ref{eq:eqwd}). These conditions,
for $j = 1, \dots, N$ are clearly satisfied  since both
$\xi_j\partial_{j} w_0$ and $\partial^2_{ij} w_0$ are even functions in $ \xi$,
while the $Z_i$'s are odd functions in $\xi$ for every $i$. It remains to
compute the $L^2$ product of the right-hand side against $Z_0$. Imposing this $L^2$ product equal to zero will define the function $\mu_0$.

We define $\mu_0$ to satisfy, at main order,
\begin{equation}\label{chmu0}
\int_{\hat \D} H_1 (z, \xi ) Z_0 (\xi ) d\xi= 0 \quad \forall z \in K_\e.
\end{equation}
Let us be more precise. We have
\begin{eqnarray*}
&&\int_{\hat \D} H_1 (z, \xi ) Z_0 (\xi ) d\xi= \mu_0(y) \,\Delta_{K}
 (\mu_0) \,\int_{\R^N} Z^2_0(\xi)d\xi\\
 && -\,| \nabla_K\mu_0 |^2\int_{\R^N}Z_0(\xi) \mathcal{T}_1 (w_0)d\xi  +\mu_0^2(y)\int_{\R^N}(\mathcal{T}_2(w_0)-\mathcal{T}_3(w_0))Z_0(\xi)d\xi\\
& &+  \mu_0^2(y)h(y)\int_{\R^N}  w_0(\xi)Z_0(\xi)d\xi \mp \int_{\R^N}w_0^p\ln (w_0)Z_0d\xi+O\left(\left(\frac{\e}{\eta^2}\right)^{\frac{N-4}{2}}\right).
\end{eqnarray*}
Define
\begin{eqnarray}\label{c1}
c_{1,N}:=\int_{\R^N} Z^2_0(\xi)d\xi
&=&\alpha_N^2\frac{(N-2)^2(N+2)}{2N(N-4)}\omega_N\, I_{N}^{N/2}>0.
\end{eqnarray}
A direct computation gives that
\begin{eqnarray*}
\int_{\R^N}Z_0(\xi) \mathcal{T}_1(w_0)d\xi =0.
\end{eqnarray*}
Moreover
\begin{eqnarray*}
\int_{\R^N}Z_0(\xi) \mathcal{T}_2(w_0)d\xi &=&{1 \over 3}\,\,  \sum\limits_{i,j}  \sum\limits_{m,l} R_{mijl}
\int_{\mathbb{R}^N} \xi_m
\xi_l   \partial^2_{ij} w_0 Z_0d\xi\nonumber\\&=&{1 \over 3}\,\,  \sum\limits_{i,j}  R_{jiij}
\int_{\mathbb{R}^N}
\xi_j    \partial_{j} w_0 Z_0d\xi,
\end{eqnarray*}
because $R_{mijl}$ is antisymmetric (i.e. $R_{mijl}=-R_{imjl}$) and
\begin{eqnarray*}\int_{\mathbb{R}^N} \xi_m
\xi_l  \partial^2_{ij} w_0 Z_0d\xi& =& \alpha_N(N-2) \int_{\mathbb{R}^N} \xi_m
\xi_i    \left(-{\delta_{ij}\over (1+|\xi|^2)^{N\over2}}+ {N \xi_i\xi_j\over(1+|\xi|^2)^{N+2\over2}}\right) Z_0d\xi
\end{eqnarray*}
and $\int_{\mathbb{R}^N}  {\xi_m\xi_l\xi_i\xi_j\over(1+|\xi|^2)^{N+2\over2}}  Z_0d\xi$
 is symmetric.
On the other hand,
\begin{eqnarray}\label{c2from}
\int_{\R^N}Z_0(\xi) \mathcal{T}_3(w_0)d\xi &=& \sum\limits_j \bigg[ \sum_s
\frac23 R_{mssj}+\sum\limits_{m,a,b}
 \big( \, {\tilde g}^{ab}  R_{mabj}- \G_{am}^b
\G_{bj}^a \big)\bigg]
\int_{\mathbb{R}^N} \xi_m
 \partial_j w_0 Z_0d\xi\nonumber
 \\&=& \sum\limits_j \bigg[ \sum_s
\frac23 R_{jssj}+\sum\limits_{j,a,b}
 \big( \, {\tilde g}^{ab}  R_{jabj}- \G_{aj}^b
\G_{bj}^a \big)\bigg]
c_{2,N}
\end{eqnarray}
where
\begin{eqnarray}\label{c2}
c_{2,N}:=
\int_{\mathbb{R}^N} \xi_j
 \partial_j w_0 Z_0d\xi= \alpha_N^2\frac{3(N-2)^2}{2N(N-4)}\omega_N\, I_{N}^{N/2}>0.
\end{eqnarray}
From the above computations we deduce
\begin{eqnarray*}\int_{\R^N}Z_0(\xi) \left(\mathcal{T}_2(w_0)- \mathcal{T}_3(w_0)\right)d\xi = c_{2,N} \left[ \sum_{i,j}
\frac13 R_{ijij}+\sum\limits_{j}\sum\limits_{a,b}
 \big( \,{\tilde g}^{ab}  R_{jajb}+ \G_{aj}^b
\G_{bj}^a \big)\right],\end{eqnarray*}

We also have
\begin{eqnarray*}
\int_{\hat \D}  \ln (w_0)w_0^p\, Z_0d\xi& = &\int_{\mathbb{R}^N}  \ln (w_0)w_0^p\, Z_0d\xi-\int_{\mathbb{R}^N\backslash\hat \D}  \ln (w_0)w_0^p\, Z_0d\xi\\
& = &\int_{\mathbb{R}^N}  \ln (w_0)w_0^p\, Z_0d\xi+O\left(\left(\frac{\e}{\eta}\right)^{\frac{N-2}{2}}\right)\\
&: = &c_{4,N}  +O\left(\left(\frac{\e}{\eta}\right)^{\frac{N-2}{2}}\right),
\end{eqnarray*}
where
\begin{eqnarray}\label{c4}
c_{4,N}:=
\frac{N}{(p+1)^2}\int_{\mathbb{R}^N}  w_0^{p+1}(\xi)d\xi=\alpha_N^{p+1}\frac{ (N-2)^3}{4N^2}\omega_N\, I_{N}^{N/2}>0.
\end{eqnarray}
Finally, set
\begin{eqnarray}\label{c3}
c_{3,N}:=
\int_{\R^N}  w_0(\xi)Z_0(\xi)d\xi= -\alpha_N^2\frac{2(N-1)(N-2)}{N(N-4)}\omega_N\, I_{N}^{N/2}<0.
\end{eqnarray}
We define $\mu_0$ to satisfy
\begin{eqnarray}\label{choiceofmu0}
 - \Delta_K\mu_0+ {\frac{c_{3,N}}{c_{1,N}}} \, {\mathcal H} (y)\mu_0\pm\frac{c_{4,N}}{c_{1,N}}\frac{1}{\mu_0}=0,\quad  {\mbox {in}} \quad  \ K.
\end{eqnarray}
where  ${\mathcal {H}}$  is the function defined in Fermi coordinates by
\begin{equation}\label{defgy}
{\mathcal H} (y) =   h(y) -\hat \Omega (y)
\end{equation}
where
$$
\hat \Omega (y) = -  {3(N-2)\over 4(N-1)}\[\sum\limits_{i,j=1}^{N}
\frac13 R_{jiji}(y)+ \sum\limits_{i=1}^{N}\sum\limits_{a,b=1}^k
 \(\tilde g^{ab}R_{iaib}(y)+ \Gamma_{ai}^b(y)
\Gamma_{bi}^a(y) \)\] .
$$
The existence of $\mu_0$ is guaranteed by our assumption. With this choice for $\mu_0$, the integral of  the right hand side
in (\ref{eq:eqw1}) against $Z_{0}$ vanishes  on $K$ and this implies
the existence of $w_{1,\e}$, thanks to Proposition \ref{linear}.
Moreover, it is straightforward to check that
$$
\|  H_1 (z, \xi ) \|_{\e , N-2 , \sigma } \leq C
$$
for some $\sigma \in (0,1)$. Proposition \ref{linear} thus gives that
\begin{equation}\label{ew1}
\| D^2_\xi w_{1,\e} \|_{\e , N-2 , \sigma } + \| D_\xi
w_{1,\e} \|_{\e , N-3 , \sigma } +\|w_{1,\e} \|_{\e, N-4 ,
\sigma}\le C \e
\end{equation}
and that there exists a positive constant $\beta$ (depending only on
$ K$ and $N$) such that for any integer $\ell$ there holds
\begin{equation}\label{eq:estw1}
    \|\nabla^{(\ell)}_{z} w_{1,\e}(z,\cdot)\|_{\e,N-4,\sigma} \leq \beta C_l
    \e \qquad \quad \sqrt{\e}z \in K
\end{equation}
where $C_l$ depends only on $l$, $p$ and $K$.

\medskip
With this choice of $\mu_{0}$ and $w_{1, \e}$ we get that
\begin{equation}
\label{mi3}
\| - {\mathcal A}_{\mu_\e , \Phi_\e } v_{1,\e} + \e\mu_0^2 hv_{2,\e}
-  \mu_0^{\mp \e\frac{N-2}{2}}\,v_{1,\e}^{p\pm \e} \|_{\e, N-2 , \sigma} \leq C\e.
\end{equation}
To prove this, we first observe that
$$
\| \mu_0^2 \Delta_{K_\e} w_{1,\e} \|_{\e, N-2 , \sigma} \leq C\e,
$$
as consequence of (\ref{eq:estw1}).
Then we
claim that $\| \e L_\e (w_{1 , \e} ) + Q_\e (w_{1,\e} ) \|_{\e , N-2 } \leq C \e$, see (\ref{mi1}) and (\ref{mi2}). Indeed, first observe that
\begin{eqnarray*}
I_2= (p\pm\e)\mu_0^{\mp\frac{N-2}{2}\e}\,w_0^{p-1\pm\e}w_{1,\e}-pw_0^{p-1}w_{1,\e} = O(\e)w_0^{p-1}w_{1,\e}.
\end{eqnarray*}
and then easily we get that $|I_2 | \leq C {\e \over (1+ |\xi | )^{N-2}}$.  Analogous consideration gives that $\| Q_\e (w_{1 , \e} ) \|_{\e , N-2} \leq C \e$. Furthermore, using the fact that
$|\xi | \leq \eta \e^{-{1\over 2}}$, we have that
$$
\left| \sum_{\ell =0}^2 {\mathcal A}_\ell  (w_{1,\e} )+ B(w_{1,\e} ) \right| \leq C  {\e \over (1+ |\xi | )^{N-2}}.
$$
Estimate (\ref{mi3}) follows from the regularity of the function $w_{1, \e}$.

\medskip
\noindent
Let $I=2$. Then  $\mu_{2,\e}=\mu_0+\e  \mu_1$, $\Phi_{\e}=\Phi_{1,\e}$, and
$v_2= w_0+w_{1,\e}+w_{2,\e}$, where $\mu_0$ and $w_{1,\e}$ have already been constructed  in the previous step.
Computing $\Xi_\e(v_{2})$ ,
\begin{eqnarray}\label{sw2}
\Xi_\e(v_2)& = & - \Delta_{\R^{N}} w_{2,\e}
  - p w_0^{p-1} w_{2,\e} + \e \mu_{2,\e}^2 \, h \, w_{2,\e}
+ H_2 (z, \xi) \nonumber \\
&+&\e^2 L_\e (w_{2,\e} )
 + Q_\e (w_{2,\e}),
\end{eqnarray}
where
\begin{eqnarray*}
H_2 (z,\xi ) &=& \e^2  \,  \left\{\mu_0\, \Delta_K
 (\mu_1) \,Z_0 +\mu_1\Delta_K(\mu_0)Z_0-  \,  \nabla_K(\mu_0)\nabla_K(\mu_1)\mathcal{T}_1(w_0)\right. \\
& &\qquad\qquad \left.+2 \mu_0\,\mu_1(\mathcal{T}_2(w_0)-\mathcal{T}_3(w_0))+  2\mu_0\, \mu_1h w_0\right.\\
&&\qquad\qquad \left.\pm \frac{(N-2)^2}{16} w_0^p\left(\ln (w_0)+1\right)\right\}\\
&&+\e^{\frac{3}{2}} \, \mu_0  \left\{- \Delta_K\Phi_{1,\e}D_\xi w_0  + {1\over 3}\,\,  \sum\limits_{i,j}  \sum\limits_{m,l} R_{mijl}
( \xi_m
\Phi_{1,\e}^l +\xi_I\Phi_{1,\e}^m)   \partial^2_{ij} w_0,\right. \\
& & \qquad\qquad\left.+   \sum\limits_j \bigg[ \sum_s
\frac23 R_{mssj}+\sum\limits_{m,a,b}
 \big( \,{\tilde g}^{ab}  R_{mabj}- \G_{am}^b
\G_{bj}^a \big)\bigg] \Phi_{1,\e}^m
 \partial_j w_0  \right\}\\
&&+ \e^{\frac{3}{2}}  \mathcal{E}_{2,\e}(y,\xi,w_0,w_{1,\e},\mu_0),
\end{eqnarray*}
and $\mathcal{E}_{2, \e}$ is a sum of
functions of the form
$$
  \left(   \mu_0 +   \partial_a \mu_0 +   \partial^2_a
\mu_0 \right) a (\sqrt{\e }z) b (\xi )
$$
and $a(\sqrt{\e }z)$ is a smooth function uniformly bounded, together
with its derivatives,  as $ \e \to 0$, while the function $b$ is
such that
$$
\sup_{\xi } (1+|\xi|^{N-2} ) |b (\xi ) |< \infty.
$$
Furthermore, we have
\begin{eqnarray*}
\e^2 L_\e (w_2)& = & -  \mu_0^{\mp\frac{N-2}{2}\e}\,\left[ (w_0+w_{1,\e} + w_{2,\e} ) ^{p+\e} \right. \\ & & \left.-(w_0 + w_{1,\e})^{p\pm\e}-(p\pm\e) (w_0 + w_{1,\e})^{p-1\pm\e}w_{2,\e}\right]\\
& - &\sum_{\ell=0}^2 {\mathcal A}_\ell w_{2,\e} - B(w_{2,\e})
\end{eqnarray*}
and the term $Q_\e (w_{2,\e} )$ in (\ref{sw2}) is a sum of quadratic terms in
$w_{2,\e}$ like
$$
(\mu_0+\e\mu_1)^{\mp \frac{N-2}{2}\e}\left[(w_0 + w_{1,\e} + w_{2,\e} )^{p\pm \e} - (w_0 + w_{1,\e} )^{p\pm \e}\right.
$$
$$
\left.  -(p\pm \e) (w_0 +
w_{1,\e} )^{p-1\pm \e}  w_{2,\e}\right].
$$

We will choose $w_{2,\e}$ to
satisfy the following equation
\begin{equation}\label{eq:eqw2}
    - \Delta_{\R^{N}} w_{2,\e}
  - p w_0^{p-1} w_{2,\e}+\e \mu_{2,\e}^2 \, h \, w_{2,\e}
=- H_2 (z, \xi),     \hbox{ in } \D  , \quad w_{2, \e} = 0 \quad {\mbox {on}} \quad \partial \hat \D.
\end{equation}
Thanks to Proposition \ref{linear}, we see that equation \eqref{eq:eqw2} is
solvable if the right-hand side is $L^2$-orthogonal to the functions $Z_j$, for $j=0, \ldots , N$.
These orthogonality conditions will define the
parameters $\mu_{1}$ and the normal section $\Phi_{1,\e}$.

\medskip
\noindent  Projection onto $Z_0$ and choice of
$\mu_{1}$:\ \   the function $\mu_1$ is asked to satisfy, at main order,
\begin{eqnarray*}
 \int_{\hat D } H_{2,\e}Z_0d\xi&=& 0.
\end{eqnarray*}
Computations similar to the ones already performed to define $\mu_0$ give  that $\mu_1$ satisfies
\begin{equation}\label{mu11}
 -\mu_0(y)\, \Delta_K
 \mu_1(y)   -\mu_1\Delta_K \mu_0(y) +2\mu_0(y)\,\mu_1(y)g(y) \pm \frac{(N-2)^2}{16}\frac{c_{4,N}}{c_{1,N}}=0
\end{equation}
 in $K$,
where $g(y)$ is given by (\ref{defgy}), and $c_{i,N}(i=1,2,3,4)$ are defined in (\ref{c1}), (\ref{c2}), (\ref{c3}) and (\ref{c4}).
According to our choose of $\mu_0$ satisfies (\ref{choiceofmu0}), then (\ref{mu11}) is equivalent to the following equation
\begin{eqnarray}\label{mu11a}
 - \Delta_K\mu_1(y)+g(y)\mu_1(y)\mp\frac{c_{4,N}}{c_{1,N}}\frac{1}{\mu_0^2(y)}\mu_1(y)=\mp\frac{(N-2)^2}{16}\frac{c_{4,N}}{c_{1,N}}\frac{1}{\mu_0(y)},\quad
\end{eqnarray}
 in  $K$.
The existence of $\mu_1$ is guaranteed by the nondegeneracy of $\mu_0$.

\medskip	
\noindent  Projection onto $Z_s$ and choice of
$\Phi_{1,\e}$. \ \    Multiplying $H_{2,\e}$ with $Z_s = \partial_s w_0$,
integrating over $\hat\D$ and using the fact $w_0$ is even in the
variable $\xi$, one obtains
\begin{eqnarray}\label{projZl}
 \int_{\hat\D}H_{2,\e}\,\partial_s w_0d\xi &=& - \e^{\frac{3}{2}} \,  \mu_0 \, \sum_j \, \Delta_K \Phi_{1,\e}^j \,\int_{\hat\D} \partial_jw_0 \partial_s w_0d\xi +\e^{\frac{3}{2}} \, \int_{\hat\D}\mathcal{E}_{2, \e}\partial_sw_0d\xi\nonumber\\[2mm]
&& +\e^{\frac{3}{2}} \, \, { \mu_0 \over 3}\,\,  \sum\limits_{i,j}  \sum\limits_{m,l} R_{mijl}\int_{\hat\D}
( \xi_m
\Phi_{1,\e}^l +\xi_l \Phi_{1,\e}^m)   \partial^2_{ij} w_0\partial_s w_0d\xi \nonumber \\[2mm]
&&+ \e^{\frac{3}{2}} \, {2  \mu_0 \over 3} \,\sum\limits_{j , m} \sum_l
 R_{mllj}  \Phi_{1,\e}^m
 \int_{\hat\D}\partial_jw_0\,\partial_sw_0 \\[2mm]
&&+ \e^{\frac{3}{2}} \,  \,\sum\limits_{j , m}
 \sum\limits_{a,b}
 \big( \, {\tilde g}^{ab}  R_{mabj}- \G_{am}^b
\G_{bj}^a \big) \Phi_{1,\e}^m
 \int_{\hat\D}\partial_jw_0\,\partial_sw_0.\nonumber
\end{eqnarray}
First of all, observe that  by oddness in $\xi$ we have that
$$
\int_{\hat\D} \pa_jw_0 \pa_sw_0=
\d_{js}\, C_0+O(\e^{\frac{N-2}{2}}), \quad C_0:= \int_{\R^N} |\pa_lw_0|^2d\xi.
$$
Thus
$$
- \e^{\frac{3}{2}} \,  \mu_0 \, \sum_j \, \Delta_K \Phi_{1,\e}^j \,\int_{\hat\D} \partial_jw_0 \partial_s w_0d\xi  = - \e^{\frac{3}{2}} \,  \mu_0   \, C_0 \, \Delta_K \Phi_{1,\e}^s  + O(\e^{\frac{N-2}{2}})
$$
and
$$
 \e^{\frac{3}{2}} \, {2  \mu_0 \over 3} \,\sum\limits_{j , m} \sum_l
 R_{mllj}  \Phi_{1,\e}^m
 \int_{\hat\D}\partial_jw_0\,\partial_sw_0 =  \e^{\frac{3}{2}} \, {2  \mu_0 \over 3} \, C_0 \, \sum\limits_{ m} \sum_l
 R_{mlls}  \Phi_{1,\e}^m.
$$

On the other hand the integral $\int_{\hat\D} \xi_m
\,\partial^2_{ij}w_0  \pa_sw_0$ is non-zero only if, either $i = j$ and
$m = s$, or $i = s$ and $j = m$, or $i = m$ and $j = s$. In the
latter case we have $R_{mijs}=0$ (by the antisymmetry of the
curvature tensor in the first two indices). Therefore, the first
term of the second line in (\ref{projZl}) becomes simply
\begin{eqnarray*}
&& 2 \ve^{3\over 2} {\mu_0 \over 3} \sum_m \sum_{l} R_{mlls} (- \int \xi_l \partial^2_{ls} w_0 \partial_s w_0 + \int \xi_s \partial^2_{ll} w_0 \partial_s w_0 )  {\Phi }^m_{1,\e}  +O(\e^{\frac{N-2}{2}})\\
&& =- 2 \ve^{3\over 2} {\mu_0 \over 3} \sum_m \sum_{l} R_{mlls}  (- \int \xi_l \partial^2_{ls} w_0 \partial_s w_0 + \int \xi_s \partial^2_{ls} w_0 \partial_l w_0 )  {\Phi }^m_{1,\e}  +O(\e^{\frac{N-2}{2}}) \\
&& =- 2 \ve^{3\over 2} {\mu_0 \over 3} \, C_0 \,  \sum_m \sum_{l} R_{mlls}   {\Phi }^m_{1,\e}  +O(\e^{\frac{N-2}{2}})
\end{eqnarray*}
where the last identity is consequence of the following fact
$$
\int \xi_l \partial^2_{sl} w_0 \partial_s w_0 = -{C_0 \over 2} ,
$$
which can be proved with a straightforward computation.

Hence formula (\ref{projZl}) becomes simply
\begin{eqnarray*}
\int_{\hat\D}H_{2,\e}\,\pa_s w_0 d\xi&=& \e^{3\over 2}  \,  \mu_0 \, C_0\,\left[ - \Delta_K\,\Phi^s_{1,\e}+\sum_m  \sum_{a,b} \bigg( \tilde g^{ab}\,R_{mabs}- \Gamma_a^b(E_m) \Gamma_b^a(E_s)  \bigg)
\,\Phi^m_{1,\e} \right]\\
&& +O(\e^{N-{1\over2}} ) \, \Phi^m_{1, \e}
+\e^{3\over 2}  \int_{\hat\D}\mathcal{E}_{2,\e}\,\pa_s w_0.
\end{eqnarray*}
We thus obtain that $\int H_{2,\e}(z,\xi,  w_0,  \dots, w_{1,\e}) Z_l =0$ at main order
 if $\Phi_{1,\e}$ satisfies an
equation of the form
\begin{equation}\label{phi1defa}
  \Delta_K \,\Phi^s_{1,\e}- \sum_m \sum_{a,b} \bigg( \tilde g^{ab}\,R_{mabs}- \Gamma_a^b(E_m) \Gamma_b^a(E_s)   \bigg)
\,\Phi^m_{1,\e} = G_{2,\e}(\sqrt{\e} z),
\end{equation}
for some expression $G_{2,\e}$ smooth on its argument. Observe that
the operator acting on $\Phi_{1,\e}$ in the left hand side is
nothing but the Jacobi operator, which is invertible by the
non-degeneracy condition on $K$. This implies the solvability of the
above equation in $\Phi_{1,\e}$.
Furthermore, equation (\ref{phi1defa}) defines $\Phi_{1, \e}$ as a
smooth function on $K$, of order $\e$, more precisely we have
\begin{equation}
\label{ePhi1e}
\| \Phi_{1,\e}\|_{\infty } +\| \partial_a
\Phi_{1,\e}\|_{\infty} +\|\partial^2_a \Phi_{1,\e}\|_{\infty
} \leq C .
\end{equation}

By our choice of $\mu_{1}$ and $\Phi_{1,\e}$ we have solvability
of equation (\ref{eq:eqw2}) in $w_{2,\e}$. Moreover, it is
straightforward to check that
\begin{eqnarray*}
|H_{2,\e} ( z, \xi ) | &\leq&  C \max\left\{\e^{2 } ,\e^{\frac{3}{2}} \right\}{1\over (1+|\xi|)^{N-2}} .
\end{eqnarray*}
Furthermore, for a given $\sigma \in (0,1)$ we have
$$
\| H_{2,\e} \|_{\e, N-2, \sigma} \leq C \e^{\frac{3}{2}} .
$$
Proposition \ref{linear} thus gives then that
\begin{equation}
\label{ew2}
 \| D^2_\xi w_{2,\e} \|_{\e , N-2 , \sigma } + \| D_\xi
w_{2,\e} \|_{\e , N-3 , \sigma } +\|w_{2,\e} \|_{\e, N-4,
\sigma}\le C \e^{\frac{3}{2}}
\end{equation}
and that there exists a positive constant $\beta$ (depending only on
$\Omega, K$ and $n$) such that for any integer $\ell$ there holds
$$
    \|\nabla^{(\ell)}_{z} w_{2,\e}(z,\cdot)\|_{\e,N-4,\sigma} \leq \beta C_\ell \,
   \e^{\frac{3}{2}},
$$
where $C_\ell$ depends only on $\ell$, $p$, $K$.

\medskip
\noindent
Arguing as in the previous step for $I=1$, we see that with this choice of $\mu_{1, \e}$, $\Phi_{1,\e}$ and $w_{2, \e}$ we get that
$$
\| - {\mathcal A}_{\mu_\e , \Phi_\e } v_{2,\e} + \e\mu_\e^2 hv_{2,\e}
-  \mu_\e^{\mp \e\frac{N-2}{2}}\,v_{2,\e}^{p\pm \e} \|_{\e, N-2 , \sigma} \leq C\e^{\frac{3}{2}}.
$$

\medskip
\noindent Expansion at an arbitrary order. \ \  We take
now an arbitrary integer $I$. Let
\begin{equation}\label{eqmu}
 \mu_{I+1,\e}:= \mu_0+\e \mu_{1}+\e^2\mu_2\cdots +\e^{I-1}  \mu_{I-1} +\e^{I}  \mu_{I} ,
\end{equation}
\begin{equation}\label{eqPhi}
\Phi_{\e}=\Phi_{1,\e}+\cdots +\Phi_{I-1,\e} + \Phi_{I, \e}
\end{equation}
and
\begin{equation}
\label{eqWWW}
v_{I+1, \e} =  w_0 (\xi ) + w_{1, \e} (z, \xi) +
\ldots + w_{I, \e} (z, \xi ) + w_{I+1, \e} (z, \xi )
\end{equation}
where $\mu_0,  \mu_{1} , \cdots ,  \mu_{I-1,\e}$, $\Phi_{1,\e},
\cdots , \Phi_{I,\e}$ and $w_{1, \e} $, .. , $w_{I, \e}$ have
already been constructed following an iterative scheme, as described
in the previous steps of the construction.

In particular one has, for any $i=1, \ldots , I-1$
\begin{equation}
\label{emuie} \| \mu_{i}\|_{\infty} +\| \partial_a
\mu_{i}\|_{\infty} +\|\partial^2_a \mu_{i}\|_{\infty
} \leq C
\end{equation}
\begin{equation}
\label{ePhiie} \| \Phi_{i,\e}\|_{\infty } +\| \partial_a
\Phi_{i,\e}\|_{\infty} +\|\partial^2_a \Phi_{i,\e}\|_{\infty
} \leq C \e^{i-1}
\end{equation}
We have
\begin{eqnarray}
\label{ewi}
 \| D^2_\xi w_{i,\e} \|_{\e , N-2 , \sigma } + \|
D_\xi w_{i,\e} \|_{\e , N-3 , \sigma } +\|w_{i,\e} \|_{\e, N-4
, \sigma} \leq  C  \e^{i-\frac{1}{2}}.
\end{eqnarray}
and, for any integer $\ell$, $\sqrt{\e}z\in K$,
\begin{equation}\label{eq:estwi}
    \|\nabla^{(\ell)}_{z} w_{i,\e}(z,\cdot)\|_{\e,N-2,\sigma} \leq \beta C_l \e^{i-\frac{1}{2}}.
\end{equation}

The new triplet $(\mu_{I} , \Phi_{I, \e} , w_{I+1 , \e} )$
will be found reasoning as in the construction of $(\mu_{1} ,
\Phi_{1} , w_{2 , \e} )$. Computing $\Xi_\e (v_{I+1 , \e })$  we get
\begin{eqnarray}\label{rivoli5a}
&&-\hat {\mathcal A}_{\mu_\e ,
\Phi_\e } v_{I+1,\e} + \e\mu_{I+1,\e}^2 \, hv_{I+1,\e}-  \mu_{I+1,\e}^{\mp\e\frac{N-2}{2}}v_{I+1,\e}^{p\pm\e}\\[3mm]
&=& - \Delta_{\R^{N}} w_{I+1,\e}
  - p w_0^{p-1} w_{I+1,\e}+\e \mu_{I+1,\e}^2 \, h \, w_{I+1,\e}
+ H_{I+1} (z, \xi)
 + Q_\e (w_{I+1,\e}) \nonumber
\end{eqnarray}
where the function $H_{I+1,\e}$ is given by
\begin{eqnarray}\label{wIepsilonl}
H_{I+1,\e}&=&  \e^{I+1}  \left\{\mu_0\, \Delta_K
 (\mu_I) \,Z_0 +\mu_I\Delta_K(\mu_0)Z_0-  \,  \nabla_K(\mu_0)\nabla_K(\mu_I)\mathcal{T}_1(w_0)\right. \nonumber\\
& &\qquad\qquad \left.+2 \mu_0\,\mu_I(\mathcal{T}_2(w_0)-\mathcal{T}_3(w_0))+  2\mu_0\, \mu_Ih w_0\right.\nonumber\\
&&\qquad\qquad \left.\pm \frac{(N-2)^2}{16} w_0^p\left(\ln (w_0)+1\right)\right\}\nonumber\\
&&+\e^{I+\frac{1}{2}} \, \mu_0  \left\{ -\Delta_K\Phi_{I,\e}D_\xi w_0  + {1\over 3}\,\,  \sum\limits_{i,j}  \sum\limits_{m,l} R_{mijl}
( \xi_m
\Phi_{I,\e}^l +\xi_I\Phi_{I,\e}^m)   \partial^2_{ij} w_0,\right. \nonumber\\
& & \qquad\qquad\left.+   \sum\limits_j \bigg[ \sum_s
\frac23 R_{mssj}+\sum\limits_{m,a,b}
 \big( \, R_{mabj}- \Gamma_{am}^b
\Gamma_{bj}^a \big)\bigg] \Phi_{I,\e}^m
 \partial_j w_0  \right\}\nonumber\\
&&+ \e^{I+\frac{1}{2}}\mathcal{E}_{I+1,\e}(y,\xi,w_0,w_{1,\e},\cdots, w_{I,\e},\mu_0,\cdots,\mu_{I-1}).
\end{eqnarray}
In (\ref{wIepsilonl}), $\mathcal{E}_{I+1, \e}$ is a sum of
functions of the form
$$
\biggl[   \left(   \mu_0 +   \partial_a \mu_0 +   \partial^2_a
\mu_0 + \ldots +  \mu_{I-1} +   \partial_a \mu_{I-1} +   \partial^2_a
\mu_{I-1} \right)
$$
$$
 +  \left(   \Phi_1 +   \partial_a \Phi_1 +   \partial^2_a
\Phi_1 + \ldots +  \Phi_{I-2} +   \partial_a \Phi_{I-2} +   \partial^2_a
\Phi_{I-2} \right)   \biggl]  a(\sqrt{\e }z) b(\xi )
$$
where $a(\sqrt{\e} z)$ is a smooth function uniformly bounded, together
with its derivatives,  as $ \e \to 0$, while the function $b$ is
such that
$$
\sup_{\xi } (1+|\xi|^{N-2} ) |b (\xi ) |< \infty.
$$
Finally the term $Q_\e (w_{I,\e} )$ in (\ref{rivoli5a}) is a sum of quadratic terms in
$w_{I+1,\e}$ like
$$
(\mu_0+\e\mu_1+\cdots+\e^{I}\mu_I)^{\mp \frac{N-2}{2}\e}\times
$$
$$\left[(w_0 + w_{1,\e} + w_{2,\e}+\cdots+w_{I+1,\e} )^{p\pm \e} \right.
$$
$$
\left.- (w_0 + w_{1,\e} +\cdots+w_{I,\e})^{p\pm \e}\right.
$$
$$
\left.  -(p\pm \e) (w_0 +
w_{1,\e} +\cdots+w_{I,\e})^{p-1\pm \e}  w_{I+1,\e}\right]
$$
and linear terms in $w_{I,\e}$ multiplied by a term of order
$\e$, like
$$
p \left( ( w_0 + w_{1,\e} + \ldots + w_{I-1, \e}  )^{p-1\pm \e} - w_0^{p-1\pm \e} \right) w_{I,\e}.
$$

Consider the following problem
\begin{eqnarray}\label{eq:eqwI}
     - \Delta_{\R^{N}} w_{I+1,\e}
  - p w_0^{p-1} w_{I+1,\e}+\e \mu_{I+1,\e}^2 \, h \, w_{I+1,\e} & & =-H_{I+1,\e} ( z , \xi )   \quad\hbox{ in } \D \nonumber \\
 \quad w_{I+1, \e} = 0 & & \quad {\mbox {on}} \quad \partial \hat \D.
\end{eqnarray}
Again by Proposition \ref{linear}, the above problem is solvable in $w_{I+1 , \e}$  if $H_{I+1,\e}$ is $L^2$-orthogonal to $Z_j$, $j=0,
1,\cdots,N$. These orthogonality conditions will define the
parameters $\mu_{I}$ and the normal section $\Phi_{I,\e}$.

\medskip
\noindent   Projection onto $Z_0$ and choice of
$\mu_{I}$. \ \  We define $\mu_{I}$  to make, at main order,
\begin{eqnarray*}
 \int_{\hat\D} H_{I+1,\ve}Z_0d\xi=0.
\end{eqnarray*}
 The above
relation defines $\mu_{I}$ as a smooth function of $\sqrt{\e} z $ in
$K$.
From estimates (\ref{ewi})  we get that
\begin{equation}
\label{emu1e}
\| \mu_{I}\|_{ \infty } +\| \partial_a
\mu_{I}\|_{ \infty } +\|\partial^2_a \mu_{I}\|_{ \infty
} \leq C
\end{equation}

\medskip
\noindent   Projection onto $Z_s$ and choice of
$\Phi_{I,\e}$. \ \  Multiplying $H_{I+1,\e}$ with $\pa_s w_0$,
integrating over $\hat\D$ and arguing as in the construction of
$\Phi_{I, \e}$, we get
\begin{eqnarray*}
 &&\left(C_0\e^{I+\frac{1}{2}} \, \mu_0 \right)^{-1}\int_{\hat\D}H_{I+1,\e}\,\pa_s w_0 d\xi=  -  \Delta_K \,\Phi^l_{I,\e}\\
 & +&
 \sum_m \sum_{a,b}  \,\bigg(\tilde g^{ab}\, R_{maal}- \Gamma_a^b(E_m) \Gamma_b^a(E_l) + O(\e^{\frac{N-2}{2}}) \bigg)
\,\Phi^m_{I,\e}+  \int_{\hat\D}\mathfrak{G}_{I+1,\e}\,\pa_l
w_0.
\end{eqnarray*}
We then conclude that $H_{I+1,\e}(z,\xi,  w_0,  \dots, w_{I,\e})$, the
right-hand side of (\ref{eq:eqwI}), is $L^2$-orthogonal to $Z_l$
($l=1,\cdots,N$) if and only if $\Phi_{I,\e}$ satisfies an
equation of the form
\begin{eqnarray}\label{phi1def}
&&  \Delta_K \,\Phi^l_{I,\e}- \sum_m \sum_{a,b} \bigg( \tilde g^{ab}\,R_{maal}- \Gamma_a^b(E_m) \Gamma_b^a(E_l) + O(\e^{\frac{N-2}{2}}) \bigg)
\,\Phi^m_{I,\e} \nonumber\\
&=& \e^{I+\frac{1}{2}}  G_{I+1,\e}(\sqrt{\e} z),
\end{eqnarray}
where $G_{I+1,\e}$ is a smooth function on $K$, uniformly bounded as $\e \to 0$. Using again the
non-degeneracy condition on $K$ we have solvability of the above
equation in $\Phi_{I,\e}$. Furthermore,  taking into account
(\ref{rivoli5a}), we get
\begin{eqnarray}
\label{ePhiIe}
\| \Phi_{I,\e}\|_{ \infty  } +\| \partial_a
\Phi_{I,\e}\|_{ \infty } +\|\partial^2_a \Phi_{I,\e}\|_{ \infty
 } \leq C \e^{I+\frac{1}{2}}.
\end{eqnarray}

By our choice of $\mu_{I+1}$ and $\Phi_{I+1,\e}$ we have
solvability of equation (\ref{eq:eqwI}) in $w_{I+1,\e}$. Moreover,
it is straightforward to check that
$$
|H_{I+1,\e} (\e z, \xi ) | \leq C \e^{I+\frac{1}{2}}{1 \over
(1+|\xi|)^{N-2}} .
$$
Furthermore, for a given $\sigma \in (0,1)$ we have
$$
\| H_{I+1,\e} \|_{\e, N-2, \sigma} \leq C \e^{I+\frac{1}{2}}.
$$
Proposition \ref{linear} gives then that
\begin{eqnarray}
\label{ewI}
\| D^2_\xi w_{I+1,\e} \|_{\e , N-2 , \sigma } + \|
D_\xi w_{I+1,\e} \|_{\e , N-3 , \sigma } +\|w_{I+1,\e} \|_{\e, N-4
, \sigma}\leq  C \e^{I+\frac{1}{2}},
\end{eqnarray}
and that there exists a positive constant $\beta$ (depending only on
$K$ and $N$) such that for any integer $\ell$ there holds, for $\sqrt{\e}z\in K$,
$$   \|\nabla^{(\ell)} w_{I+1,\e}(z,\cdot)\|_{\e,N-2,\sigma} \leq \beta C_l
   \e^{I+\frac{1}{2}}.
$$

\medskip
\noindent
With this choice of $\mu_{I}$, $\Phi_{I,\e}$ and $w_{I+1 , \e}$ we obtain that
\begin{eqnarray*}
\| -{\mathcal A}_{\mu_\e ,
\Phi_\e } v_{I+1,\e} + \e\mu_{I+1,\e}^2 \, hv_{I+1,\e}- \mu_\e^{\mp\e\frac{N-2}{2}}v_{I+1,\e}^{p\pm\e} \|_{\e , N-2 , \sigma}\leq C \e^{I+\frac{1}{2}}.
\end{eqnarray*}

\medskip
\noindent
This concludes our construction.

\end{proof}

\section{A global approximation}\label{s:linear}
\smallskip

Let us recall that if $u$ is a solution to problem (\ref{p}), and we define
$$ u(x) =(1+\alpha_\e) \e^{-{N-2 \over 4}}
\tilde{u}(\e^{-{1 \over 2}} x ).
$$
Then $\tilde{u}$ satisfies the following equation
\begin{equation} \label{changea}
-\Delta_{g^\e} \tilde{u} +\e h \tilde{u}=  \tilde{u}^{{N+2\over N-2} \pm\e} \quad
\mbox{ in }
\M_\e ; \quad
 \tilde{u}>0 \mbox{ in }\M_\e,
\end{equation}
 where $\Delta_{g^\e}$ denotes the Laplace-Beltrami operator on $\M_\e$ is given by
 \begin{equation*} \label{laplacee}
\Delta_{g^\e}=\frac{1}{\sqrt{\mbox{det}\ g^\e}}\partial_A(\sqrt{\mbox{det}\ g^\e}(g^\e)^{AB}\partial_B)
 \end{equation*}
 here indices $A$ and $B$ run between $1$ and $n=N+k$, and $g^\e$ is the scaled metric on $\M_\e$ whose coefficient are defined by
\begin{equation*} \label{ge}
 g^\e_{\alpha,\beta}(z,x)=g_{\alpha,\beta}(\sqrt{\e} z,\sqrt{\e} x)
 \end{equation*}
 where $g_{\alpha,\beta}$ are the coefficients of the metric $g$ on $\M$.

Let $\mu_\e (y)  $, $\Phi_\e (y)$ and $v_{I+1 , \e}$ be the functions whose existence and
properties have been established in Lemma \ref{Construction}. We
define locally around $K_\e := {K\over \sqrt{\e}}$  the function
\begin{equation} \label{Vdef}
\tilde{U}_\e (z, x):= \, \mu_{\e}^{-{N-2 \over 2}} (\sqrt{\e} z) \, v_{I+1 , \e}
\left( z, \,  \frac{x-\sqrt{\e} \Phi_\e (\sqrt{\e} z)}{ \mu_\e (\sqrt{\e} z)} \right) \chi_\e (|(x-\sqrt{\e} \Phi_\e (\sqrt{\e} z) )|)
\end{equation}
where $z \in K_\e$. The function $\chi_\e$
is a smooth cut-off function with
\begin{equation}\label{magaly}
 \chi_\e (r) =
\left\{
 \begin{array}{ll}
     1, & \hbox{for} \quad r \in [0,2 \e^{-\gamma} ] \\[3mm]
    0, & \hbox{for} \quad  r \in [3 \e^{-\gamma} , 4\e^{-\gamma} ],
  \end{array}
\right.
\end{equation}
and
$$
|\chi_\e^{(l)} (r) | \leq C_l \e^{l \gamma}, \quad \forall l\geq 1,
$$
for some $\gamma \in ({1\over 2} , 1)$ to be fixed later.

We will use the notation
\begin{equation} \label{defTTa}
\tilde{u}=\widetilde{{\mathcal T}}_{\mu_\e , \Phi_\e} (\tilde{v})
\end{equation}
if and only if $\tilde{u}$ and $\tilde{v}$ satisfy
\begin{equation*}
 \tilde{u}=\mu_{\e}^{-{N-2 \over 2}} (\sqrt{\e} z) \, \tilde{v}
\left( z, \,  \frac{x-\sqrt{\e} \Phi_\e (\sqrt{\e} z)}{ \mu_\e (\sqrt{\e} z)} \right) .
\end{equation*}

The function $\tilde{U}_\e$ is well defined in a small neighborhood of
$K_\e$. We will look for a solution to (\ref{changea})
of the form
$$
\tilde{u}_\e = \tilde{U}_\e +\phi.
$$
Thus $\phi$ satisfies the following problem
\begin{equation}\label{nonlinearproblem}
    -\Delta_{g^\e} \phi + \e h \phi- (p\pm\e) \tilde{U}_\e^{p\pm\e-1} \phi =S_\e (\tilde{U}_\e) + N_\e (\phi)  \quad  \text{ in } \mathcal{M}_\e,
\end{equation}
where
\begin{equation}
\label{eomegaeps}
S_\e (\tilde{U}_\e )= \Delta_{g^\e} \tilde{U}_\e - \e h\tilde{U}_\e +
\tilde{U}_\e^p
\end{equation}
and
\begin{equation}
\label{Nomegaeps}
N_\e (\phi) =  (\tilde{U}_\e + \phi )^{p\pm\e} - \tilde{U}_\e^p - (p\pm\e)
\tilde{U}_\e^{p\pm\e-1} \phi.
\end{equation}
Define
$$
L_\e (\phi) = -\Delta_{g^\e} \phi + \e h \phi- (p\pm\e) \tilde{U}_\e^{p\pm\e-1} \phi.
$$
We shall solve the Non-Linear Problem
(\ref{nonlinearproblem}) by using a fixed point argument based on the
contraction Mapping Principle. To do so, we first establish some
invertibility properties of the linear problem
$$
L_\e (\phi) = f \quad {\mbox {in}} \quad \mathcal{M}_\e ,
$$
with $f\in L^2 (\mathcal{M}_\e )$. We do this  in two steps. First
we  study  the above problem in a strip close to the scaled
manifold $K_\e= {K\over \sqrt{\e}}$.
Let $\gamma \in ({1\over 2} , 1)$ be the number fixed before in
(\ref{magaly}) and consider
\begin{equation}
\label{omegaepsilongamma}
\mathcal{M}_{\e, \gamma} := \{ x \in \mathcal{M}_\e
\, : \, {\mbox {dist}}_{g^\e} (x,K_\e ) <2 \e^{-\gamma} \}.
\end{equation}

We are first interested in solving the following problem: given $f
\in L^2 (\mathcal{M}_{\e , \gamma})$
\begin{equation} \label{lineare}
    -\Delta_{g^\e} \phi + \e h \phi- (p\pm \e) \tilde{U}_\e^{p\pm \e-1} \phi =f\quad \text{ in } \mathcal{M}_{\e , \gamma}.
\end{equation}

We have the validity of the following result.
\begin{proposition}\label{teouffa}
There exist a constant $C>0$ and a sequence $\e_l = \e \to 0$ such
that, for any $f \in L^2 (\mathcal{M}_{\e , \gamma} )$ there exists a
solution $\phi \in H^1_{g^\e}(\mathcal{M}_{\e , \gamma} ) $ to Problem (\ref{lineare}) such that
\begin{equation}
\label{uffa1}
\| \phi \|_{H^1_{\e}} \leq C \e^{- \max \{2, \frac{k}{2}\}} \| f
\|_{L^2 (\mathcal{M}_{\e , \gamma} )}.
\end{equation}
\end{proposition}

The proof will be given in Section \ref{linearas}.

\medskip

Using this, we can get the existence of solution to the linear problem in the whole domain
$\mathcal{M}_\e$.

\begin{proposition}\label{glinear}
There exist a sequence $\e_l \to 0$ and a positive
constant $C >0$, such that, for any $f \in L^2 (\mathcal{M}_{\e_l} )$,
there exists a solution $\phi \in H^1 (\mathcal{M}_{\e_l} )$ to the
equation
$$
L_{\e_l }\phi = f \quad \hbox{in } \mathcal{M}_{\e_l }.
$$
Furthermore,
\begin{equation}
\label{gigio}
\| \phi \|_{H^1_{g^{\e}} (\mathcal{M}_{\e_l} )} \leq C\, \e_l^{-\max \{
2, \frac{k}{2} \}} \| f \|_{L^2 (\mathcal{M}_{\e_l})}.
\end{equation}
\end{proposition}

\begin{proof}
By contradiction, assume that for all $\e \to 0$ there exists a
solution $(\phi_\e , \lambda_\e )$, $\phi_\e \not= 0$, to
\begin{equation}
\label{uno}
L_\e (\phi_\e )= \lambda_\e \phi_\e \quad {\mbox {in}} \quad
\mathcal{M}_\e ,
\end{equation}
with
\begin{equation}\label{due}
|\lambda_\e | \e^{-\max \{ 2, \frac{k}{2}\} } \to 0 , \quad {\mbox
{as}} \quad \e \to 0.
\end{equation}
Let $\eta_\e $ be a smooth cut off function (like the one defined in
(\ref{magaly})) so that
$$
\eta_\e = 1 \quad\hbox{ if dist$_{g^\e}(z, K_\e )
<{\e^{-\gamma} \over 2}$\qquad and} \quad \eta_\e =0 \quad \hbox{if dist$_{g^\e}(z, K_\e ) >
\e^{-\gamma} $.}
$$
In particular one has that
 $|\nabla_{K_\e} \eta_\e| \leq c
\e^{\gamma}$ and  $|\Delta_{K_\e} \eta_\e |\leq c \e^{2\gamma}$  in the
whole domain.

Define $\tilde \phi_\e = \phi_\e \eta_\e$. Then $\tilde \phi_\e$
solves
\begin{equation}\label{tre}
   L_\e (\tilde \phi_\e ) = \lambda_\e \tilde \phi_\e -\nabla_{K_\e} \eta_\e \nabla_{K_\e} a \phi_\e - \Delta_{K_\e} \eta_\e \phi_\e   \text{ in } \mathcal{M}_{\ve,\gamma},
\end{equation}
where $\mathcal{M}_{\e , \gamma}$ is the set defined in
(\ref{omegaepsilongamma}). We now apply Proposition (\ref{teouffa}), that
guarantees the existence of a sequence $\e_l \to 0$ and a constant
$c$ such that
\begin{equation}
\label{quattro} \| \tilde \phi_{\e_l} \|_{H^1_{\e_l}} \leq c
\e_l^{-\max \{2,\frac{k}{2}\}} \left[ \lambda_{\e_l} \| \tilde \phi_{\e_l}
\|_{L^2} + \| \nabla_{K_\e} \eta_{\e_l} \nabla_{K_\e}  \phi_{\e_l} \|_{L^2} + \|
\Delta_{K_\e} \eta_{\e_l}  \phi_{\e_l} \|_{L^2} \right].
\end{equation}
Observe now that, in the region where $\nabla_{K_\e} \eta_{\e_l} \not= 0$
and $\Delta_{K_\e} \eta_{\e_l} \not= 0$, the function $\tilde{U}_{\e_l}$ can be
uniformly bounded $|\tilde{U}_\e (y) |\leq c \e $, with a positive constant
$c$, fact that follows directly from (\ref{Vdef}) and (\ref{boh1}).
Furthermore, since we are assuming (\ref{due}), we see that in the
region we are considering, namely where $\nabla_{K_\e} \eta_{\e_l} \not= 0$
and $\Delta_{K_\e} \eta_{\e_l} \not= 0$, the function $\phi_{\e_l}$
satisfies the equation $$- \Delta_{K_\e} \phi_{\e_l} + \e_l^2 a_{\e_l} (y)
\phi_{\e_l} = 0$$ for a certain smooth function $a_{\e_l}$, which is
uniformly positive and bounded as $\e_l \to 0$. Elliptic estimates
give that, in this region, $|\phi_{\e_l} |\leq c
e^{-\e_l^{\gamma'}}$, and $|\nabla_{K_\e} \phi_{\e_l} | \leq c
e^{-\e_l^{\gamma'}}$ for some $\gamma'>0$ and $c>0$. Inserting this
information in (\ref{quattro}), it is easy to see that
$$
\| \tilde \phi_{\e_l} \|_{H^1_{\e_l}} \leq c \e_l^{-\max \{2,\frac{k}{2}\}}
\lambda_{\e_l} \| \tilde \phi_{\e_l} \|_{H^1_{\e_l}}  (1+ o(1))
$$
where $o(1) \to 0$ as $\e_l \to 0$. Taking into account (\ref{due})
the above inequality gives a contradiction with the fact that, for
all $\e$, the function $\phi_\e$ is not identically zero. This
concludes the proof.
\end{proof}

\begin{proof}[Proof of the main Theorem]
By Proposition \ref{glinear}, $\phi \in H^1_{g^\e} (\mathcal{M}_\e )$ is a
solution to (\ref{nonlinearproblem}) if and only if
$$
\phi = L_\e^{-1} \left( S_\e (\tilde{U}_\e ) + N_\e (\phi )
\right).
$$
Notice that
\begin{equation}
\label{nonno1} \| N_\e (\phi ) \|_{L^2 (\mathcal{M}_\e )} \leq C
\begin{cases}
    \| \phi \|_{H^1 (\mathcal{M}_\e )}^p & \text{ for } p\leq 2, \\
    \| \phi \|_{H^1 (\mathcal{M}_\e )}^2 & \text{ for } p>2
  \end{cases} \quad\qquad  \| \phi \|_{H^1_{g^\e} (\mathcal{M}_\e )} \leq 1
\end{equation}
and
\begin{eqnarray}
\label{nonno2}
&&\| N_\e (\phi_1 ) - N_\e (\phi_2 ) \|_{L^2
(\mathcal{M}_\e
)}  \nonumber\\
&&\leq C \,\begin{cases}
   \left(  \| \phi_1  \|_{H^1_{g^\e} (\mathcal{M}_\e )}^{p-1} + \| \phi_2  \|_{H^1_{g^\e} (\mathcal{M}_\e )}^{p-1} \right)
   \| \phi_1 -\phi_2 \|_{H^1_{g^\e} (\mathcal{M}_\e )} & \text{ for } p\leq 2, \\[3mm]
    \left(  \| \phi_1  \|_{H^1_{g^\e} (\mathcal{M}_\e )} + \| \phi_2  \|_{H^1_{g^\e} (\mathcal{M}_\e )} \right) \| \phi_1 -\phi_2 \|_{H^1_{g^\e} (\mathcal{M}_\e )}  & \text{ for } p>2
  \end{cases},
\end{eqnarray}
for any $\phi_1$, $\phi_2$ in $H^1_{g^\e} (\mathcal{M}_\e )$ with $\| \phi_1
\|_{H^1_{g^\e} (\mathcal{M}_\e )}$, $ \| \phi_2 \|_{H^1_{g^\e} (\mathcal{M}_\e )} \leq 1.$

Defining $T_\e :H^1_{g^\e} (\mathcal{M}_\e ) \to H^1_{g^\e} (\mathcal{M}_\e )$ as
$$
T_\e (\phi ) = L_\e^{-1} \left( S_\e (\tilde{U}_\e ) + N_\e (\phi
) \right)
$$
we will show that $T_\e$ is a contraction in some small ball in $H^1_{g^\e} (\mathcal{M}_\e )
$. A direct consequence of (\ref{bf4}), we have
$$
\|S_\e (\tilde{U}_\e )\|_{L^2(\mathcal{M}_\e)}\leq C\e^{I+\frac{1}{2}}.
$$
Using this inequality and by
(\ref{nonno1}),
(\ref{nonno2}) and (\ref{gigio}), we obtain
$$
\| T_\e (\phi ) \|_{H^1_{g^\e} (\mathcal{M}_\e )} \leq C \e^{- \max \{ 2 , \frac{k}{2} \}}
\begin{cases}
   \left( \e^{I+\frac{1}{2} }+  \| \phi  \|_{H^1_{g^\e} (\mathcal{M}_\e )}^p \right) & \text{ for } p\leq 2, \\
    \left( \e^{I+\frac{1}{2}} + \| \phi  \|_{H^1_{g^\e} (\mathcal{M}_\e )}^2 \right)   & \text{ for } p>2.
  \end{cases}
$$
Now we choose integers $d$ and $I$ so that
$$
d> \begin{cases}
   {\max \{ 2 ,\frac{ k}{2} \} \over p-1}& \text{ for } p\leq 2, \\
    \max \{ 2 ,\frac{ k}{2} \}  & \text{ for } p>2
  \end{cases} \quad I > d-\frac{1}{2} + \max \{ 2 , \frac{k}{2} \}.
  $$
  Thus one easily gets that $T_\e$ has a unique fixed point in set
  $${\mathcal B} = \{ \phi \in H^1_{g^\e} (\mathcal{M}_\e )\, : \, \| \phi \|_{H^1_{g^\e} (\mathcal{M}_\e )} \leq \e^d \},
  $$
  as a direct application of the contraction mapping Theorem. This concludes the proof.
\end{proof}

\section{A linear problem: proof of Proposition \ref{teouffa}}\label{linearas}
\smallskip

The quadratic functional associated to  problem (\ref{lineare})
given by
\begin{equation}\label{functional1}
E(\phi ) = {1\over 2} \int_{\mathcal{M}_{\e,\gamma}}
(|\nabla_{g^\e} \phi |^2 + \e h \phi^2 - (p\pm \e) \tilde{U}_\e^{p\pm \e-1} \phi^2 )
\end{equation}
for functions $\phi \in H^1_{g^\e}(\mathcal{M}_{\e,\gamma}) $.

Let
$(y,x)\in \mathbb{R}^{k+N}$ be the local coordinates along $K_\e$ introduced in (\ref{eq:fc}). With abuse of notation we will denote
\begin{equation}
\label{bb}
  \phi (\mathfrak F ( y,  x))= \phi(z, x ),\quad \mbox{with}\ \ y=\sqrt{\e}z.
\end{equation}
Since the original variable $(z, x)\in \mathbb{R}^{k+N}$
are only local coordinates along $K_\e$ we let the variable $(z, x)$
vary in the set $\mathcal{C}_\e$ defined by
\begin{equation}\label{ddomain}
\mathcal{C}_\e = \{ (z,x) \ / \ \sqrt{\e} z\in \ K,\quad    | x|  <  \e^{-\gamma} \}.
\end{equation}
We write $\mathcal{C}_\e=\frac{1}{\sqrt{\e}} K\times\hat {\mathcal{C}_\e}$ where
\begin{equation}\label{dddomain}
\hat {\mathcal{C}_\e} = \{ x \ / \  |x|  <  \e^{-\gamma} \}.
\end{equation}
Observe  that $\hat {\mathcal{C}_\e}$ approaches, as $\e \to 0$, the whole
space $\mathbb{R}^N$.

In these new local coordinates,  the energy density associated to
the energy $E$ in (\ref{functional1}) is given by
\begin{equation} \label{e}
\frac12 \left[|\nabla_{g^\e} \phi|^2+\e h \phi^2 - (p\pm\e) \tilde{U}_\e^{p\pm \e-1}
\phi^2  \right] \sqrt{\det(g^\e)},
\end{equation}
where $\nabla_{g^\e}$ denotes the gradient in the new variables and
where $g^\e$ is the metric in the coordinates
$(z, x)$. Using the expansions contained in the proof of Lemma \ref{scaledlaplacian}, we have that, if $(z,x)$ vary in $\mathcal{C}_\e$,
then, the energy functional (\ref{functional1}) in the new variables
(\ref{bb}) is given by
\begin{eqnarray}\label{energydensity}
E  ( \phi) & = &\int_{K_\e \times \hat{\mathcal{C}_\e}} \left(\frac12 (
|\nabla_x \phi|^2 +\e h \phi^2 - (p\pm \e) \tilde{U}_\e^{p\pm \e-1} \phi^2 )  \right)
\sqrt{\det(g^\e)} \, dz \, dx\nonumber\\
&& -\frac{\e}{6}\int_{K_\e \times \hat{ \mathcal{C}_\e}} R_{islj} x_lx_s \,\partial_{i}\phi\partial_{j}\phi\,\sqrt{\det(g^\e)} \, dz \, dx\\
 &&
 +\frac12 \int_{K_\e \times \hat{\mathcal{C}_\e}} |\nabla_{K_\e}\phi|^2\,\sqrt{\det(g^\e)} \, dz \, dx+
 \int_{K_\e \times \hat{\mathcal{C}_\e}} B(\phi,\phi)\,\sqrt{\det(g^\e)} \, dz \, dx, \nonumber
\end{eqnarray}
where we denoted by $B(\phi,\phi)$ a quadratic term in $\phi$ that can be
expressed in the following form
\begin{eqnarray}\label{defBB}
 B(\phi,\phi)= O \left(   \e^{\frac{3}{2}} | x  |^3   \right)\partial_{i} \phi \partial_{j} \phi  +{\e }\,|\nabla_{K_\e}\phi|^2\,
  O(\sqrt{\e} |x|) +\partial_j \phi \partial_{\bar a}\phi \left(\mathcal{O}(\sqrt{\e} |x
|)\right)
\end{eqnarray}
and we used the Einstein convention over repeated indices.
Furthermore we use the notation $\partial_a =
\partial_{y_a} $ and $\partial_{\bar a} = \partial_{z_a}$. A detailed proof of expansion (\ref{energydensity}) can be found in \cite{demamu}.

\medskip
\noindent
Given a function $\phi \in H^1_{g^{\e}}(\mathcal{M}_{\e,\gamma})$, we decompose it
as
\begin{equation}\label{decomp}
\phi= \left[{\delta \over \, \mu_\e} \widetilde{{\mathcal T}}_{\, \mu_\e ,
\, \Phi_\e} ( Z_0) + \sum_{j=1}^{N} {d^j \over \, \mu_\e}
\widetilde{{\mathcal T}}_{\, \mu_\e , \, \Phi_\e} ( Z_j)
 + {e \over \, \mu_\e}   \widetilde{{\mathcal T}}_{\, \mu_\e , \, \Phi_\e} (Z) \right] \bar \chi_\e   + \phi^\bot
\end{equation}
where the expression $\widetilde{{\mathcal T}}_{\, \mu_\e , \, \Phi_\e} (
v)$ is defined in (\ref{defTTa}), the functions $Z_0$ and $Z_j$ are
already defined in (\ref{lezetas}) and where $Z$ is the eigenfunction, with $\int_{\mathbb{R}^N} Z^2 =1$,
corresponding to the unique positive eigenvalue $\lambda_0$ in $L^2
(\mathbb{R}^N)$ of the problem
\begin{equation}\label{lambda0}
\Delta_{\mathbb{R}^N} \phi + p w_0^{p-1} \phi = \lambda_0 \phi
\quad {\mbox {in}} \quad \mathbb{R}^N.
\end{equation}
It is worth mentioning that $Z (\xi )$ is even and it has
exponential decay of order $O(e^{-\sqrt{\lambda_0} |\xi|} )$ at
infinity. The function $\bar \chi_\e$ is a smooth cut off function
defined by
\begin{equation}
\label{chibar} \bar \chi_\e (x) = \hat \chi_\e \left(  \left|\left({  x-
\sqrt{\e}\Phi_\e \over \,  \mu_\e}\right) \right| \right),
\end{equation}
with $\hat \chi(r) = 1$ for $r \in (0,{3\over 2} \e^{-\gamma} )$,
and $\chi(r)=0$ for $r>2\e^{-\gamma}$. Finally, in (\ref{decomp}) we
have that $\delta = \delta (\sqrt{\e} z)$, $d^j = d^j (\sqrt{\e} z)$ and $e= e(\sqrt{\e}
z)$ are function defined in $K$ such that $\forall z\in K_\e$
\begin{equation}
\label{orth1}
\int_{\hat{\mathcal{C}}_\e} \phi^\bot {\widetilde{{\mathcal T}}}_{\, \mu_\e , \, \Phi_\e} (Z_0) \bar \chi_\e d x =
\int_{\hat{\mathcal{C}}_\e} \phi^\bot {\widetilde{{\mathcal T}}}_{\, \mu_\e ,
\, \Phi_\e} (Z_j) \bar \chi_\e = \int_{\hat{\mathcal{C}}_\e}
\phi^\bot {\widetilde{{\mathcal T}}}_{\, \mu_\e , \, \Phi_\e} (Z) \bar
\chi_\e=0.
\end{equation}
We will denote by $(H_\e^1)^\bot$ the subspace of the functions in
$H_\e^1$ that satisfy the orthogonality conditions (\ref{orth1}).

A direct computation shows that
$$
\delta (\sqrt{\e} z) = {\int \phi  {\widetilde{{\mathcal T}}}_{\, \mu_\e ,
\Phi_\e} (Z_0) \over \,   \mu_\e \int Z_0^2} (1+ O(\e )) + O(\e )
(\sum_j d^j (\sqrt{\e} z) + e (\sqrt{\e} z)), \quad
$$
$$
d^j (\sqrt{\e} z) = {\int \phi  {\widetilde{{\mathcal T}}}_{\, \mu_\e , \, \Phi_\e}
(Z_j) \over \,   \mu_\e \int Z_j^2} (1+ O(\e ))  + O(\e ) (\delta
(\sqrt{\e} z) + \sum_{i\not= j}  d^i (\sqrt{\e} z) + e (\sqrt{\e} z)),
$$
and
$$
e(\sqrt{\e} z) = {\int \phi  {\widetilde{{\mathcal T}}}_{\, \mu_\e , \, \Phi_\e} (Z)
\over \,   \mu_\e \int Z^2} (1+ O(\e))  + O(\e) (\delta (\sqrt{\e} z)
+\sum_j d^j (\sqrt{\e} z) ).
$$
Observe that, since $\phi \in H_{g^\e}^1$, one easily get that the
functions $\delta $, $d^j$ and $e$ belong to the Hilbert space
\begin{equation}\label{H1K}
{\mathcal H}^1 (K) = \{ \zeta \in {\mathcal L}^2 (K) \, : \,
\partial_a \zeta \in {\mathcal L}^2 (K),\quad a=1,\cdots,k \}.
\end{equation}

Thanks to the above decomposition (\ref{decomp}), we have the validity
of the following expansion for $E(\phi)$.

\bigskip

Observe that in the region we are considering the function $\tilde{U}_\ve$
is nothing but $\tilde{U}_\ve= {\widetilde{{\mathcal T}}}_{\, \mu_\e , \, \Phi_\e}
(v_{I+1 , \e})$, where $v_{I+1, \e}$ is the function whose existence
and properties are proven in Lemma \ref{Construction}. For the
argument in this part of our proof it is enough to take $I=3$, and
for simplicity of notation we will denote by $\hat w$ the function
$v_{I+1 , \e}$ with $I=3$. Referring to (\ref{bf4}) we have
\begin{equation}\label{hatw}
\hat w (z, \xi) = w_0 (\xi) + \sum_{i=1}^4 w_{i,\e} (z,
\xi)
\end{equation}
where $w_0$ is defined by (\ref{roma1}) and
\begin{equation}\label{hatw1}
\| D^2_\xi w_{i+1,\e} \|_{\ve , N-2 , \sigma } + \|
D_\xi w_{i+1,\e} \|_{\ve , N-3 , \sigma } +\|w_{i+1,\e} \|_{\ve, N-4
, \sigma}\le C \ve^{i+\frac{1}{2}}
\end{equation}
and, for any integer $\ell$
$$
\|\nabla^{(\ell)}_{y} w_{i+1,\e}(y,\cdot)\|_{\ve,N-2,\sigma} \leq \beta C_l \ve^{i+\frac{1}{2}}
\qquad \quad y = \sqrt{\ve} z \in K
$$
for any $i=0, 1, 2, 3$.

\begin{theorem} \label{teo4.1}
Let $\gamma= 1-\sigma$, for some $\sigma >0$ and small.
Assume we write $\phi \in H^1_\e$ as in (\ref{decomp}) and let $d=
(d^1 , \ldots , d^{N})$. Then, there exists $\e_0>0$ such that,
for all  $0<\e <\e_0$, the following expansion holds true
\begin{equation}\label{EEE}
E(\phi ) = E(\phi^\bot ) + \e^{-\frac{k}{2}} \left[ P_\e (\delta )
+Q_\e (d ) + R_\e (e)  \right]  +{\mathcal G} (\phi^\bot , \delta ,
d , e).
\end{equation}
In (\ref{EEE})
\begin{equation}\label{q0e}
P_\e (\delta ) = P (\delta ) + P_1 (\delta)
\end{equation}
with
\begin{eqnarray*}
 P (\delta )  &=&{A_\e \over 2} \int_K \e
|\nabla_K ( \delta (1 + o(\e)  \beta_1^\e (y) ) )|^2\nonumber\\
&& +\e
\int_K  \delta^2\left(2c_{3,N}h  -c_{2,N} \sum\limits_j \bigg[ \sum_s
\frac23 R_{mssj}+\sum\limits_{m,a,b}
 \big( \,{\tilde g}^{ab} R_{mabj}- \G_{am}^b
\G_{bj}^a \big)\bigg]\right)\nonumber\\
&& \mp\e
c_{4,N}\int_K  \frac{\delta^2}{\mu_0}
\end{eqnarray*}
where
$A_\e$ a real number such that $\lim\limits_{\e \to 0 } A_\e = c_{1,N}:=
\int_{\mathbb{R}^N} Z_0^2 $,  and $c_{2,N},c_{3,N},c_{4,N}$ are given in (\ref{c2})-(\ref{c3}),
$\beta_1^\e$ is an explicit smooth function defined on $K$ which is
uniformly bounded as $\e \to 0$; furthermore,  $P_1 (\delta ) $ is a
small compact perturbation in ${\mathcal H}_g^1 (K) $ whose shape is a
sum of quadratic functional in $\delta$ of the form
$$
 \e^{2} \int_K b(y) |\delta |^2
$$
where $b(y)$ denotes a generic explicit function, smooth and
uniformly bounded, as $\e \to 0$, in $K$. In (\ref{EEE}),
\begin{equation}
\label{qje}
Q_\e (d ) = Q(d ) + Q_1 (d)
\end{equation}
with
\begin{equation} \label{Qj}
Q (d)=  {\e  \over 2} C_\e \left( \int_K |\nabla_K ( d (1+ o(\e^2
) \beta_2^\e (y)  ) )|^2 + \int_K (\tilde g^{ab}R_{mabl} - \Gamma_a^c (E_m)
\Gamma_c^a (E_l) ) d^m d^l \right)
\end{equation}
where $C_\e$ is a real number such that $\lim_{\e \to 0} C_\e = C:=
\int_{\mathbb{R}^N_+} Z_1^2$, $\beta_2^\e$ is an explicit smooth function
defined on $K$ which is uniformly bounded as $\e \to 0$ and the
terms $R_{maal} $ and $\Gamma_a^c (E_m)$ are smooth functions on $K.$ Furthermore,
$Q_1 (d)$ is a small compact perturbation  in ${\mathcal H}^1 (K) $
whose shape is a sum of quadratic functional in $d$ of the form
$$
\e^{3} \int_K b(y) d^i d^j
$$
where again   $b(y)$ is a generic explicit function, smooth and
uniformly bounded, as $\e \to 0$, in $K$. In (\ref{EEE}),
\begin{equation}
\label{qe0}
R_\e (e) = R(e) + R_1 (e)
\end{equation}
\begin{equation}
\label{Q} R(e) = \e^{-\frac{k}{2}} \left[ {D_\e \over 2} \left( \ve^2 \int_K
|\nabla_K ( e (1+ e^{-{\lambda_0 \over 2} \e^{-\gamma}} \beta_3^\e
(y)  ) )|^2 -\lambda_0  \int_K  e^2 \right)\right]
\end{equation}
with $D_\e$ a real number so that  $\lim_{\e \to 0 } D_\e = D:=
\int_{\mathbb{R}^N} Z^2 $, $\beta_3^\e$ an explicit smooth function in
$K$, which is uniformly bounded as $\e \to 0$,  and $\lambda_0$ the positive number defined in (\ref{lambda0}). Furthermore, $R_1$
is a small compact perturbation  in ${\mathcal H}^1 (K) $ whose shape is a sum of quadratic functional in $e$ of the form
$$
\e^{3} \int_K b(y) e^2
$$
where again   $b(y)$ is a generic explicit function, smooth and
uniformly bounded, as $\e \to 0$, in $K$. Finally in (\ref{EEE})
$$
{\mathcal G} : (H^1_{g^\e} )^\bot \times ({\mathcal H}^1 (K) )^{N+1} \to
\mathbb{R}
$$
is a continuous and differentiable functional with respect to the
natural topologies, homogeneous of degree $2$
$$
{\mathcal G} (t \phi^\bot , t \delta , t d , t e ) = t^2 {\mathcal
M} ( \phi^\bot , \delta ,  d ,  e ) \quad \forall t.
$$
The derivative of ${\mathcal G}$ with respect to each one of its
variable is given by a small multiple of a linear operator in
$(\phi^\bot , \delta , d , e)$ and it satisfies
$$
\| D_{(\phi^\bot , \delta , d)} {\mathcal G}(\phi_1^\bot , \delta_1
, d_1 , e_1) - D_{(\phi^\bot , \delta , d)} {\mathcal G}(\phi_2^\bot ,
\delta_2 , d_2 , e_2) \| \leq C \e^{\gamma (N-3)} \times
$$
\begin{equation}
\label{lips} \left[ \| \phi_1^\bot - \phi_2^\bot \| + \e^{-\frac{k}{2}} \|
\delta_1 - \delta_2 \|_{{\mathcal H}^1 (K)} + \e^{-\frac{k}{2}} \| d_1 - d_2
\|_{({\mathcal H}^1 (K))^{N-1}} + \e^{-\frac{k}{2}} \| e_1 - e_2 \|_{{\mathcal
H}^1 (K)}\right].
\end{equation}
 Furthermore, there exists a constant $C>0$ such that
\begin{equation}
\label{estMM}
\left| {\mathcal G } ( \phi^\bot , \delta ,  d ,  e )
\right| \leq C \ve^{2} \left[ \| \phi^\bot \|^2 + \e^{-\frac{k}{2}} \left( \|
\delta \|_{{\mathcal H}^1 (K)}^2  + \| d \|_{{\mathcal H}^1 (K)}^2 +
\| e\|_{{\mathcal H}^1 (K)}^2   \right) \right].
\end{equation}
\end{theorem}

\medskip
\begin{proof}

\medskip
\noindent
STEP 1. We claim that there exists $\e_0>0 $ such that for all
$0<\e<\e_0$, we have
\begin{equation} \label{leQ0}
E\left({\delta \over \, \mu_\e} {\widetilde{\mathcal T}}_{\, \mu_\e ,\, \Phi_\e
} ( Z_0 ) \bar \chi_\e\right)= \e^{-\frac{k}{2}} P_\e (\delta ),
\end{equation}
\begin{equation} \label{leQj}
E\left({ d^j \over \, \mu_\e} {\widetilde{\mathcal T}}_{\, \mu_\e ,\, \Phi_\e }
( Z_ j) \bar \chi_\e\right)= \e^{-\frac{k}{2}} Q_\e (d^j),
\end{equation}
\begin{equation} \label{leQ}
E\left({e\over \, \mu_\e} {\widetilde{\mathcal T}}_{\, \mu_\e ,\, \Phi_\e } (Z)
\bar \chi_\e\right)= \e^{-\frac{k}{2}} R_\e (e).
\end{equation}
Define
\begin{eqnarray}
\label{defFunctionalF}
 F (u ) :&=&\int_{K_\e \times \hat{\mathcal{C}}_\e}
\left(\frac12  |\nabla_x u|^2 +\frac12 \e hu^2 - \frac{1}{p+1\pm\e} u^{p+1\pm\e}
\right) \sqrt{\det(g^\e)} \, dz \, dx \nonumber\\
&&-\frac{\e}{6}\int_{K_\e \times \hat{\mathcal{C}}_\e} R_{islj}\  x_l\ x_s \,\partial_{i} u \partial_{j} u\sqrt{\det(g^\e)} \, dz \, dx\\
 &&+\frac12 \int_{K_\e \times \hat{\mathcal{C}}_\e}  |\nabla_{K_\e}u|^2\sqrt{\det(g^\e)} \, dz \, dx
  + \int_{K_\e \times \hat{\mathcal{C}}_\e} B(u,u)\sqrt{\det(g^\e)} \, dz \,
  dx.\nonumber
\end{eqnarray}

\smallskip

\noindent
To prove (\ref{leQ0})), we write for
small $t \not= 0$
\begin{eqnarray}\label{tere1}
&&\left[ DF({\widetilde{\mathcal T}}_{\, \mu_\e + t \delta , \, \Phi_\e }
(\hat w ) \bar \chi_\e) - DF({\widetilde{\mathcal T}}_{\, \mu_\e  ,
\Phi_\e } (\hat w ) \bar \chi_\e )\right] ({\delta \over \mu_\e+t\delta} {\widetilde{\mathcal T}}_{\, \mu_\e+t\delta , \, \Phi_\e} (Z_0 ) \bar
\chi_\e)\nonumber\\
&=& -2t E ({\delta \over \, \mu_\e} {\widetilde{\mathcal T}}_{\, \mu_\e ,
\, \Phi_\e} (Z_0 ) \bar \chi_\e ) (1+ O(t)).
\end{eqnarray}
On the other hand,  for any $\psi$
\begin{equation}
\label{tere2} \left[ DF({\widetilde{\mathcal T}}_{\, \mu_\e + t \delta ,
\Phi_\e } (\hat w )\bar \chi_\e) - DF({\widetilde{\mathcal T}}_{\, \mu_\e  ,
\, \Phi_\e } (\hat w ) \bar \chi_\e )\right] ( \psi) =
\mathfrak{a} (t) - \mathfrak{a} (0)  + \mathfrak{b}(t) +
\mathfrak{c}(t)
\end{equation}
where
\begin{eqnarray*}
\mathfrak{a}(t) &= &\int_{K_\e \times \hat{\mathcal{C}}_\e}\left[ \left( \nabla_x
{\widetilde{\mathcal T}}_{\, \mu_\e + t \delta , \, \Phi_\e } (\hat w ) \bar
\chi_\e \right) \nabla_x \psi + \e h {\widetilde{\mathcal T}}_{\, \mu_\e + t
\delta , \, \Phi_\e } (\hat w ) \bar \chi_\e \psi \right.\\
&&\quad \qquad\left.-
\left({\widetilde{\mathcal T}}_{\, \mu_\e + t \delta , \, \Phi_\e } (\hat w
) \bar \chi_\e \right)^{p\pm \e} \psi\right]\sqrt{\det(g^\e)} \, dz \, dx
\\
&&-\frac{\e}{6} \int_{K_\e \times \hat{\mathcal{C}}_\e} R_{islj} \ x_lx_s \partial_i
\left( {\widetilde{\mathcal T}}_{\, \mu_\e + t \delta , \, \Phi_\e } (\hat w
)\bar \chi_\e \right)
\partial_j \psi\sqrt{\det(g^\e)} \, dz \, dx,
\end{eqnarray*}
\begin{eqnarray*}
\mathfrak{b}(t)&=&\int_{K_\e \times \hat{\mathcal{C}}_\e} \partial_{\bar a} (
{\widetilde{\mathcal T}}_{\, \mu_\e + t \delta , \, \Phi_\e } (\hat w )\bar
\chi_\e )
\partial_{\bar a} \psi \sqrt{\det(g^\e)} \, dz \, dx\\
&&- \int_{K_\e \times \hat{\mathcal{C}}_\e} \partial_{\bar a} ( {\widetilde{\mathcal T}}_{
\mu_\e  , \, \Phi_\e } (\hat w ) \bar \chi_\e ) \partial_{\bar a}
\psi\sqrt{\det(g^\e)} \, dz \, dx
\end{eqnarray*}
and
\begin{eqnarray*}
\mathfrak{c}(t)&=& \int_{K_\e \times \hat{\mathcal{C}}_\e} B({\widetilde{\mathcal T}}_{
\mu_\e + t \delta , \, \Phi_\e } (\hat w ) \bar \chi_\e , \psi) \sqrt{\det(g^\e)} \, dz \, dx\\
&&-
\int_{K_\e \times \hat{\mathcal{C}}_\e} B({\widetilde{\mathcal T}}_{ \mu_\e  , \bar
\Phi_\e } (\hat w ) \bar \chi_\e , \psi)\sqrt{\det(g^\e)} \, dz \, dx.
\end{eqnarray*}
We now compute $\mathfrak{a} (t) - \mathfrak{a} (0) $ with $\psi = {\delta \over
\mu_\e+t\delta} {\widetilde{\mathcal T}}_{\, \mu_\e + t \delta , \, \Phi_\e } (Z_0 )
\bar \chi_\e $.
Performing the change of variables $x =
(\mu_\e + t \delta )  \xi + \, \sqrt{\e} \Phi_\e$ and
using the expansion of ${\delta \over \, \mu_\e+t\delta}=\frac{\delta}{\mu_e}-\frac{\delta^2}{\mu_e^2}t+O(t^2)$, we see   that
\begin{eqnarray*}
&&t^{-1} \left[\mathfrak{a}(t)-\mathfrak{a}(0) \right]  \\
&=&- \e\left\{
\int_{K_\e} \delta^2 \left(c_{3,N}h  -c_{2,N} \sum\limits_j \bigg[ \sum_s
\frac23 R_{mssj}+\sum\limits_{m,a,b}
 \big( \, {\tilde g}^{ab} R_{mabj}- \G_{am}^b
\G_{bj}^a \big)\bigg]\right)  \mp c_{4,N}\int_{K_\e} \frac{\delta^2}{\mu_\e^2}   \right\} \\
&&\times\big(1+ O(t)\big) \bigg(1
+O(\ve^{\gamma (N-4)})\bigg).
\end{eqnarray*}

On the other hand, by the definition of the function $\mathfrak{b}(t)$ above,
 a Taylor expansion gives
$$
\mathfrak{b}(t) = -t\left(  \int_{K_\e \times \hat {\mathcal{C}}_\e} |\nabla_{K_\e} ({\delta \over \, \mu_\e} {\widetilde{\mathcal T}}_{\, \mu_\e ,
\Phi_\e} ( Z_0 ) \bar \chi_\e  )|^2 dzdx  \right) \times  (1+ O(t))  .
$$
Observe now that \begin{eqnarray*}
\partial_{\bar a} \left( {\delta \over \, \mu_\e} {\widetilde{\mathcal T}}_{\, \mu_\e , \, \Phi_\e} (Z_0 ) \bar \chi_\e \right)
&=& (\partial_{\bar a} \delta ) {1 \over \, \mu_\e} {\widetilde{\mathcal
T}}_{\, \mu_\e , \, \Phi_\e} (Z_0 )  \bar \chi_\e  + \delta
\partial_{\bar a} ( {1 \over \, \mu_\e} {\widetilde{\mathcal T}}_{\, \mu_\e
, \, \Phi_\e} (Z_0 ) \bar \chi_\e  )
\\
&=& \sqrt{\e} ( \partial_a \delta ) {1 \over \, \mu_\e} {\widetilde{\mathcal
T}}_{\, \mu_\e , \, \Phi_\e} (Z_0 )\bar \chi_\e   + \sqrt{\e} \delta (
\partial_a \, \mu_\e ) \partial_{\, \mu_\e} ({1 \over
\mu_\e} {\widetilde{\mathcal T}}_{\, \mu_\e , \, \Phi_\e} (Z_0 ) \bar
\chi_\e )\\
 &&+ \sqrt{\e} \delta ( \partial_a \, \Phi_\e ) \partial_{
\Phi_\e}  ({1 \over \, \mu_\e} {\widetilde{\mathcal T}}_{\, \mu_\e , \bar
\Phi_\e} (Z_0 ) \bar \chi_\e ).
\end{eqnarray*}
Since $ \int \left( {1 \over \, \mu_\e} {\widetilde{\mathcal T}}_{\, \mu_\e
, \, \Phi_\e} (Z_0 ) \bar \chi_\e  \right)^2 \, dx=    \int \left( {1 \over \,  \mu_\e} {\widetilde{\mathcal T}}_{\, \mu_\e
, \, \Phi_\e} (Z_0 ) \bar \chi_\e  \right)^2 \, dx=  c_{1,N} (1+ o(\e)
)$, we conclude that
\begin{equation}
\label{III1} \mathfrak{b}(t) = - t \e^{-\frac{k}{2}} \left[ A_\e \e  \int_K
|\nabla_K ( \delta (1+ o(\e  ) \beta_1^\e (y)  )) |^2 \right]
\end{equation}
where $A_\e \in \mathbb{R}$, $\lim_{\e \to 0} A_\e = c_{1,  N}= \int_{\mathbb{R}^N} Z_0^2
$ and $\beta_1^\e$ is an explicit smooth function in $K$, which is
uniformly bounded as $\e \to 0$.   Finally we observe that the last
term $\mathfrak{c}(t)$ defined above is of lower order, and can be
absorbed in the terms already described.

\medskip

\noindent
Proof of (\ref{leQj}).  Let
$d $ be the vector field along $K$ defined by $$d (\sqrt{\e} z) = (d^1 (\sqrt{\ve}
z ) , \ldots , d^{N} (\sqrt{\ve} z) ).$$ For any $t$ small and $t\not= 0$,
we have (see (\ref{defFunctionalF}))
 \begin{eqnarray*}
 &&\left[ DF({\widetilde{\mathcal T}}_{\, \mu_\ve , \, \Phi_\ve + t d} (\hat w) \bar
\chi_\e) - DF ({\widetilde{\mathcal T}}_{\, \mu_\ve , \, \Phi_\ve } (\hat w)
\bar \chi_\e )\right] [\varphi]   \\
&=& t D^2 F ({\widetilde{\mathcal T}}_{\, \mu_\ve , \, \Phi_\ve } (\hat w)
\bar \chi_\e ) \left[\sum_l {   d^l \over \, \mu_\ve} {\widetilde{\mathcal
T}}_{\, \mu_\ve , \, \Phi_\ve} (Z_l ) \bar \chi_\e \right]
[\varphi] \,\big(1+ O(t)\big) \big(1+ O(\ve)\big)
\end{eqnarray*}
for any function $\varphi \in H^1_{g^\ve}$. Choosing
$\varphi = { d^j \over \, \sqrt{\e}\mu_\ve} {\widetilde{\mathcal T}}_{\, \mu_\ve ,
\, \Phi_\ve} (Z_j) \bar \chi_\e$   we write
\begin{eqnarray}
&&\label{leQj2} \left[ DF ({\widetilde{\mathcal T}}_{\, \mu_\ve , \, \Phi_\ve
+ t d} (\hat w) \bar \chi_\e) - DF ({\widetilde{\mathcal T}}_{\, \mu_\ve ,
\, \Phi_\ve } (\hat w) \bar \chi_\e )\right] [{ d^j \over \sqrt{\e}
\mu_\ve} {\widetilde{\mathcal T}}_{\, \mu_\ve , \, \Phi_\ve} (Z_j) \bar
\chi_\e]\nonumber\\
&=& \mathfrak{a}_2(t) - \mathfrak{a}_2(0) + \mathfrak{b}_2(t)
+ \mathfrak{c}_2(t)
\end{eqnarray}
where we have set, for $\psi = { d^j \over \,\sqrt{\e} \mu_\ve} {\widetilde{\mathcal
T}}_{\, \mu_\ve , \, \Phi_\ve} (Z_j) \bar \chi_\e$,
\begin{eqnarray*}
\mathfrak{a}_2(t)& =& \int \left( \nabla_X {\widetilde{\mathcal T}}_{\, \mu_\e
, \, \Phi_\e +t d } (\hat w )  \bar \chi_\e \right) \nabla_X \psi
+ \ve {\widetilde{\mathcal T}}_{\, \mu_\e  , \, \Phi_\e +t d } (\hat w )
\bar \chi_\e \psi - \left({\widetilde{\mathcal T}}_{\, \mu_\e  , \, \Phi_\e
+t d}
(\hat w \bar \chi_\e )\right)^{p\pm\e} \psi\\
& &-\frac{\e}{6} \int R_{islj} \ x_lx_s
\partial_i \left( {\widetilde{\mathcal T}}_{\, \mu_\e  , \, \Phi_\e +t d }
(\hat w ) \bar \chi_\e \right) \partial_j \psi,
\end{eqnarray*}
$$
\mathfrak{ b}_2(t) =\int \partial_{\bar a} ( {\widetilde{\mathcal T}}_{
\mu_\e  , \, \Phi_\e +t d} (\hat w ) \bar \chi_\e) \partial_{\bar
a} ( {\widetilde{\mathcal T}}_{\, \mu_\e  , \, \Phi_\e +t d } (\hat w ) \bar
\chi_\e) - \int \partial_{\bar a} ( {\widetilde{\mathcal T}}_{\, \mu_\e  ,
\, \Phi_\e } (\hat w ) \bar \chi_\e) \partial_{\bar a} ( {\widetilde{\mathcal
T}}_{\, \mu_\e , \, \Phi_\e } (\hat w ) \bar \chi_\e)
$$
and
$$
\mathfrak{c}_2 (t)= \int B({\widetilde{\mathcal T}}_{\bar  \mu_\e  , \bar
\Phi_\e +t d } (\hat w ) \bar \chi_\e , \psi) - \int B({\widetilde{\mathcal
T}}_{\bar  \mu_\e , \, \Phi_\e } (\hat w ) \bar \chi_\e , \psi).
$$
Define
$$
{\mathcal R}_{ml}=\bigg( (\tilde g)^{ab}\,R_{mabl}-\G_a^c(E_m)
\G_c^a(E_l) \bigg).
$$
Performing the change of variables $x=
\, \mu_\e    \xi + \, \sqrt{\e} (\Phi_\e+td)$, we get
\begin{eqnarray*}
&&t^{-1} [\mathfrak{a}_2(t) - \mathfrak{a}_2(0) ]\\
&=&   \e  \left\{
\int { d^j \over \, \mu_\ve} \bigg[
\nabla \hat w \nabla Z_j+ \e  \, \mu_\e^2 h \hat w   - \mu_{\e}^{\mp\frac{N-2}{2}\e} \hat w^{p\pm\e}
Z_j\bigg]  \right.
\\
&&\times   ({R_{mijl} \over 6} + {{\mathcal R}_{lm} \over 2} ) [(
\mu_\ve\xi_m +\, \sqrt{\e}\Phi_{\ve, m} ) d^l +  (\, \mu_\ve \xi_l +\, \sqrt{\e}\Phi_{\ve, l}) d^m ]
\\
&&\left. - \int { d^j \over \, \mu_\ve} {R_{ilsr} \over 6} [(\mu_\ve \xi_s +\, \sqrt{\e}\Phi_{\ve, s}) d^l
+  (\mu_\ve \xi_l +\, \sqrt{\e}\Phi_{\ve, l}) d^s ]
\partial _i \hat w \partial_r Z_j \right\}\\
&&\times  (1+ O(\ve )) (1+ O(t)).
\end{eqnarray*}
Integration by parts in the $\xi$ variables,  using the fact that
$\hat {\mathcal C}_\ve \to \mathbb{R}^N$ as $\ve \to 0$,  $R_{irll}=0$, we
deduce that
\begin{eqnarray}\label{leQj3}
t^{-1} [\mathfrak{a}_2(t) - \mathfrak{a}_2(0) ] &= &\ve   \left\{ -
C \int_{K_\e} ({R_{miij} \over 6} + {{\mathcal R}_{mj} \over 2} ) d^j
d^m   +C \int_{K_\e}  {R_{jrrm} \over 3} d^m  d^j \right\} \times\nonumber\\
&\times& \big(1+ O(\ve )\big) \,\big(1+ O(t)\big)
\\&=& \ve^{-\frac{k}{2}} \ve \left[
- C \int {{\mathcal R}_{mj} \over 2}  d^j d^m + O(\e) Q(d) \right]
\big(1+ O(t)\big)\nonumber
\end{eqnarray}
where here we have set
$$
C= \int_{\R^N} Z_1^2 \qquad \hbox{ and } \quad Q(d):=\int_K \pi(y)
d^i d^j
$$
for some smooth and uniformly bounded (as $\e \to 0$) function
$\pi(y)$.  To estimate the term $\mathfrak{b}_2$ above we argue as
in (\ref{III1}), we get that
\begin{equation}
\label{leQj4} t^{-1} \mathfrak{b}_2(t) = -\ve^{-\frac{k}{2}} \left[ \e C_\e
\int_K |\nabla_K (d^j (1+ \beta_2^\e (y) o(\e) ))|^2  \right]
(1+ O(t)).
\end{equation}
Finally we observe that the last term $\mathfrak{c}_2(t)$ is of
lower order, and can be absorbed in the terms described in
(\ref{leQj3}) and (\ref{leQj4}).

\medskip
\noindent
 Proof of (\ref{leQ}). To get the
expansion in (\ref{leQ}), we compute
\begin{equation}\label{tere3}
E({e\over \, \mu_\e} {\widetilde{\mathcal T}}_{\, \mu_\e ,
\, \Phi_\e} (Z)) = I + II + III
\end{equation}
where
\begin{eqnarray*}
I&=& \int_{K_\e \times \hat{\mathcal{C}}_\e} {e^2 \over  \mu_\e^2}
\left(\frac12 ( |\nabla_x {\widetilde{\mathcal T}}_{ \mu_\e , \Phi_\e} (Z) |^2
+\e h {\widetilde{\mathcal T}}_{\mu_\e ,  \Phi_\e} (Z)^2 - (p\pm\e) \tilde{U}_\e^{p\pm\e-1} {\widetilde{\mathcal
T}}_{ \mu_\e ,  \Phi_\e} (Z)^2 ) \right) \times \\
& &  \sqrt{\det g^\e} dz dx
\\
&&-\frac{\e}{6}\int_{K_\e \times \hat{\mathcal{C}}_\e} {e^2 \over \, \mu_\e^2}
R_{islj}x_sx_l \,\partial_{i} {\widetilde{\mathcal T}}_{
\mu_\e , \, \Phi_\e} (Z) \partial_{j}{\widetilde{\mathcal T}}_{\, \mu_\e , \,
\Phi_\e} (Z) \,\sqrt{\det g^\e} \, dz \, dx,
\\[3mm]
II&=& \frac12 \int_{K_\e \times \hat{\mathcal{C}}_\e} |\nabla_{K} \left(
{e\over \, \mu_\e} {\widetilde{\mathcal T}}_{\, \mu_\e , \, \Phi_\e} (Z) \right)|^2
 \,\sqrt{\det g^\e} \, dz
\, dx\\[3mm]
\mbox{and}\\
III&=& \int_{K_\e \times \hat{\mathcal{C}}_\e} B({e\over \, \mu_\e}
{\widetilde{\mathcal T}}_{\, \mu_\e , \, \Phi_\e} (Z), {e\over \, \mu_\e}
{\widetilde{\mathcal T}}_{\, \mu_\e , \, \Phi_\e} (Z))\,\sqrt{\det(g^\e)} \,
dz \, dx.
\end{eqnarray*}
Using the change of variables $x = \, \mu_\e   \xi + \sqrt{\e}
\Phi_\e$  in $I$, we have
$$
I= \int {1\over 2} {e^2 \over \, \mu_\e^2} \bigg[ |\nabla
Z|^2 - p \hat w^{p-1} Z^2 +\e \, \mu_\e^2 h Z^2 \bigg] \bigg(1+ \e
O(e^{-|\xi|})\bigg).
$$
Then, recalling the definition of $\lambda_0$ in (\ref{lambda0}), we get
\begin{equation} \label{tere4}
I=\e^{-\frac{k}{2}} \left[  -{\lambda_0 \over 2} D \int_K e^2 + \e Q(e) \right]
\end{equation}
where we have set
$$
D= \int_{\mathbb{R}^N} Z^2(\xi)\,d\xi \quad \hbox{and }\qquad Q(e):= \int_K \tau(y)
e^2\,dy,
$$
for some  smooth and uniformly bounded, as $\ve \to 0$,
function $\tau$. On the other hand, using a direct computation and
arguing as in (\ref{III1}), we get
\begin{equation} \label{ttere4}
II= {D_\e \over 2} \int_{K_\e} |\nabla_{K_\e} e  + e^{-\lambda_0
\ve^{-\gamma}} \beta_3^\e (\e z) e |^2 =
 \ve^{-\frac{k}{2}}
 \bigg[ {D_\e \over 2} \ve \int_{K} |\nabla_{K} (e (1+ e^{-\lambda'\ve^{-\gamma}} \beta_3^\e (y) )) |^2 \bigg]
\end{equation}
where $\beta_3^\e$ is an explicit smooth function on $K$, which is
uniformly bounded as $\e \to 0$, while $\lambda'$ is a positive real
number. Finally we observe that the last term $III$ is of lower
order, and can be absorbed in the terms described in (\ref{tere4}) and
(\ref{ttere4}). This concludes the proof of (\ref{leQ}).

\medskip
\noindent
STEP 2.
We write
\begin{eqnarray*}
{\mathcal G} (\phi^\bot , \delta , d, e) &=& E(\phi ) - E(\phi^\bot
) -
 E({\delta \over \, \mu_\e} {\widetilde{\mathcal T}}_{\, \mu_\e ,\, \Phi_\e } ( Z_0 ) \bar \chi_\e)
 - \sum_{j=1}^{N} E({d_j \over \, \mu_\e} {\widetilde{\mathcal T}}_{\, \mu_\e ,\, \Phi_\e } ( Z_j ) \bar
 \chi_\e)\\
&&- E({e \over \, \mu_\e} {\widetilde{\mathcal T}}_{\, \mu_\e ,\, \Phi_\e
} ( Z ) \bar \chi_\e).
\end{eqnarray*}
Thus it is clear that the term ${\mathcal G}$ recollects all the
mixed terms in the expansion of $E(\phi)$. Indeed, if we define
\begin{eqnarray*}
m(f_1,f_2)& =& \int_{K_\e \times \hat{\mathcal{C}}_\e} \left(\nabla_x f_1 \nabla_x
f_2  - (p\pm\e) U_\e^{p\pm\e-1} f_1f_2    \right) \sqrt{\det(g^\e)} \, dz \, dx \\
& & +\e \int_{K_\e \times \hat{\mathcal{C}}_\e}hf_1f_2dzdx
\\
&&-\frac{\e}{6}\int_{K_\e \times \hat{\mathcal{C}}_\e} \,R_{islj}x_sx_l\partial_{i}f_1\partial_{j}f_2\,\sqrt{\det(g^\e)} \, dz \, dx\\
&&+ \int_{K_\e \times \hat{\mathcal{C}}_\e} \pa_{\bar a}f_1\,\pa_{\bar
a}f_2\,\sqrt{\det(g^\e)} \, dz \, dx+ \int_{K_\e \times \hat{\mathcal{C}}_\e}
 B(f_1,f_2)\,\sqrt{\det(g^\e)} \, dz \, dx\\
 &:=& m_1(f_1,f_2)+m_2(f_1,f_2)+m_3(f_1,f_2)+m_4(f_1,f_2)+m_5(f_1,f_2),
\end{eqnarray*}
for $f_1$ and $f_2$ in $H^1_\e$, then
\begin{eqnarray}\label{auxi}
{\mathcal G}(\phi^\bot , \delta , d , e) &=& m (\phi^\bot , {\delta
\over \, \mu_\e} {\widetilde{\mathcal T}}_{\, \mu_\e ,\, \Phi_\e } ( Z_0 )
\bar \chi_\e ) +\sum_j  m (\phi^\bot , {d_j \over \, \mu_\e}
{\widetilde{\mathcal T}}_{\, \mu_\e ,\, \Phi_\e } ( Z_j ) \bar \chi_\e
)\nonumber
\\
&+& m (\phi^\bot , {e \over \, \mu_\e} {\widetilde{\mathcal T}}_{\, \mu_\e
,\, \Phi_\e } ( Z_0 ) \bar \chi_\e ) + \sum_j m ({\delta \over \bar
\mu_\e} {\widetilde{\mathcal T}}_{\, \mu_\e ,\, \Phi_\e } ( Z_0 ) \bar
\chi_\e , {d^j \over \, \mu_\e} {\widetilde{\mathcal T}}_{\, \mu_\e ,\bar
\Phi_\e } ( Z_j ) \bar \chi_\e )\nonumber
\\
&+&\sum_{i\not= j} m ({d^j \over \, \mu_\e} {\widetilde{\mathcal T}}_{\bar
\mu_\e ,\, \Phi_\e } ( Z_j ) \bar \chi_\e , {d_i \over \bar
\mu_\e} {\widetilde{\mathcal T}}_{\, \mu_\e ,\, \Phi_\e } ( Z_i ) \bar
\chi_\e )
\\
&+& m ({\delta \over \, \mu_\e} {\widetilde{\mathcal T}}_{\, \mu_\e ,\bar
\Phi_\e } ( Z_0 ) \bar \chi_\e , {e \over \, \mu_\e} {\widetilde{\mathcal
T}}_{\, \mu_\e ,\, \Phi_\e } ( Z ) \bar \chi_\e )\nonumber\\
& +& \sum_j m ({d^j \over \, \mu_\e} {\widetilde{\mathcal T}}_{\, \mu_\e ,\,
\Phi_\e } ( Z_j ) \bar \chi_\e , {e \over \, \mu_\e} {\widetilde{\mathcal
T}}_{\, \mu_\e ,\, \Phi_\e } ( Z ) \bar \chi_\e ).\nonumber
\end{eqnarray}
One can see clearly that ${\mathcal G}$ is homogeneous of degree $2$
and that its first derivatives with respect to its variables is a
linear operator  in $(\phi^\bot , \delta, d , e)$. We will then show
the validity of estimate (\ref{estMM}). In a very similar way one
shows the validity of (\ref{lips}). To prove (\ref{estMM}), we should
treat each one of the above terms. Since the computations are very
similar, we will limit ourselves to treat the term
$$
m:= m ({\delta
\over \, \mu_\e} {\widetilde{\mathcal T}}_{\, \mu_\e ,\, \Phi_\e } ( Z_0 )
\bar \chi_\e , {d^j \over \, \mu_\e} {\widetilde{\mathcal T}}_{\, \mu_\e
,\, \Phi_\e } ( Z_j ) \bar \chi_\e ).
$$
This term can be written as
\begin{equation}
m= \sum_{i=1}^5 m_i(f_1,f_2)
\end{equation}
with $f_1= {\delta \over \, \mu_\e} {\widetilde{\mathcal T}}_{\, \mu_\e ,\bar
\Phi_\e } ( Z_0 ) \bar \chi_\e$ and $f_2=  {d^j \over \, \mu_\e}
{\widetilde{\mathcal T}}_{\, \mu_\e ,\, \Phi_\e } ( Z_j ) \bar \chi_\e$.
Using the fact that the function $Z_0$ solves
$$\Delta Z_0+ p w_0^{p-1} Z_0 = 0 \quad \hbox{in}\quad \R^N,
$$
with  $\int_{\mathbb{R}^N}
\partial_{\xi_N} Z_0 Z_j = 0$ and integrating by parts in the $x$
variable (recalling the expansion of $\sqrt{{\mbox {det}} g^\e}$),
one gets
\begin{eqnarray*}
m_1 &=& \left\{ \int {\delta d^j \over \, \mu_\e^2} [-\Delta
{\widetilde{\mathcal T}}_{\, \mu_\e , \, \Phi_\e } (Z_0) - (p\pm\e) U_\e^{p\pm\e-1}
{\widetilde{\mathcal T}}_{\, \mu_\e , \, \Phi_\e} (Z_0 ) ] \bar \chi_\e^2
{\widetilde{\mathcal T}}_{\, \mu_\e , \, \Phi_\e} (Z_j ) \sqrt{{\mbox {det}}
g^\e} \right.
\\
&+& \left. \int {\delta d^j \over \, \mu_\e^2} \partial_{\xi_N}
\left( {\widetilde{\mathcal T}}_{\, \mu_\e , \, \Phi_\e} (Z_0) \bar \chi_\e
\right) {\widetilde{\mathcal T}}_{\, \mu_\e , \, \Phi_\e} (Z_j )
\,\frac{1}{\, \mu_\e}\,(\e \,tr(H)+O(\e^2))\,\bar \chi_\e \right\}
(1+ o(1) )
\end{eqnarray*}
where $o(1) \to 0$ as $\e \to 0$. Thus, a H\"older inequality yields
$$
| m_1 | \leq C \e^{-\frac{k}{2}}  \e^{\gamma (N-2)} \| \delta \|_{{\mathcal
L}^2 (K)} \| d^j \|_{{\mathcal L}^2 (K)}.
$$
On the other hand, using the orthogonality condition $\int_{\R^N}
Z_0 Z_j =0$, we get
$$
|m_2 | \leq C \e \e^{-\frac{k}{2}} (\int_{|\xi|>\e^{-\gamma}} Z_0 Z_j )   \|
\delta \|_{{\mathcal L}^2 (K)} \| d^j \|_{{\mathcal L}^2 (K)} \leq C
\e^{-k} \e^{1+\gamma (N-3) } \| \delta \|_{{\mathcal L}^2 (K)} \|
d^j \|_{{\mathcal L}^2 (K)}.
$$
Now, since $\int_{\R^N} \xi_N \partial_i Z_0 \partial_l Z_j = 0$,
for any $i, j, l=1, \ldots , N-1$, one gets
\begin{eqnarray*}
|m_3| &\leq &C \e \e^{-\frac{k}{2}} \left( \int_{|\xi|>\e^{-\gamma}} \xi_N
\partial_i Z_0 \partial_l Z_j \right) \| \delta \|_{{\mathcal L}^2
(K)} \| d^j \|_{{\mathcal L}^2 (K)} \\
&\leq& C \e^{-\frac{k}{2}} \e^{1+\gamma (N-2)} \| \delta \|_{{\mathcal L}^2
(K)} \| d^j \|_{{\mathcal L}^2 (K)}.
\end{eqnarray*}
A direct computation on the term $m_4$ gives
\begin{eqnarray*}
|m_4 | &\leq& C \e^{-\frac{k}{2}} \left\{ \e^2  (\int_{|\xi|>\e^{-\gamma}} Z_0
Z_j ) \| \partial_a \delta \|_{{\mathcal L}^2 (K)} \| \partial_a d^j
\|_{{\mathcal L}^2 (K)} \right.
\\
&+& \e (\int_{|\xi|>\e^{-\gamma}} Z_0 Z_j ) ( \|  \delta
\|_{{\mathcal L}^2 (K)} \| \partial_a d^j \|_{{\mathcal L}^2 (K)}  +
\|  \partial_a \delta \|_{{\mathcal L}^2 (K)} \| d^j \|_{{\mathcal
L}^2 (K)} )
\\
&+& \left.  (\int_{|\xi|>\e^{-\gamma}} Z_0 Z_j ) \|  \delta
\|_{{\mathcal L}^2 (K)} \|  d^j \|_{{\mathcal L}^2 (K)} \right\}
\\
&\leq& C \e^{-\frac{k}{2}} \e^{\gamma (N-3)} [ \| \delta \|_{{\mathcal H}^1
(K)}^2 + \| d^j \|_{{\mathcal H}^1 (K)}^2 ].
\end{eqnarray*}
Since $|m_5| \leq C \sum_{j=1}^4 |m_j|$ we conclude that
$$
|m| \leq C \e^{-\frac{k}{2}} \e^{\gamma (N-3)} [ \| \delta \|_{{\mathcal H}^1
(K)}^2 + \| d^j \|_{{\mathcal H}^1 (K)}^2 ].
$$
Each one of the terms appearing in (\ref{auxi}) can be estimated to
finally get the validity of (\ref{estMM}). This conclude the proof of
Theorem (\ref{teo4.1}).
\end{proof}

\medskip

Now, we are going to prove Proposition \ref{teouffa}.

\begin{proof}
We define the energy functional associated to Problem
(\ref{lineare})
$$
{\mathcal E}: (H^1_{g^{\e}})^\bot \times ( {\mathcal H}^1 (K) )^{N+2} \to
\R
$$
by
\begin{equation}\label{functional2}
{\mathcal E} (\phi^\bot , \delta , d , e) =
E(\phi ) -{\mathcal L}_f (\phi )
\end{equation}
where $E$ is the functional in (\ref{functional1}) and ${\mathcal L}_f
(\phi )$ is the linear operator given by
$$
{\mathcal L}_f (\phi ) = \int_{\mathcal{M}_{\e , \gamma}} f \phi.
$$
Observe that
$$
{\mathcal L}_f (\phi ) = {\mathcal L}_f^1 (\phi^\bot ) + \e^{-\frac{k}{2}}
\left[ {\mathcal L}_f^2 (\delta ) + {\mathcal L}_f^3 (d ) +{\mathcal
L}_f^4 (e) \right]
$$
where ${\mathcal L}_f^1 :H^1_{g^\e} \to \mathbb{R}$, $ {\mathcal L}_f^2 ,
{\mathcal L}_f^4 \, : \, {\mathcal H}^1 (K) \to \mathbb{R}$ and  ${\mathcal
L}_f^3 \, : \, ( {\mathcal H}^1 (K) )^{N} \to \mathbb{R}$ with
$$
{\mathcal L}_f^1 (\phi^\bot ) = \int_{\mathcal{M}_{\e , \gamma}} f
\phi^\bot , \quad \e^{-\frac{k}{2}} {\mathcal L}_f^2 (\delta ) =
\int_{\mathcal{M}_{\e , \gamma}} f {\delta \over \, \mu_\e} {\widetilde{\mathcal
T}}_{\, \mu_\e ,\, \Phi_\e } ( Z_0 ) \bar \chi_\e
$$
$$
\e^{-\frac{k}{2}} {\mathcal L}_f^3 (d ) = \sum_{j=1}^{N}\int_{\mathcal{M}_{\e ,
\gamma}} f {d^j \over \, \mu_\e} {\widetilde{\mathcal T}}_{\, \mu_\e ,\bar
\Phi_\e } ( Z_j ) \bar \chi_\e \quad {\mbox {and}} \quad \e^{-\frac{k}{2}}
{\mathcal L}_f^4 (e ) = \int_{\mathcal{M}_{\e , \gamma}} f {e \over \bar
\mu_\e} {\widetilde{\mathcal T}}_{\, \mu_\e ,\, \Phi_\e } ( Z ) \bar \chi_\e
.
$$
Finding a solution $\phi \in H^1_{g^\e}$ to Problem (\ref{lineare})
reduces to finding a critical point $(\phi^\bot , \delta , d , e) $
for ${\mathcal E}$. This will be done in several steps.

\bigskip

\noindent \textrm{\textbf{Step 1}}. We claim that there exist
$\sigma >0$ and $\e_0$ such that for all $\e \in (0,\e_0)$ and all
$\phi^\bot \in (H^1_{g^\e})^\bot$ then
\begin{equation}\label{punto1}
E(\phi^\bot ) \geq \sigma \| \phi^\bot \|^2_{L^2}.
\end{equation}

\medskip
In fact, using the local change of variables (\ref{bb}),
together with the expansion of energy $E$ in (\ref{energydensity}), we see that, for
sufficiently small $\e>0$
$$
E(\phi^\bot ) \geq {1\over 4} E_0 (\phi^\bot ),
$$
with
$$
E_0 (\phi^\bot ) = \int_{K_\e \times \hat{\mathcal C}_\e}
\left[ |\nabla_x \phi^\bot |^2 - (p\pm\e) U_\e^{p\pm\e-1} (\phi^\bot)^2
\right]\sqrt{\det(g^\e)}dzdx
$$
for any $\phi^\bot = \phi^\bot (\sqrt{\e} z , x),$ with $ z \in K_\e = {1
\over \sqrt{\e}}K$. The set $\hat {\mathcal C}_\e$ is defined in
(\ref{dddomain}) and the function $U_\e$ is given by (\ref{Vdef}). We
recall that $\hat {\mathcal C}_\e \to \R^N$ as $\e \to 0$.

We will establish (\ref{punto1}) showing that
\begin{equation}\label{pas1}
E_0 (\phi^\bot ) \geq \sigma \| \phi^\bot \|^2_{L^2}
\quad \forall \phi^\bot.
\end{equation}
To do so, we first observe that if we scale in the $z$-variable,
defining $\varphi^\bot (y,x) = \phi^\bot ({y\over \sqrt{\e} }, x)$, the
relation (\ref{pas1}) becomes
\begin{equation}\label{pas2}
E_0 (\varphi^\bot ) \geq \sigma \| \varphi^\bot
\|^2_{L^2}.
\end{equation}
Thus we are led to show the validity of (\ref{pas2}). We argue by
contradiction. Assume that   for any $n$, there exist
$\e_n \to 0 $ and $\varphi_n^\bot \in (H_{g^{\e_n}}^1)^\bot $ such that
\begin{equation}\label{shen}
E_0 (\varphi_n^\bot ) \leq \frac1n \| \varphi_n^\bot
\|^2_{L^2}.
\end{equation}
Without loss of generality we can assume that the sequence $(\|
\varphi_n^\bot\|)_n$ is bounded, as $n \to \infty$.
Hence, up to subsequences, we have that
$$
\varphi_n^\bot \rightharpoonup \varphi^\bot \quad \hbox{in}
\quad H^1 (K \times \R^N )\qquad  \hbox{and } \quad \varphi_n^\bot
\to \varphi^\bot \quad \hbox{in} \quad L^2 (K \times \R^N ).
$$
Furthermore, using the estimate in (\ref{boh1}) we get that
$$
\sup_{y \in K , x \in \R^N} \left| (1+ |x|)^{N-4}  \left[ U_\e
({y\over \sqrt{\e}} , x ) - \mu_0^{-{N-2 \over 2}} (y) w_0 ({  x -\sqrt{\e}
\Phi_1 (y) \over\mu_0 (y) } ) \right]
\right| \to 0,
$$
as $\e \to 0$, where $\mu_0$ and $\Phi_1 $ are the smooth explicit
function defined in (\ref{choiceofmu0}) and (\ref{phi1def}).

Passing to the limit as $n \to \infty$ in (\ref{shen}) and applying
dominated convergence Theorem, we get
\begin{equation}\label{pas3}
\int_{K \times \R^N} \left[ |\nabla_x \varphi^\bot
|^2 - p \left( ( \mu_0)^{-{N-2 \over 2}} (y) w_0 ({ x - \sqrt{\e}\Phi_1 (y)
\over  \mu_0 (y) } \right)^{p-1} (\varphi^\bot)^2 \right] dy dx \leq
 0.
\end{equation}
Furthermore, passing to the limit in the orthogonality condition we
get, for any $y \in K$
\begin{equation}
\label{shen1}
\int_{\R^N} \varphi^\bot (y, x) Z_0 ({x -
\sqrt{\e}\Phi_1 (y) \over  \mu_0 (y) } ) dx = 0 ,
\end{equation}
\begin{equation}
\label{shen2}
\int_{\R^N} \varphi^\bot (y, x) Z_j ({x -
\sqrt{\e}\Phi_1 (y) \over  \mu_0 (y) }  ) dx = 0 , \quad
j=1, \ldots N
\end{equation}
and
\begin{equation}
\label{shen3} \int_{\R^N} \varphi^\bot (y, x) Z ({x - \sqrt{\e}\Phi_1
(y) \over  \mu_0 (y) } ) dx = 0 .
\end{equation}
We thus get a contradiction with (\ref{pas3}), since for any function
$\varphi^\bot$ satisfying the orthogonality conditions
(\ref{shen1})--(\ref{shen3}) for any $y \in K$ one has
$$
\int_{K \times \R^N } \left[ |\nabla_x \varphi^\bot |^2 - p \left(
\mu_0^{-{N-2 \over 2}} (y) w_0 ({x - \sqrt{\e}\Phi_1 (y) \over  \mu_0 (y)
} ) \right)^{p -1} (\varphi^\bot)^2 \right] dy
dx > 0.
$$

\bigskip
\noindent \textrm{\textbf{Step 2}}. For all $\e>0$ small, the
functional $P_\e (\delta )$ defined in (\ref{q0e}) is continuous and
differentiable in
 ${\mathcal H}^1 (K)$. Furthermore,
 $P_\e$ is a small perturbation in
$({\mathcal H}^1 (K) )^{N }$ of
\begin{eqnarray*}
 P (\delta )  &=&{A \over 2}  \e \left[ \int_K
|\nabla_K \delta |^2 + a_N
\int_K  {\mathcal H} \delta^2 \pm \e
b_{N}\int_K  \frac{\delta^2}{\mu_0} \right]
\end{eqnarray*}
where $A= \int_{\mathbb{R}^N} Z_0^2
$  and  ${\mathcal H} (y)$ defined in (\ref{defgy}).
Since we are assuming that  $\mu_0$ is a  nondegenerate solution to the following problem,
\begin{eqnarray*}
 - \Delta_K\mu +a_N {\mathcal H} (y)\mu \pm \frac{b_N}{\mu}=0,\quad \   \  in \ K,
\end{eqnarray*}
the operator $P$ is invertible.
Thus, for each $f \in L^2 (\Omega_{\e
, \gamma} )$,
$$
\delta  \in {\mathcal H}^1 (K) \longrightarrow \R , \quad \delta
\longmapsto P_\e (\delta ) - {\mathcal L}_f^3 (\delta )
$$
has a unique critical point $\delta $, which satisfies
$$
\ve^{-{k\over 2}} \| \delta \|_{{\mathcal H}^1 (K)} \leq
\widetilde{\sigma} \ve^{-2} \| f \|_{L^2 (\Omega_{\e , \gamma} )}
$$
for some proper $\widetilde{\sigma} >0$.

\bigskip

\noindent \textrm{\textbf{Step 3}}. For all $\e >0$ small, the
functional $Q_\e$ defined in (\ref{qje}) is a small perturbation in
$({\mathcal H}^1 (K) )^{N }$ of the quadratic form $\e Q_0 (d)$,
defined by
$$
\e  Q_0 (d)=  {\e  \over 2} C \left[ \int_K |\nabla_K d |^2 +
\int_K (\tilde g^{ab}\,R_{mabl} - \Gamma_a^c (E_m) \Gamma_c^a (E_l) ) d^m d^l
\right] $$ with $C:= \int_{\R^N} Z_1^2$ and the terms $R_{maal} $
and $\Gamma_a^c (E_m)$ are smooth functions on $K.$ Recall that the
non-degeneracy assumption on the minimal submanifold $K$ is
equivalent to the invertibility of the operator $Q_0 (d)$.

A consequence, for each $f \in L^2 (\Omega_{\e
, \gamma} )$,
$$
d \in ({\mathcal H}^1 (K))^{N} \longrightarrow \R , \quad d
\longmapsto Q_\e (d) - {\mathcal L}_f^3 (d)
$$
has a unique critical point $d$, which satisfies
$$
\ve^{-{k\over 2}} \| d \|_{({\mathcal H}^1 (K))^{N}} \leq
\widetilde{\sigma} \ve^{-2} \| f \|_{L^2 (\Omega_{\e , \gamma} )}
$$
for some proper $\widetilde{\sigma} >0$.

\bigskip

\noindent \textrm{\textbf{Step 4}}. Let $f \in L^2 (\Omega_{\e ,
\gamma})$ and assume that $e$ is a given (fixed) function in
${\mathcal H}^1 (K)$. We claim that for all $\e >0$ small enough,
the functional ${\mathcal Q} : (H_\e^1)^\bot \times ({\mathcal H}^1
(K) )^N \to \R$
$$
(\phi^\bot , \delta , d) \to {\mathcal E} (\phi^\bot , \delta , d ,
e)
$$
has a critical point $(\phi^\bot , \delta , d) $. Furthermore there
exists a positive constant $C$, independent of $\e$, such that
\begin{equation} \label{punto4}
\| \phi^\bot \| + \e^{-{k\over 2}} \bigg[ \| \delta \|_{{\mathcal
H}^1 (K)} + \| d \|_{({\mathcal H}^1 (K))^{N}} \bigg] \leq C
\e^{-2} \bigg[ \|f \|_{L^2 (\Omega_{\e , \gamma} )} + \e^{-{k\over
2}} \ve^{2} \| e \|_{{\mathcal H^1 (K) } }\bigg].
\end{equation}

To prove the above assertion, we first consider the functional
$$
{\mathcal Q}_0 (\phi^\bot , \delta , d) = {\mathcal Q} (\phi^\bot , \delta , d , e) - {\mathcal G} (\phi^\bot , \delta , d , e)
$$
where ${\mathcal G} $ is the functional that recollects all mixed
terms, as defined in (\ref{EEE}). A direct consequence of Step 1, Step
2 and Step 3 is that ${\mathcal Q}_0$ has a critical point
$(\phi^\bot = \phi^\bot (f) , \delta = \delta (f) , d = d(f))$,
namely the system
$$
D_{\phi^\bot } E (\phi^\bot ) = D_{\phi^\bot} {\mathcal L}_f^1
(\phi^\bot ), \quad \e^{-{k\over 2}} D_\delta P_\e (\delta ) =
D_\delta {\mathcal L}_f^2 (\delta ), \quad \e^{-{k\over 2}} D_d Q_\e
(d) = D_d {\mathcal L}_f^3 (d)
$$
is uniquely solvable in $(H_\e^1 )^\bot \times ({\mathcal H}^1 (K)
)^N $ and furthermore
$$
\| \phi^\bot \|_{H^1_\e} +\e^{-{k\over 2}} \| \delta \|_{{\mathcal
H}^1 (K)} + \e^{-{k\over 2}} \| d \|_{({\mathcal H}^1 (K) )^{N-1}}
\leq C \e^{-2} \| f \|_{L^2 (\mathcal{M}_{\e , \gamma })}
$$
for some constant $C>0$, independent of $\e$.

If we now consider the complete functional ${\mathcal Q}$, a
critical point of ${\mathcal Qx}$ shall satisfy the system
\begin{equation}\label{syst2}
  \begin{cases}
  D_{\phi^\bot } E(\phi^\bot ) = D_{\phi^\bot} {\mathcal L}_f^1 (\phi^\bot ) + D_{\phi^\bot } {\mathcal G} (\phi^\bot , \delta , d , e) \\
 D_\delta P_\e (\delta ) = D_\delta {\mathcal L}_f^2 (\delta ) + D_{\delta } {\mathcal G} (\phi^\bot , \delta , d , e)
    \\   D_d Q_\e (d ) = D_d {\mathcal L}_f^3 (d ) + D_{d } {\mathcal G} (\phi^\bot , \delta , d ,
    e).
  \end{cases}
\end{equation}
On the other hand, as we have already observed in Theorem
\ref{teo4.1}, we have
$$
\| D_{(\phi^\bot , \delta , d)} {\mathcal G}(\phi_1^\bot , \delta_1
, d_1 , e_1) - D_{(\phi^\bot , \delta , d)} {\mathcal G}(\phi_2^\bot ,
\delta_2 , d_2 , e_2) \| \leq C \e^{2} \times
$$
$$
\left[ \| \phi_1^\bot - \phi_2^\bot \| + \e^{-{k\over 2}} \|
\delta_1 - \delta_2 \|_{{\mathcal H}^1 (K)} + \e^{-{k\over 2}} \|
d_1 - d_2 \|_{({\mathcal H}^1 (K))^{N}} + \e^{-{k\over 2}} \| e_1
- e_2 \|_{{\mathcal H}^1 (K)}\right].
$$
Thus the contraction mapping Theorem guarantees the existence of a
unique solution $(\bar \phi^\bot , \bar \delta , \bar d)$ to
(\ref{syst2}) in the set
$$
\| \phi^\bot \|_{H_\e^1} + \e^{-{k\over 2}} \| \delta \|_{{\mathcal
H}^1 (K)} + \e^{-{k\over 2}} \| d \|_{({\mathcal H}^1 (K))^{N}}
\leq C \left[ \e^{-2} \| f \|_{L^2 (\Omega_{\e , \gamma} )} + \e^{2}
\e^{-{k\over 2}} \| e \|_{{\mathcal H}^1 (K)} \right].
$$
Furthermore, the solution $\bar \phi^\bot = \bar \phi^\bot (f,e) $,
$\bar \delta =\bar  \delta (f,e )$ and $\bar d= \bar d(f,e)$ depend
on $e$ in a smooth and non-local way.

\bigskip

\noindent {\rm \textbf{Step 5}}. Given $f\in L^2 (\Omega_{\e ,
\gamma} )$, we replace the critical point $(\bar \phi^\bot = \bar
\phi^\bot (f,e), \bar \delta =\bar  \delta (f,e ), \bar d= \bar
d(f,e) )$ of ${\mathcal Q}$ obtained in the previous step into the
functional ${\mathcal E}(\phi^\bot , \delta , d , e)$ thus getting a
new functional depending only on $e \in {\mathcal H}^1 (K)$, that we
denote by ${\mathcal F}_\e(e)$, given by
\begin{eqnarray*}
{\mathcal F}_\e(e) &=& \e^{-{k\over 2}} [ R_\e (e) - {\mathcal L}_f^4 (e) ] +
E(\bar \phi^\bot (e) ) -\ve^{-{k\over 2}} {\mathcal L}_f^1 (\bar \phi^\bot
(e)) +\e^{-{k\over 2}} [ P_\e (\bar \delta (e)) - {\mathcal L}_f^2 (\bar
\delta (e) )]
\\
&+& \e^{-{k\over 2}} [ Q_\e (\bar d (e)) - {\mathcal L}_f^3 (\bar d (e)) ] +
{\mathcal G} (\bar \phi^\bot (e) , \bar \delta (e) , \bar d (e) ,
e).
\end{eqnarray*}
The rest of the proof is devoted to show that there exists a
sequence $\e = \e_l \to 0$ such that
\begin{equation}
\label{effe1} D_e {\mathcal F}_\e(e) = 0
\end{equation}
is solvable. Using the fact that $(\bar \phi^\bot , \bar \delta ,
\bar d )$ is a critical point for ${\mathcal Q}$ (see Step 4 for the
definition), we have that
\begin{equation}
\label{effe11} D_e {\mathcal F}_\e (e) = \e^{-{k\over 2}} D_e [ R_\e (e) -
{\mathcal L}_f^4 (e) ]  + D_e {\mathcal G} (\bar \phi^\bot (e) ,
\bar \delta (e) , \bar d (e) , e).
\end{equation}
Define
\begin{equation}
\label{defLe} {\mathcal L}_\e := \e^{-{k\over 2}} D_e  R_\e (e)  + D_e
{\mathcal G} (\bar \phi^\bot (e) , \bar \delta (e) , \bar d (e) ,
e),
\end{equation}
regarded as self adjoint in ${\mathcal L}^2 (K)$. The work to solve
the equation $D_e {\mathcal F}_\e (e)=0$ consists in showing the
existence of a sequence $\e_l \to 0$ such that $0$ lies suitably far
away from the spectrum of ${\mathcal L}_{\e_l}$.

We recall now that the map
$$
(\phi^\bot , \delta , d , e ) \to D_e {\mathcal G} (\phi^\bot ,
\delta , d , e)
$$
is a linear operator in the variables $\phi^\bot , \delta , d$, while it is constant in $e$. This is contained in the result of Theorem \ref{teo4.1}. If we
furthermore take into account that the terms $\bar \phi^\bot $, $\bar \delta $ and $\bar d $ depend smoothly and in a non-local way through $e$, we conclude that, for
any $e \in {\mathcal H}^1 (K)$,
\begin{equation}
\label{effe2} D_e {\mathcal G} (\bar \phi^\bot (e) , \bar \delta (e)
, \bar d (e) , e )[e] = \e^{\gamma (N-3)} \e^{-{k\over 2}} \int_K \left( \e
\eta_1 (e) \partial_a e + \eta_2 (e) e \right)^2
\end{equation}
where $\eta_1$ and $\eta_2$ are  non local operators in $e$, that
are bounded, as $\e \to 0$, on bounded sets of ${\mathcal L}^2 (K)$.
Thanks to the result contained in Theorem \ref{teo4.1} and the above
observation, we conclude that the quadratic from
$$
\Upsilon_\e (e) := \e^{-{k\over 2}}  D_e R_\e (e) [e] + D_e {\mathcal G}
(\bar \phi^\bot (e) , \bar \delta (e) , \bar d (e) , e) [e]
$$
can be described as follows
\begin{equation}
\label{Uptilde} \tilde \Upsilon_\e (e) = \e^{k\over 2} \Upsilon_\e (e ) =
\Upsilon^0_\e (e)
 - {\bar \lambda_0 } \int_K e^2 + \e \Upsilon_\e^1 (e)
\end{equation}
where
\begin{equation}
\label{Up0} \Upsilon_\e^0 (e) = \e^2 \int_K  (1+ \e^{\gamma (N-3)}
\eta_1 (e) ) \left| \partial_a \left( e (1+ e^{-\e^{-\lambda'}}
\beta_3^\e (y)  ) \right) \right|^2.
\end{equation}
In the above expression $\bar \lambda_0 $ is the positive number defined
by
$$
\bar \lambda_0 = (\int_{\R^N} Z_1^2 )\, \lambda_0,
$$
$\Upsilon_e^1 (e)$ is a compact quadratic form in ${\mathcal H}^1
(K)$, $\beta_3^\e$ is a smooth and bounded (as $\e \to 0$) function
on $K$, given by (\ref{Q}). Finally, $\eta_1 $ is a non local operator
in $e$, which is uniformly bounded, as $\e \to 0$
 on bounded sets of ${\mathcal L}^2 (K)$.

Thus, for any $\e >0$, the eigenvalues of
$$
{\mathcal L}_\e e = \lambda e , \quad e \in {\mathcal H}^1 (K)
$$
are given by a sequence $\lambda_j (\e)$, characterized by the
Courant-Fisher formulas
\begin{equation}\label{courantFisher}
\lambda_j (\e )= \sup_{dim (M)= j-1} \inf_{e \in M^\bot \setminus \{0
\}} {\tilde \Upsilon_\e (e ) \over \int_K e^2 } = \inf_{dim (M)= j }
\sup_{e \in M\setminus \{0 \}} {\tilde \Upsilon_\e (e) \over \int_K
e^2 }.
\end{equation}
The proof of Theorem \ref{teouffa} and of the inequality (\ref{uffa1})
will follow then from Step 4 and formula (\ref{punto4}), together with
the validity of the following

\bigskip

\begin{lemma} \label{gafa}
There exist a sequence $\e_l \to 0 $ and a constant $c>0$ such that,
for all $j$, we have
\begin{equation}\label{autovalore}
|\lambda_j (\e_l  ) | \geq c \e_l^k.
\end{equation}
\end{lemma}

For the proof of Lemma we refer to \cite{demamu}.

\end{proof}

\section{Proof of Proposition \ref{linear}}\label{luigi}
\smallskip

 The proof of this Proposition will be divided into several steps.

{\bf Step 1}. \ \ Let us assume that $\phi$ solves (\ref{eq:eqwd}).
We  claim that there exists $C>0$ such that
\begin{equation}
\label{mar2} \|\phi \|_{\ve, r- 2}\le C \|h\|_{\ve,r}.
\end{equation}
By contradiction, assume that there exist sequences $\ve_n
\to 0$, $h_n$ with $\| h_n \|_{\ve_n , r} \to 0$ and solutions
$\phi_n$  to (\ref{eq:eqwd}) with $\| \phi_n \|_{\ve_n ,  r -2} =1$.

Let $z_n \in K_{\e_n}$ and $\xi_n$ be such that
$$
|\phi_n (\e_n z_n , \xi_n )| = \sup  |\phi_n (y, \xi )|.
$$
We may assume that, up to subsequences, $(\e_n z_n ) \to \bar y $ in
$K$. Furthermore, we have by assumption that
 $|\xi_n | \leq \eta \ve_n^{-{1\over 2}}$.

Let us now assume that there exists a positive constant $R$ such
that $|\xi_n |\leq R$. In this case, up to subsequences, one gets
that $\xi_n \to \xi_0$. Consider the functions
$$\tilde \phi_n ( z, \xi ) = \phi_n ( z , \xi + \xi_n ), \quad {\mbox {for}} \quad (z, \xi) \in K_{\e_n} \times
\{ \xi \in \R^N \, : \, |\xi | \leq \eta' \e_n^{-{1\over 2}} \} $$
for some $\eta' >0$.  This is a
sequence of uniformly bounded functions, that converges uniformly over compact sets of $K \times
\R^N$ to a function $\tilde \phi$ solution to
$$
    - \Delta \tilde \phi - p w_0^{p-1} \tilde \phi =0 \quad   \hbox{ in } \R^{N}
  $$
Since the orthogonality conditions pass to the limit, we get that
furthermore
$$
    \int_{\R^{N} } \tilde \phi (y, \xi ) Z_j (\xi ) \, d\xi = 0 \quad {\mbox {for all}} \quad y \in K, \quad  j=0, \ldots N.
    $$
These facts imply that $\tilde \phi \equiv 0$, that is a
contradiction.

Assume now that $\lim\limits_{n \to \infty} |\xi_n | = \infty$. Consider the scaled function
$$
\tilde \phi_n (z, \xi ) = \phi_n (z, |\xi_n | \xi + \xi_n )
$$
defined on the set
$$
\tilde \D=\left\{(z,\xi) :\  z \in K_{\e_n}, \, |\xi|<\frac{\eta}{\sqrt{\ve_n}|\xi_n|}-\frac{\xi_n}{|\xi_n|},\
\right\}.
$$
Thus
$\tilde \phi_n$ satisfies the equation
$$
\Delta \tilde \phi_n + p \,C_N {|\xi_n |^2 \over (1+| \, |\xi_n | \xi
+ \xi_n |^2 )^2 } \tilde \phi_n -|\xi_n |^2 \ve_n  a \tilde \phi_n =
|\xi_n |^2 h (z, |\xi_n | \xi +\xi_n )\ \ \mbox{in}\ \tilde \D.
$$
Consider first the case in which  $\lim\limits_{n \to \infty} \ve_n |\xi_n
|^2 = 0$. Under our assumptions, we have that $\tilde \phi_n$ is
uniformly bounded and it converges locally over compact sets to
$\tilde \phi$ solution to
$$
\Delta \tilde \phi = 0, \quad  |\tilde
\phi | \leq C |\xi |^{2-r}\quad {\mbox {in}} \quad \R^{N}\setminus \{ 0 \}.\quad
$$
Since $4<r<N$, we conclude that $\tilde \phi \equiv 0 $, which is a
contradiction.

Consider now the other possible case, namely that  $$\lim_{n \to
\infty} \ve_n |\xi_n |^2 = \beta >0.$$ Then,
$$
\tilde \D\to \mathcal{S}:=\{\xi \, : \, |\xi| \in [0,L)\} \quad \mbox{as}\ n\to\infty
$$
where $L$ is some positive constant.
Furthermore, up to subsequences, we
get that $\tilde \phi_n$ converges uniformly over compact sets to
$\tilde \phi$ solution to
$$
\Delta \tilde \phi - \beta a \tilde \phi = 0, \quad |\tilde \phi | \leq C |\xi |^{2-r} \quad  {\mbox {in}}
\  \mathcal{S}, \quad \tilde \phi = 0 \quad {\mbox {on}} \quad \partial \mathcal{S}.
$$
Multiplying equation by $\tilde\phi$, and integrating it over $ \mathcal{S}$ only in $\xi$, we get
$$
\int_{\mathcal{S}}(|\nabla \tilde\phi|^2+\beta a \tilde\phi^2)d\xi=0.
$$
Thus we conclude that $\tilde \phi \equiv 0$, which is a
contradiction. The proof
of (\ref{mar2}) is completed.

\smallskip
{\bf Step 2}. \ \  We shall now show that there exists $C>0$ such
that, if $\phi$ is a solution to (\ref{eq:eqwd}), then
\begin{equation}
\label{mar3} \| D^2_\xi \phi \|_{\ve , r } + \| D_\xi \phi \|_{\ve ,
r -1  } +\|\phi \|_{\ve, r- 2 }\le C \|h\|_{\ve,r }.
\end{equation}

For $z \in K_\ve$, we have that $\phi$ solves $ -\Delta \phi = \tilde h $
in $ |\xi |< \eta \ve^{-{1\over 2}} $ where $|\tilde h |\leq {C
\over (1+|\xi|^{r } )}$, for some constant $C>0$. Elliptic estimates give that $| \phi | \leq
{C \over (1+|\xi|^{r -2 } )}$.

Let us now fix a point $e \in \R^N$ and a positive number $R>0$.
Perform the change of variables $ \tilde \phi (z,t) = \phi (z, Rt
+3Re)$, so that
$$
\Delta \tilde \phi = {1\over R^{r-2}} \tilde h \quad {\mbox {in}}
\quad |t|\leq 1
$$
where $|\tilde h | \leq {c \over |t+3e|^{r}}$. Elliptic estimates
give then that $\| R^{r-2} D^2 \tilde \phi \|_{L^\infty(0,1)} \leq C
\| \tilde h \|_{L^\infty (B(0,2))}$, inequality that translates into
$$
\| R^{r} D^2 \phi \|_{L^\infty (B(3Re,R))} \leq C \| (1+|\xi |)^{r}
h \|_{L^\infty (|\xi |\leq \eta \e^{-{1\over 2}} )}.
$$
This inequality finally gives
$$
\| (1+|\xi |)^{r} D^2 \phi \|_{L^\infty (|\xi |\leq \eta
\e^{-{1 \over 2}} )} \leq C \| (1+|\xi |)^{r} h \|_{L^\infty (|\xi
|\leq \eta \e^{-{1\over 2}} )}.
$$
Arguing in a similar way, one gets the internal weighted estimate
for the first derivative of $\phi$
$$
\| (1+|\xi |)^{r-1} D \phi \|_{L^\infty (|\xi |\leq \eta
\e^{-{1\over 2}} )} \leq C \| (1+|\xi |)^{r} h \|_{L^\infty (|\xi
|\leq \eta \e^{-{1\over 2}} )}.
$$
By using the
representation formula for solution $\phi$ to the above equation, we
see that $ | \phi | \leq C \ve^{{r - 2 \over 2}}$ in $|\xi |<\eta
\e^{-{1\over 2}}$. Furthermore, elliptic estimates give that in this
region $ | D\phi | \leq C \ve^{{r - 1 \over 2}}$ and $ | D^2 \phi |
\leq C \ve^{{r  \over 2}}$. This concludes the proof of (\ref{mar3}).

\smallskip

{\bf Step 3}. \ \  We shall now show that there exists $C>0$ such
that, if $\phi$ is a solution to (\ref{eq:eqwd}), then
\begin{equation}
\label{mar4} \| D^2_\xi \phi \|_{\ve , r , \sigma } + \| D_\xi \phi
\|_{\ve , r -1 , \sigma } +\|\phi \|_{\ve, r- 2 , \sigma}\le C
\|h\|_{\ve,r ,\sigma} .
\end{equation}
From
elliptic regularity, we have that if  $ \|h\|_{\ve,r ,\sigma} \leq
C$ then $\|\phi \|_{\ve, r- 2 , \sigma} \leq C$. Thus, we write that
$\phi$ solves $ -\Delta \phi = \tilde h $ in $ |\xi |< \eta
\ve^{-{1\over 2}} $ where $\|\tilde h\|_{\ve,r ,\sigma} \leq C$.

Arguing as in the previous step, we fix a point $e \in \R^N$ and a
positive number $R>0$. Perform the change of variables $ \tilde \phi
(z,t) = \phi (z, Rt +3Re)$, so that
$$
\Delta \tilde \phi = {1\over R^{r-2}} \tilde h \quad {\mbox {in}}
\quad |t|\leq 1
$$
where $|\tilde h | \leq {c \over |t+3e|^{r}}$. Elliptic estimates
give then that $\| R^{r-2} D^2 \tilde \phi \|_{C^{0,\sigma}
(B(0,1))} \leq C \| \tilde h \|_{L^\infty (B(0,2))}$. This implies
that
$$
R^{r-2} \| D_\xi^2 \tilde \phi \|_{L^\infty (B_1)} + R^{r-2} [D^2
\tilde \phi ]_{\sigma, B(0,1)} \leq C.
$$
In particular, we have for any $z \in K_\ve$, that
$$
R^{r-2} \sup_{y_1 , y_2 \in B(0,1)} {|D^2 \tilde \phi (z, y_1) - D^2
\tilde \phi (z, y_2) | \over |y_1 - y_2 |^\sigma} \leq C.
$$
This inequality gets translated in term of $\phi$ as
$$
R^{r+\sigma} \sup_{\xi_1 , \xi_2 \in B(\xi,1)} {|D^2  \phi (z,
\xi_1) - D^2  \phi (z, \xi_2) | \over |\xi_1 - \xi_2 |^\sigma} \leq
C.
$$
In a very similar way, one gets the estimate on $D\phi$.
This concludes the
proof of (\ref{mar4}).

\smallskip

{\bf Step 4}. \ \ Differentiating equation (\ref{eq:eqwd}) with respect to the
$z$ variable $l$ times and using elliptic regularity estimates, one
proves that
\begin{equation}
\label{est1} \| D^l_y \phi \|_{\ve , r- 2 , \sigma }  \le C_l \left(
\sum_{k\leq l} \|D^k_y h\|_{\ve,r , \sigma}\right)
\end{equation}
for any given integer $l$.

\smallskip

\smallskip

{\bf Step 5}. \ \ Now we shall prove the existence of the solution
$\phi$ to problem (\ref{eq:eqwd}). We consider the Hilbert space $\mathcal{H}$ defined as the subspace of functions $\psi$ which are in $H^1(\D)$ such that $\psi = 0 $ on $\partial \hat \D $, and
$$
\int_{\hat\D} \psi (\e z, \xi ) Z_j (\xi ) \, d\xi = 0  \  {\mbox {for all}} \quad  z \in K_\e, \quad  j=0, \ldots N.
$$
Define a bilinear form in $\mathcal{H}$ by
$$
B(\phi,\psi):=\int_{\hat\D}\psi L\phi.
$$
Then problem (\ref{eq:eqwd}) gets weakly formulated as that of finding $\phi\in \mathcal{H}$ such that
$$
B(\phi,\psi)=\int_{\hat\D}h\psi\quad \forall\ \psi\in \mathcal{H}.
$$
By the Riesz representation theorem, this is equivalent to solve
\begin{eqnarray*}
\phi= T(\phi)+\tilde{h}
\end{eqnarray*}
with $\tilde{h}\in \mathcal{H}$ depending linearly on $h$, and $T:  \mathcal{H} \rightarrow \mathcal{H}$ being a compact operator.
Fredholm's alternative guarantees that there is a unique solution to problem (\ref{eq:eqwd}) for any $h$ provided that
\begin{eqnarray}\label{linear7}
\phi= T(\phi)
\end{eqnarray}
has only the zero solution in $\mathcal{H}$. Equation (\ref{linear7}) is equivalent to problem (\ref{eq:eqwd}) with $h=0$.
If $h=0$, the estimate in (\ref{est0a}) implies that $\phi=0$.

This concludes the proof of Proposition \ref{linear}.

 \end{document}